\documentclass[11pt]{article}
\usepackage{listings}
\usepackage{pdflscape}
\usepackage{amsfonts}
\usepackage{epsfig}
\usepackage{lscape}
\usepackage{amsmath,amssymb,amsthm,mathtools}
\usepackage{graphicx}
\usepackage{color}
\usepackage{placeins}
\usepackage{url}
\usepackage{cases,empheq}
\usepackage{longtable}
\usepackage{supertabular}

\allowdisplaybreaks[4]
\usepackage{algorithm,algpseudocode,setspace,float}
\usepackage{booktabs,multirow,array,threeparttable}
\usepackage{enumitem}
\usepackage{lipsum}
\usepackage{hyperref}
\usepackage{cleveref}
\usepackage{svg}

\usepackage[numbers,sort]{natbib}
\usepackage{makecell}
\usepackage{rotating,tabularx}
\usepackage{array,bm}
\usepackage{float,framed,subfigure,multicol}  
\usepackage{adjustbox}
\usepackage{nicematrix}
\usepackage{changepage} 
\usepackage{tikz}
\usepackage{fullpage}



\theoremstyle{plain}
\newtheorem{theorem}{Theorem}[section]
\newtheorem{claim}[theorem]{Claim}
\newtheorem{lema}[theorem]{Lemma}
\newtheorem{coro}[theorem]{Corollary}
\newtheorem{assumption}{Assumption}
\theoremstyle{remark}
\newtheorem{remark}[theorem]{Remark}
\theoremstyle{definition}

\numberwithin{equation}{section}

\newcommand{\mb}[1]{\mathbb{#1}}
\newcommand{\mf}[1]{\mathbf{#1}}
\newcommand{\mr}[1]{\mathrm{#1}}
\newcommand{\mc}[1]{\mathcal{#1}}

\newcommand{\pminf}{(-\infty,+\infty]}

\newcommand{\iidsm}{\overset{\mathrm{i.i.d.}}{\sim}}
\newcommand{\proj}{\mr{proj}}
\newcommand{\prox}{\mr{prox}}
\newcommand{\dist}{\mr{dist}}
\newcommand{\dom}{\mr{dom}}

\newcommand{\olx}{\overline{x}}

\newcommand{\olf}{\overline{f}}

\newcommand{\olL}{\overline{L}}

\newcommand{\olrho}{\overline{\rho}}
\newcommand{\oltheta}{\overline{\theta}}

\newcommand{\oltau}{\overline{\tau}}
\newcommand{\oleta}{\overline{\eta}}
\newcommand{\olalpha}{\overline{\alpha}}
\newcommand{\alphaol}{\alpha_{l0}}
\newcommand{\alphaou}{\alpha_{u0}}
\newcommand{\olvphi}{\overline{\varphi}}
\newcommand{\tlx}{\tilde{x}}

\newcommand{\tlF}{\tilde{F}}
\newcommand{\tlm}{\tilde{m}}
\newcommand{\tly}{\tilde{y}}
\newcommand{\tlxi}{\tilde{\xi}}

\newcommand{\tlrho}{\tilde{\rho}}
\newcommand{\hrho}{\hat{\rho}}
\newcommand{\hatL}{\hat{L}}
\newcommand{\lipsmf}{L_{\mathrm{sm}}}
\newcommand{\mephi}[1]{e_{\left(\overline{\rho}\right)^{-1}\phi}^{#1}}
\newcommand{\Norm}[1]{\left\Vert #1\right\Vert} 
\newcommand{\norm}[1]{\left\vert #1\right\vert} 
\newcommand{\lrangle}[1]{\left\langle #1\right\rangle} 
\newcommand{\lrbrackets}[1]{\left\{ #1\right\}}
\newcommand{\lrrbrackets}[1]{\left( #1\right)}

\newcommand{\conv}{\mbox{conv}}

\newcommand{\glip}{\mathsf{G}_{\mathrm{Lip}}}
\newcommand{\losslin}{F_{\mathrm{lin}}}
\newcommand{\losslog}{F_{\mathrm{log}}}

\newcommand{\upperroman}[1]{\uppercase\expandafter{\romannumeral#1}}

\usepackage{caption}
\begin{document}

\title{\bf On the Convergence Analysis of an Inexact Preconditioned Stochastic Model-Based Algorithm \footnotemark[1]}
\author{Chenglong Bao\footnotemark[2], \quad Yancheng Yuan\footnotemark[3], \quad Shulan Zhu\footnotemark[4]}
\date{\today}
\maketitle

\renewcommand{\thefootnote}{\fnsymbol{footnote}}
\footnotetext[1]{{\bf Funding:} The work of Chenglong Bao was supported by the National Key R\&D Program of China under grant 2021YFA1001300, and the National Natural Science Foundation of China under grant 12271291. The work of Yancheng Yuan was supported by the RGC Early Career Scheme (Project No. 25305424).
}
\footnotetext[2]{Yau Mathematical Sciences Center, Tsinghua University, Beijing, China ({\tt clbao@tsinghua.edu.cn}).}
\footnotetext[3]{Department of Applied Mathematics, The Hong Kong Polytechnic University, Hung Hom, Hong Kong (\textbf{Corresponding author}. {\tt yancheng.yuan@polyu.edu.hk}).}
\footnotetext[4]{Department of Mathematical Sciences, Tsinghua University, Beijing, China ({\tt zhusl22@mails.tsinghua.edu.cn}).}
\renewcommand{\thefootnote}{\arabic{footnote}}

\begin{abstract}
     This paper focuses on investigating an inexact stochastic model-based optimization algorithm that integrates preconditioning techniques for solving stochastic composite optimization problems. The proposed framework unifies and extends the fixed-metric stochastic model-based algorithm to its preconditioned and inexact variants. Convergence guarantees are established under mild assumptions for both weakly convex and convex settings, without requiring smoothness or global Lipschitz continuity of the objective function. By assuming a local Lipschitz condition, we derive nonasymptotic and asymptotic convergence rates measured by the gradient of the Moreau envelope. Furthermore, convergence rates in terms of the distance to the optimal solution set are obtained under an additional quadratic growth condition on the objective function. Numerical experiment results demonstrate the theoretical findings for the proposed algorithm. 
\end{abstract}

\newcommand{\keywords}[1]{%
  \par\vspace{0.5\baselineskip}%
  \noindent{\textbf{Keywords:} #1}
}


\section{Introduction}\label{sec:introduction}
We consider the following stochastic composite optimization problem:
\begin{align}\label{equ:compositeoptimization}
    \min_{x\in \mb{R}^d}\ \phi(x)\triangleq F(x) + r(x),\quad \mbox{where}\quad F(x)\triangleq \mb{E}_{s\sim P}\left[f(x;s)\right],\tag{ComOpt}
\end{align}
where $\mc{S}$ is a sample space, $P$ represents a distribution over $\mc{S}$, $r\colon \mb{R}^d\to \pminf$ is a proper and closed function, and for $P$-almost $s\in \mc{S}$, the component function $f(\cdot;s)\colon \mb{R}^d\to \pminf$ is proper and closed, and the composite component function $\varphi(\cdot;s) \triangleq f(\cdot;s)+r(\cdot)$ is weakly convex. The weakly convex problems in the form of (\ref{equ:compositeoptimization}) cover a wide class of nonconvex optimization problems with numerous applications in machine learning and statistics, including robust phase retrieval, robust PCA, and blind deconvolution \cite{duchi2019solving,charisopoulos2021low,ling2015self}. In particular, it is well known that the composition of a Lipschitz continuous convex function with an $L$-smooth function is weakly convex \cite[Lemma~4.2]{drusvyatskiy2019efficiency}. Recently, it has attracted growing attention for researchers to design provably fast algorithms for solving (\ref{equ:compositeoptimization}), and significant progress has been achieved for designing stochastic model-based algorithms. However, the analysis of these algorithms still relies on some restrictive assumptions, which may not be satisfied in real applications. In this paper, we will provide a new analysis of the convergence and convergence rates of the inexact stochastic model-based algorithm for solving (\ref{equ:compositeoptimization}) under mild and more practical assumptions. Next, we move on to briefly discuss the motivation and some related works. 

\subsection{Motivation and Related Works}\label{subsec:motivationrelatedworks}
Among the numerous stochastic first-order methods for solving problem (\ref{equ:compositeoptimization}), the stochastic gradient descent (SGD) method and its variants remain popular in modern machine learning and deep learning \cite{robbins1951stochastic,duchi2011adaptive,kingma2015adam}. However, both theoretical and empirical
studies have demonstrated that SGD suffers from instability and faces significant challenges in stepsize selection \cite{moulines2011non,asi2019stochastic}, which will reduce its practical efficiency due to time-consuming stepsize tuning. The stochastic proximal point algorithm (sPPA) has gained growing attention recently in stochastic optimization for its superior stability compared to SGD \cite{bertsekas2011incremental,ryu2014stochastic,asi2019stochastic}. Convergence guarantees for sPPA and its variants have been developed in several settings. For the convex case, \citet{zhu2025tight} extended the stability results for sPPA in \cite{asi2019stochastic} to the inexact variants under mild conditions, and derived convergence rate bounds under a local Lipschitz condition on the component functions and a quadratic growth condition on the objective function. In the weakly convex setting, \citet{davis2019stochastic} established convergence guarantees under certain assumptions for the stochastic model-based algorithm based on the so-called model functions $f_x(\cdot;s)$ satisfying
\begin{align*}
    \mb{E}_{s\sim P}[f_x(x;s)] = f(x)\quad \mbox{and}\quad \mb{E}_{s\sim P}[f_x(y;s)-f(y)]\leq \frac{\eta}{2}\Norm{y-x}_2^2\quad \mbox{for all}\ x,y
\end{align*}
for some constant $\eta\in \mb{R}_{+}$. Building on the stochastic model-based framework, it can unify the convergence analysis of several widely used methods, including SGD, sPPA, and the stochastic prox-linear algorithm \cite{duchi2018stochastic}.

\paragraph*{Stability analysis in the weakly convex setting.} The convergence guarantees derived in \cite{davis2019stochastic} apply to weakly convex functions, but their analysis relies on a global Lipschitz condition on the model function, which can be restrictive and fail to hold for important applications. For instance, the model functions employed in \cite{davis2019stochastic} for solving the phase retrieval problem fail to satisfy this global Lipschitz condition. To extend the theoretical results to a broader class of problems, it is therefore necessary to relax the global Lipschitz requirement to a more practical local version. Under a local Lipschitz condition, a main difficulty for establishing the convergence guarantee lies in ensuring the stability of the stochastic algorithm \cite{gao2024stochastic}. Following \cite[Definition~3.1]{asi2019stochastic}, a stochastic algorithm is said to be stable if the iterates generated with square summable stepsizes remain bounded almost surely. Some sufficient conditions have been proposed to establish such stability for the stochastic model-based algorithm with square summable stepsizes. In the weakly convex setting, when the regularizer satisfies the $(\lambda,\sigma)$-regular coercivity 
\begin{equation}\label{equ:regcoercive}
    \begin{split}
        \lim_{\Vert x\Vert_2\to \infty}\frac{r(x)}{\Vert x\Vert^\sigma} &= \infty,\quad \mbox{and}\\ r(x)&\geq r(y)\quad \mbox{for}\ \Norm{x}_2,\Norm{y}_2\ \mbox{sufficiently large}\ \mbox{and}\ \Norm{y}_2\leq \lambda\Norm{x}_2,
    \end{split}
\end{equation}
 the following generalized Lipschitz condition is introduced in \cite{duchi2018stochastic} to guarantee the boundedness of the iterates generated by the stochastic proximal subgradient method:
\begin{equation}\label{equ:generalizedloclipnu}
    \norm{f(x;s)-f(y;s)}\leq L_f\left(1+\Norm{x}_2^\nu\right) \Norm{x-y}_2\quad \mbox{for all}\ x,y.
\end{equation}
However, their analysis requires $\nu < \sigma - 1$, which can be restrictive in applications such as phase retrieval, where $\nu = 1$. In this case, the regularizer must grow strictly faster than the squared $\ell_2$-norm to ensure boundedness. Beyond this condition, based on robust regularization stepsize strategies, and assuming that the Moreau envelope of the objective function has bounded sublevel sets, \citet{gao2024stochastic} established the stability of the stochastic model-based algorithm under a more flexible generalized Lipschitz condition on the model function $f_x(\cdot;s)$:
\begin{align*}
    \norm{f_x(y;s) - f_x(z;s)} \leq L_f(s)\glip\left(\Norm{x}_2\right)\Norm{y-z}_2\quad \mbox{for all}\ x,y,z,
\end{align*}
where $\glip\colon \mb{R}_{+}\to \mb{R}_{+}$ is a monotonically nondecreasing function. Nevertheless, their analysis still requires the regularizer $r$ to be globally Lipschitz continuous over its domain $\dom(r)$, thereby excluding regularizers involving the squared $\ell_2$-norm. 

For the convex case, \citet{nguyen2018sgd} proved the stability of SGD under the Lipschitz smoothness assumption on the component function $f(\cdot;s)$. Furthermore, \citet{zhu2025tight} demonstrated the stability of the inexact stochastic proximal point algorithm (isPPA) under the assumptions that (\romannumeral1) the subdifferential and expectation operators commute, and (\romannumeral2) the subgradients of composite component functions are uniformly bounded over the optimal solution set, i.e., 
\begin{align*}
    \mb{E}_{s\sim P}\left[\varphi^\prime\left(x^{*};s\right)\right] = \phi^\prime\left(x^{*}\right)\quad \mbox{and}\quad & \mb{E}_{s\sim P}\left[\Norm{\varphi^\prime\left(x^{*};s\right)}_2^2\right]\leq \sigma_\phi^2
\end{align*}
holds for all $x^{*}\in \arg\min_{x\in \mb{R}^d}\phi(x)$, where $\varphi^\prime\left(x^{*};s\right)\in \partial \varphi\left(x^{*};s\right)$ and $\phi^\prime\left(x^{*}\right)\in \partial \phi\left(x^{*}\right)$. Note that the stability of isPPA in the convex setting does not depend on any Lipschitz assumptions on $f(\cdot;s)$ or $r$. Extending such analysis to the stochastic model-based algorithm under weak convexity, however, remains challenging. 

Moreover, the analyses of \cite{duchi2018stochastic} and \cite{gao2024stochastic} for the weakly convex case fundamentally differ from the analysis in \cite{zhu2025tight} for the convex case. 

\paragraph*{Inexact subproblem solving and preconditioning techniques.} A significant advantage of sPPA and the stochastic prox-linear algorithm over the stochastic subgradient method lies in their robustness to the choice of stepsizes, as demonstrated empirically in \cite{davis2019stochastic}. However, the subproblem arising in these two algorithms may not admit closed-form solutions, particularly when the batch size is large. In such cases, solving each subproblem exactly becomes time-consuming (or impractical), making it essential to consider the stochastic model-based algorithm that allows for inexact subproblem solutions. To further improve the overall computational efficiency, preconditioning techniques have recently been explored to obtain approximate subproblem solutions efficiently. In particular, the deterministic PPA combined with preconditioning techniques has attracted growing attention due to its broad applications in splitting algorithms for convex minimization and monotone inclusion problems \cite{rockafellar2021advances,bredies2022degenerate,rockafellar2023generic,sun2024accelerating}. Compared with its determinstic counterparts, the application of preconditioning techniques to the sPPA remains underexplored. Recently, \citet{milzarek2024semismooth} considered the preconditioning technique and designed an algorithim, termed SNSPP, for solving finite-sum composite optimization problems of the form
\begin{equation}\label{equ:compositefinitesum}
    \begin{split}
        \min_{x\in \mb{R}^d} \phi(x)\triangleq \frac{1}{n} \sum_{i=1}^n f_i\left(A_i x\right) + r(x),
    \end{split}\tag{ComFSum}
\end{equation}
where $A_i\in \mb{R}^{n_i\times d}$ is given, each component functions $f_i\colon \mb{R}^{n_i}\to \mb{R}$ is continuously differentiable, and the regularizer $r\colon \mb{R}^d\to \mb{R}$ is proper, closed and convex. The framework of SNSPP depends on the following assumptions: (\romannumeral1) each function $f_i$ is $L_i$-Lipschitz smooth and $\gamma_i$-weakly convex, (\romannumeral2) the conjugate function $\hat{f}_i^\star$ is essentially differentiable with a semismooth gradient mapping $\nabla \hat{f}_i^\star$ on $\mbox{int}(\dom(\hat{f}_i^\star))$, where $\hat{f}_i(\cdot) \triangleq f_i(\cdot) + \frac{\gamma_i}{2}\Norm{\cdot}_2^2$, and (\romannumeral3) the proximal mapping $\prox_{\alpha r}(\cdot)$ is semismooth on $\mb{R}^d$ for all $\alpha\in \mb{R}_{++}$. 
Under these assumptions, \citet{milzarek2024semismooth} established convergence guarantees for both weakly convex and strongly convex cases. When viewed as a variable-metric sPPA, SNSPP employs preconditioners of the form 
\begin{align*}
    M_k \triangleq I_d + \frac{\alpha_k}{\vert \mc{S}_k \vert}\sum_{i\in \mc{S}_k} \gamma_i A_i^\top A_i,
\end{align*}
where $\alpha_k$ is the stepsize and $\mc{S}_k$ denotes the randomly sampled minibatch at iteration $k$. SNSPP was tailored specifically to finite-sum problems of the form (\ref{equ:compositefinitesum}). This fact, together with the requirement of Lipschitz smoothness for each $f_i$, limits its applicability to more general stochastic composite optimization settings. 

\subsection{Contributions}
In this paper, we consider the model-based \textbf{i}nexact \textbf{s}tochastic \textbf{p}reconditioned \textbf{p}roximal \textbf{p}oint \textbf{a}lgorithm (ispPPA) for solving the stochastic nonconvex composite optimization problem (\ref{equ:compositeoptimization})  and provide the convergence analysis without requiring global Lipschitz assumptions. The main contributions of this paper are summarized as follows:
\begin{itemize}
    \item We incorporate preconditioning techniques into the stochastic model-based optimization framework and propose practically implementable stopping criteria for solving each subproblem inexactly.
    \item We establish sufficient conditions for the stability of the model-based ispPPA in both weakly convex and convex settings, while allowing for nondifferentiability and local Lipschitz continuity of both the model functions and the regularizer.
    \item We derive both nonasymptotic and asymptotic convergene rates of the model-based ispPPA under a local Lipschitz condition on the model functions for the weakly convex case.
\end{itemize}

\subsection{Notations and Blanket Assumptions}\label{subsec:notations}
Denote by $[n] \triangleq \lrbrackets{1,\cdots,n}$. For any $x\in \mb{R}^d$, denote its $q$-norm ($q \geq 1$) as $\Norm{x}_q: = (|x_1|^q + \cdots + |x_d|^q)^{1/q}$. For any self-adjoint and positive-definite linear mapping $M\colon \mb{R}^d\to \mb{R}^d$, the associated inner products and norms are defined by 
\begin{align*}
    \lrangle{x,y}_{M}\triangleq \lrangle{x,My} = \lrangle{Mx,y}\ \mbox{and}\ \Norm{x}_M^2 \triangleq \lrangle{x,Mx}\quad \mbox{for all}\ x,y\in \mb{R}^d.
\end{align*}
For any $x\in \mb{R}^d$ and $\mc{X}\subset \mb{R}^d$, let $\dist_{M}(x,\mc{X}) \triangleq \inf_{y\in \mc{X}}\Norm{x-y}_M$ be the distance in the $M$-norm from $x$ to $\mc{X}$, and set $\proj_M(x,\mc{X})$ to be the projection in the $M$-norm of $x$ onto $\mc{X}$ if $\mc{X}$ is nonempty and closed convex. For the fixed-metric case, i.e., $M = I_d$, we omit the subscript and simply use $\dist(x,\mc{X})$ and $\proj(x,\mc{X})$ to denote the associated distance and projection. Let $p\colon \mb{R}^d \to \pminf$ be a proper closed function and let $x\in \dom(p)$. When $p$ is locally Lipschitz continuous on a neighborhood $U$ of $x$, the B-subdifferential of $p$ at $x$ is defined as 
\begin{align*}
    \partial_B p(x) \triangleq \lrbrackets{v\ \vert\ \exists\ \{x^v\}\in U_p\ \mbox{with}\ \lim_{\nu\to \infty} x^v = x\ \mbox{and}\ \lim_{\nu\to \infty}\nabla p\left(x^\nu\right)= v},
\end{align*}
where $U_p$ denotes the set of all differentiable points of $p$ in $U$. The Clarke subdifferential of $p$ at $x$ is defined as the convex hull of the B-subdifferential, i.e., $\partial_C p(x) = \conv(\partial_B p(x))$. The regular subdifferential of $p$ at $x$ is defined as 
\begin{align*}
    \partial p(x) \triangleq \lrbrackets{v\in \mb{R}^d\ \vert\ \lim\inf_{x\neq y\to x} \frac{p\left(y\right) - p\left(x\right) - \lrangle{v,y-x}}{\Norm{y-x}_2}\geq 0}.
\end{align*}
A function $p\colon \mb{R}^d \to \pminf$ is said to be proper, closed and $\oltau$-weakly convex if function $p(\cdot) + \frac{\oltau}{2}\Norm{\cdot}_2^2$ is proper, closed and convex. The following lemma is from \cite[Lemma~2.1]{davis2019stochastic}.
\begin{lema}\label{lema:weaklyconvex}
    Let $p\colon \mb{R}^d\to \pminf$ be a proper closed function. Then the following three assertions are equivalent.
    \begin{enumerate}
        \item $p$ is $\oltau$-weakly convex.
        \item It holds that for all $x,y\in \mb{R}^d$ and $\theta\in [0,1]$, 
        \begin{align*}
            p\left(\theta x + \left(1-\theta\right)y\right)\leq \theta p\left(x\right) + \left(1-\theta\right)p\left(y\right) + \frac{\oltau \theta\left(1-\theta\right)}{2}\Norm{x-y}_2^2.
        \end{align*}
        \item Denote by $\partial p(x)$ the regular subdifferential of $p$ at $x$. It holds that for all $x\in \mb{R}^d$ and $v\in \partial p(x)$, 
        \begin{align*}
            p\left(y\right)\geq p\left(x\right) + \lrangle{v, y-x} - \frac{\oltau}{2}\Norm{y-x}_2^2\quad \mbox{for all}\ y\in \mb{R}^d.
        \end{align*}
    \end{enumerate}
\end{lema}
For any proper, closed and $\oltau$-weakly convex function $p\colon \mb{R}^d \to \pminf$, $p$ is Clarke regular \cite[Proposition~4.4.15]{cui2021modern}, i.e., $\partial_C p(x) = \partial p(x)$ for any $x\in \mb{R}^d$. For any $\alpha>0$, the Moreau envelope $e_{\alpha p}^M(\cdot)$ and the proximal mapping $\prox_p^M(\cdot)$ associated with $p$ in the $M$-norm are defined by 
\begin{align*}
        e_{\alpha \phi}^M\left(x\right) &\triangleq \min_{y\in \mb{R}^d}\lrbrackets{p(y) + \frac{1}{2\alpha}\Norm{y-x}_M^2},\\
        \prox_{\alpha p}^M(x) &\triangleq \arg\min_{y\in \mb{R}^d}\lrbrackets{p(y) + \frac{1}{2\alpha}\Norm{y-x}_M^2}\quad \mbox{for all}\ x\in \mb{R}^d.
\end{align*}
Similarly, when $M = I_d$, we omit the superscript and simply use $e_{\alpha p}(\cdot)$ and $\prox_{\alpha p}(\cdot)$ to denote the associated Moreau enevlope and proximal mapping. 
Since all norms are equivalent in finite dimensional vector spaces, there exists constant $\mu_M\in (0,1]$ and $L_M\in [1,\infty)$ such that 
\begin{align*}
    \mu_M \Norm{x}_M \leq \Norm{x}_2\leq L_M\Norm{x}_M\quad \mbox{for all}\ x\in \mb{R}^d.
\end{align*}
Then, for any $\oltau$-weakly convex function $p$ and positive-definite linear mapping $M$, if $\alpha\in (0,\tau^{-1})$ with $\tau \triangleq L_M^2\cdot\oltau$, by following the argument in the proof of \cite[Proposition~3.3]{hoheisel2010proximal}, one can show that the proximal mapping $\prox_{\alpha p}^M(\cdot)$ is single-valued and $(1-\alpha \tau)^{-1}$-Lipschitz continuous in $M$-norm, and the Moreau envelope $e_{\alpha p}^M(\cdot)$ is differentiable in $\ell_2$-norm with gradient given by\footnote{See Appendix \ref{apxsubsubsec:moreauenvelopediff} for a detailed proof.} 
\begin{equation}\label{equ:gradmegenral}
    \nabla e_{\alpha p}^M\left(x\right) = \alpha^{-1}M\left(x - \prox_{\alpha p}^M\left(x\right)\right).
\end{equation}
Consider the stochastic composite optimization problem (\ref{equ:compositeoptimization}). Let $\mc{X}^{*}\triangleq \arg\min_{x\in \mb{R}^d}\phi(x)$ denote the optimal solution set, and set $\phi^{*} \triangleq \inf_{x\in \mb{R}^d} \phi(x)$. We adopt the following stochastic model function from \cite{davis2019stochastic}.
\begin{assumption}\label{asm:stomodel}
    There exists a scalar $\overline{\eta}\in \mb{R}$, an open convex set $V$ containing $\dom(r)$ and a measurable function $(x,y,s)\mapsto f_x(y,s)$, defined on $V\times V\times \mc{S}$, such that for $P$-almost $s\in \mc{S}$, 
    \begin{align*}
        f_x\left(x,s\right) = f(x,s)\quad \mbox{and} \quad \norm{f_x\left(y,s\right)-f\left(y,s\right)} \leq \frac{\oleta}{2}\Norm{y-x}_2^2
    \end{align*}
    holds for all $x,y\in V$.
\end{assumption}
For each $s\in \mc{S}$, the composite component model function $\varphi_x(\cdot;s)\colon \mb{R}^d \to \pminf$ is defined by 
\begin{align*}
    \varphi_x(y;s) \triangleq f_x(y;s) + r(y)\quad \mbox{for all}\ x,y\in V.
\end{align*}
For each $k\in \mb{Z}_{+}$, we use $S_k^{1:m}\triangleq \{S_k^i\}_{i=1}^m\subset \mc{S}$ to denote a random minibatch of size $m$ and define 
\begin{subequations}
    \begin{align}
        \overline{f}_{x_k}(x;S_k^{1:m}) &\triangleq \frac{1}{m}\sum_{i=1}^m f_{x_k}(x;S_k^i)\quad \mbox{for all}\ x\in V,\label{equ:minibatchf}\\
        \olvphi_{x_k}\left(x;S_k^{1:m}\right) &\triangleq \frac{1}{m}\sum_{i=1}^m \varphi_{x_k}(x;S_k^i) = \overline{f}_{x_k}\left(x;S_k^{1:m}\right) + r(x)\quad \mbox{for all}\ x\in V.\label{equ:minibatchphi}
    \end{align}
\end{subequations}
To solve problem (\ref{equ:compositeoptimization}), we consider the model-based ispPPA method shown in Algorithm \ref{algo:ispppa}, which reduces to the stochastic model-based algorithm \cite{davis2019stochastic} when $M_k \equiv I_d$ and $\epsilon_k \equiv 0$, and simplifies to the isPPA under the setting $M_k \equiv I_d$ and $f_x(\cdot;s) = f(\cdot;s)$. 
\begin{algorithm}[htbp]
    \caption{Model-Based Inexact Stochastic Preconditioned Proximal Point Algorithm (ispPPA)}
    \label{algo:ispppa}
    \begin{algorithmic}[1]
        \Require initial point $x_1$, stepsizes $\{\alpha_k\}_{k\in \mb{Z}_{+}}$, accuracy parameters $\{\epsilon_k\}_{k\in \mb{Z}_{+}}$, minibatch size $m$, maximum iteration count $K$.
        \For{$k = 1,2,\cdots, K$}
            \State Draw a random minibatch $S_k^{1:m}$ with $S_k^i \iidsm P$.
            \State Update 
            \begin{equation}\label{equ:sippm}
                x_{k+1}\overset{\epsilon_k}{\approx} \arg\min_{x\in \mb{R}^d}\left\{\overline{f}_{x_k}(x;S_k^{1:m}) + r(x) + \frac{1}{2\alpha_k}\Vert x - x_k\Vert_{M_k}^2\right\}
            \end{equation}
            where $\overline{f}_{x_k}(x;S_k^{1:m})$ is defined in (\ref{equ:minibatchf}).
        \EndFor
        \State\textbf{Output}: $x_K$; $x_{K_{*}}$ drawn from $\{x_k\}_{k\in [K]}$  with probability $\mb{P}\{K_{*} = k\} = \frac{\alpha_k}{\sum_{i=1}^K \alpha_i}$.
    \end{algorithmic}
\end{algorithm}
Unless otherwise stated, we assume that the linear mapping $M_k\colon \mb{R}^d\to \mb{R}^d$ is self-adjoint and positive-definite for each $k\in \mb{Z}_{+}$, and such mapping $M_k$ is always referred to as the preconditioner in this work. The symbol ``$\overset{\epsilon_k}{\approx}$'' indicates that the next iterate $x_{k+1}$ is obtained by approximately solving the following subproblem:
\begin{align*}
    \min_{x\in \mb{R}^d}\ \overline{f}_{x_k}(x;S_k^{1:m}) + r(x) + \frac{1}{2\alpha_k}\Vert x - x_k\Vert_{M_k}^2
\end{align*}
until $x_{k+1}$ satisfies an accuracy $\epsilon_k$ specified by some stopping criterion (to be defined in Section \ref{sec:convergenceispp}). To establish the convergence properties of the model-based ispPPA (Algorithm \ref{algo:ispppa}), we make the following assumptions. 
\begin{assumption}\label{asm:mk}
    Fix any $k\in \mb{Z}_{+}$. Let $\mc{F}_k\triangleq\sigma(S_1^{1:m},\cdots,S_k^{1:m})$ be the $\sigma$-field generated by the first $k$ random minibatch $\{S_i^{1:m}\}_{i=1}^k$, and denote the conditional expectation $\mb{E}_k[\cdot]\triangleq\mathbb{E}[\cdot\ |\ \mathcal{F}_{k-1}]$. There exist scalars $\mu_\infty\in (0,1]$ and $L_\infty\in [1,\infty)$ such that the self-adjoint positive-definite linear mapping $M_k\colon \mb{R}^d\to \mb{R}^d$ satisfies 
    \begin{align*}
            \mu_k &\Norm{x}_{M_k}\leq \Norm{x}_2 \leq L_k \Norm{x}_{M_k}\quad \mbox{for all}\ x\in \mb{R}^d,
    \end{align*}
    where $\mu_k \in [\mu_\infty,\infty)$ and $L_k \in (0,L_\infty]$. Moreover, for each $k\in \mb{Z}_{+}$, the $M_{k-1}$-norm and $M_k$-norm are related to each other by 
    \begin{align*}
        \rho_k^{-2}\leq \mb{E}_k \left[\cfrac{\Norm{x}_{M_k}^2}{\Norm{x}_{M_{k-1}}^2}\right]\leq \rho_k^2\quad \mbox{for all}\ 0\neq x\in \mb{R}^d, 
    \end{align*}
    where $\rho_k\in [1,\infty)$ and $\nu_k\triangleq \prod_{i=0}^{k-1}\rho_i \nearrow \nu_\infty\triangleq \prod_{k=0}^\infty \rho_k < \infty$. Here we define $M_0 \triangleq I_d$ along with $\mu_0 \triangleq 1$, $L_0 \triangleq 1$ and $\rho_0 \triangleq 1$.
\end{assumption}
\begin{assumption}\label{asm:sppaall}
    The function $F\colon \mb{R}^d\to \mb{R}$ is locally Lipschitz and $r\colon \mb{R}^d\to \pminf$ is a proper closed function. For $P$-almost all $s\in \mc{S}$, $f(\cdot;s)\colon \mb{R}^d\to \pminf$ is a proper closed function satisfying $\mr{dom}(f(\cdot;s))\supset \mr{dom}(r)$, and $\varphi_x(\cdot;s)$ is $\oltau$-weakly convex for all $x\in V$ with constant $\oltau\in \mb{R}_{+}$ and the open convex set $V$ defined in Assumption \ref{asm:stomodel}.
\end{assumption}
\begin{assumption}\label{asm:sppa1}
    For the open convex set $V$ specified in Assumption \ref{asm:stomodel} and any bounded subset $U\subset V$, there exists a measurable function $\overline{L}_{f,U}\colon \mc{S}\to \mb{R}_{+}$, satisfying $\sqrt{\mb{E}_{s\sim P}[\overline{L}_{f,U}(s)^2]}\leq \overline{L}_F(U)$ for some constant $\overline{L}_F(U)\in \mb{R}_{+}$ depending on $U$, such that for $P$-almost all $s\in \mc{S}$, 
    \begin{align*}
        \left\vert f_x(x;s) - f_x(y;s)\right\vert\leq \overline{L}_{f,U}(s)\Vert x - y\Vert_2
    \end{align*} 
    holds for all $x,y\in U$.
\end{assumption}
\begin{assumption}\label{asm:sppa2}
    The objective function $\phi$ satisfies the quadratic growth condition on $\mc{X}^{*}$ globally, that is, there exists $\overline{c}_1>0$ such that 
    \begin{align*}
        \phi(x)\geq \phi^{*} + \overline{c}_1 \mr{dist}\left(x,\mc{X}^{*}\right)^2
    \end{align*}
    holds for all $x\in \mb{R}^d$.
\end{assumption}
\begin{remark}\label{remark:lipqgmknorm}
    Suppose that Assumption \ref{asm:mk} holds for linear mappings $\{M_k\}_{k\in \mb{Z}_{+}}$. 
    \begin{enumerate}
        \item Let Assumption \ref{asm:stomodel} hold. Denote by $\eta\triangleq L_\infty^2\cdot\oleta$. Then for $P$-almost all $s\in \mc{S}$, the estimate 
        \begin{align*}
            \norm{f_x\left(y,s\right)-f\left(y,s\right)} \leq \frac{\eta}{2}\Norm{y-x}_{M_k}^2\quad \mbox{for all}\ x,y\in V
        \end{align*}
        holds for all $k\in \mb{Z}_{+}$.
        \item Let Assumption \ref{asm:sppa1} hold. Fix any bounded subset $U\subset V$. Denote by $L_{f,U}(s)\triangleq L_\infty\cdot\overline{L}_{f,U}(s)$ and $L_F(U)\triangleq L_\infty \cdot \overline{L}_F(U)$. Then, we have $\sqrt{\mb{E}_{s\sim P}[L_{f,U}(s)^2]}\leq L_F(U)$ and for $P$-almost all $s\in \mc{S}$, the estimate
        \begin{align*}
            \left\vert f_x(x;s) - f_x(y;s)\right\vert\leq L_{f,U}(s)\Vert x - y\Vert_{M_k}\quad \mbox{for all}\ x,y\in U
        \end{align*}
        holds for all $k\in \mb{Z}_{+}$.
        \item Let Assumption \ref{asm:sppa2} hold. Denote by $c_1\triangleq \mu_\infty^2\cdot\overline{c}_1$. Then, we have 
        \begin{align*}
            \phi(x)\geq \phi^{*} + c_1 \mr{dist}_{M_k}(x,\mc{X}^{*})^2\quad \mbox{for all}\ x\in \mb{R}^d
        \end{align*}
        for all $k\in \mb{Z}_{+}$. 
    \end{enumerate}
    It follows that under Assumption \ref{asm:mk}, both the quadratic error and the local Lipschitz continuity condition on the component model functions $f_x(\cdot;s)$, as well as the quadratic growth condition on the objective function $\phi$, are satisfied in the $M_k$-norm with constant independent of $k$.
\end{remark}
Table \ref{tab:relationthmasm} summarizes the connections between the aforementioned assumptions and the main convergence results to be established in Section \ref{sec:convergenceispp}. 
\begin{table}[!htbp]
    \centering
    \caption{Relationships between the main convergence results and the associated assumptions.}
    \label{tab:relationthmasm}
    \vspace{0.3em}
    \setlength{\tabcolsep}{4pt}
    \renewcommand{\arraystretch}{1.3}
    \begin{tabular*}{\textwidth}{@{\extracolsep\fill}cccc@{}}
        \toprule[1.1pt]
        \multirow{2}{*}{} & \multirow{2}{*}{\begin{tabular}[c]{@{}c@{}}Stability\\ ($\sup_{k\in \mb{Z}_{+}}\Vert x_k\Vert_2 < \infty$ a.s.)\end{tabular}} & \multicolumn{2}{c}{Convergence Rate} \\
        \cmidrule(lr){3-4}
         &  & Without Quadratic Growth  & With Quadratic Growth \\
        \midrule[1.05pt]
        \multirow{3}{*}{Weakly Convex} & \begin{tabular}[c]{@{}c@{}}\textbf{Theorem \ref{thm:almostsureboundednessadpstep1}}\\ (Asm \ref{asm:stomodel}-\ref{asm:sppaall}, \ref{asm:sublevelset}, \ref{asm:generalizedlip2})\end{tabular} & \multirow{6}{*}{\begin{tabular}[c]{@{}c@{}}\textbf{Theorem \ref{thm:sippmdiminishing}}\\ (Asm \ref{asm:stomodel}-\ref{asm:sppa1} \&\\ $\sup_{k\in \mb{Z}_{+}}\Vert x_k\Vert_2 < \infty$)\\ \textbf{Theorem \ref{thm:sippmdiminishing2}}\\ (Asm \ref{asm:stomodel}-\ref{asm:sppa1} \&\\ $\sup_{k\in \mb{Z}_{+}}\Vert x_k\Vert_2 < \infty$)\end{tabular}} & \multirow{6}{*}{\begin{tabular}[c]{@{}c@{}}\textbf{Theorem \ref{thm:sippmconstantquadratic}}\\ (Asm \ref{asm:stomodel}-\ref{asm:sppa2} \&\\ $\sup_{k\in \mb{Z}_{+}}\Vert x_k\Vert_2 < \infty$)\\ \textbf{Theorem \ref{thm:sippmdiminishingquadratic}}\\ (Asm \ref{asm:stomodel}-\ref{asm:sppa2} \&\\ $\sup_{k\in \mb{Z}_{+}}\Vert x_k\Vert_2 < \infty$)\end{tabular}} \\
         & \begin{tabular}[c]{@{}c@{}}\textbf{Theorem \ref{thm:almostsureboundednessadpstep2}}\\ (Asm \ref{asm:stomodel}-\ref{asm:sppaall}, \ref{asm:sublevelset}, \ref{asm:generalizedlip3})\end{tabular} &  &  \\
        \cmidrule(lr){1-2}
        Convex & \begin{tabular}[c]{@{}c@{}}\textbf{Theorem \ref{thm:sippstabilityconvex}}\\ (Asm \ref{asm:stomodel}-\ref{asm:sppaall}, \ref{asm:sppbdoptset} \& \\ $\oleta=\oltau=0$)\end{tabular} &  & \\
        \bottomrule[1.1pt]
    \end{tabular*}
\end{table}
Abbreviations used in the table are as follows: 
\begin{itemize}
    \item[*] \textbf{Asm \ref{asm:stomodel}}: $\oleta$-quadratic error of model function $f_x(\cdot;s)$.
    \item[*] \textbf{Asm \ref{asm:mk}}: preconditioner $M_k$.
    \item[*] \textbf{Asm \ref{asm:sppaall}}: $\oltau$-weak convexity of $\varphi_x(\cdot;s) \triangleq f_x(\cdot;s) + r(\cdot)$.
    \item[*] \textbf{Asm \ref{asm:sppa1}}, \textbf{\ref{asm:generalizedlip2}}, \textbf{\ref{asm:generalizedlip3}}: local Lipschitz condition defined in Assumption \ref{asm:sppa1}, \ref{asm:generalizedlip2}, \ref{asm:generalizedlip3}.
    \item[*] \textbf{Asm \ref{asm:sppa2}}: quadratic growth of $\phi$.
    \item[*] \textbf{Asm \ref{asm:sublevelset}}: $\phi$ bounded from below \& bounded sublevel sets.
    \item[*] \textbf{Asm \ref{asm:sppbdoptset}}: uniformly bounded subdifferential on the optimal solution set.
    \item[*] \textbf{Stability}: $\sup_{k\in \mb{Z}_{+}}\Vert x_k\Vert_2 < \infty$ with probability $1$.
\end{itemize}

\section{Analysis of Model-Based Inexact Stochastic Preconditioned Proximal Point Algorithm}\label{sec:convergenceispp}
In this section, we establish the convergence properties of the model-based ispPPA, including the stability, almost sure convergence, and convergence rate guarantees for the weakly convex case. 

Fix any $k\in \mb{Z}_{+}$. Given $\epsilon_k\geq 0$, we consider the following three criteria and say that $x_{k+1}$ is obtained according to (\ref{equ:sippm}) if it satisfies one of these criteria:
\begin{subequations}
    \begin{align}
        &\left\Vert x_{k+1} - \mr{prox}^{M_k}_{\alpha_k\olvphi_{x_k}(\cdot;S_k^{1:m})}(x_k)\right\Vert_{M_k} \leq \epsilon_k,\tag{SCA}\label{equ:criteriaa}\\
        &\Phi_{\alpha_k,x_k}(x_{k+1};S_k^{1:m}) - \Phi_{\alpha_k,x_k,S_k^{1:m}}^{*}\leq \frac{1-\tau \alpha_k}{2\alpha_k}\epsilon_k^2,\tag{SCB}\label{equ:criteriab}\\
        &\mr{dist}_{M_k^{-1}}\left(0,\partial \Phi_{\alpha_k,x_k}(x_{k+1};S_k^{1:m})\right)\leq \frac{1-\tau \alpha_k}{\alpha_k}\epsilon_k,\tag{SCC}\label{equ:criteriac}
    \end{align}
\end{subequations}
where constant $\tau \triangleq L_\infty^2\cdot\oltau$, function $\Phi_{\alpha_k,x_k}(\cdot;S_k^{1:m})\colon \mb{R}^d\to \pminf$ is defined by\footnote{Let Assumption \ref{asm:stomodel} and \ref{asm:sppaall} hold, and let $x\in V$. For each $s\in \mc{S}$ and all $y\in \mb{R}^d$, we define 
\begin{align*}
    \tilde{\varphi}_x(y;s) \triangleq \begin{cases}
        \varphi_x(y;s)\quad &\mbox{if}\ y\in \dom(r),\\
        \infty\quad &\mbox{otherwise}.
    \end{cases}
\end{align*}
Then, for $P$-almost all $s\in \mc{S}$, $\tilde{\varphi}_x(\cdot;s)\colon \mb{R}^d\to \pminf$ is proper, closed and $\oltau$-weakly convex, and satisfies $\tilde{\varphi}_x(y;s) = \varphi_x(y;s)$ for all $y\in V$, since $f_x(\cdot;s)$ is defined on $V$ and $V\supset \dom(r)$. For brevity, we write $\varphi_{x}(\cdot;s)$ for the extended-value function $\tilde{\varphi}_{x}(\cdot;s)$.} 
\begin{align*}
    \Phi_{\alpha_k,x_k}\left(x;S_k^{1:m}\right) &\triangleq \olvphi_{x_k}\left(x;S_k^{1:m}\right) + \frac{1}{2\alpha_k}\left\Vert x - x_k\right\Vert_{M_k}^2\quad \mbox{for all}\ x\in \mb{R}^d,
\end{align*}
and $\Phi_{\alpha_k,x_k,S_k^{1:m}}^{*}\triangleq \min_{x\in \mb{R}^d}\ \Phi_{\alpha_k,x_k}\left(x;S_k^{1:m}\right)$ denotes the minimum value. For simplicity, we denote the proximal point $\mr{prox}^{M_k}_{\alpha_k \olvphi_{x_k}(\cdot;S_k^{1:m})}(x_k)$ by $\tilde{x}_{k+1}$. Under Assumption \ref{asm:mk} and \ref{asm:sppaall}, if $\alpha_k\in (0,\tau^{-1})$, then by leveraging the fact that $\Norm{\cdot}_2 \leq L_\infty \Norm{\cdot}_{M_k}$, $\olvphi_{x_k}(\cdot;S_k^{1:m})$ is $\oltau$-weakly convex and $\frac{1}{\alpha_k}M_k(x_k - \tilde{x}_{k+1})\in \partial \olvphi_{x_k}(\tilde{x}_{k+1};S_k^{1:m})$, we can show that
\begin{equation}\label{equ:relationinexact}
    \begin{split}
        \text{Criterion (\ref{equ:criteriac})}\quad \Rightarrow \quad \text{Criterion (\ref{equ:criteriab})} \quad \Rightarrow \quad \text{Criterion (\ref{equ:criteriaa})}.
    \end{split}
\end{equation}
Therefore, it is sufficient to establish the convergence results for the model-based ispPPA with $\{x_k\}$ satisfying criterion (\ref{equ:criteriaa}). Due to (\ref{equ:relationinexact}), the derived conclusions also hold for the model-based ispPPA with $\{x_k\}$ satisfying either criterion (\ref{equ:criteriab}) or criterion (\ref{equ:criteriac}). 

Unless otherwise stated, we will use the following notations for simplicity. For any $\beta\in \mb{R}$, the function  $\varsigma_\beta\colon \mb{R}_{++}\to \mb{R}$ is defined by 
\begin{equation}\label{equ:varphi}
    \varsigma_\beta(t)\triangleq \begin{cases}
        \frac{t^\beta - 1}{\beta}\quad &\mbox{if}\ \beta \neq 0,\\
        \ln(t)\quad &\mbox{if}\ \beta = 0,
    \end{cases}
\end{equation}
for all $t\in \mb{R}_{++}$. For the open convex set $V$ specified in Assumption \ref{asm:stomodel} and any bounded subset $U\subset V$, we assume without loss of generality that $\mathbb{E}_{s\sim P}[L_{f,U}(s)^2]>0$ and set
\begin{equation}\label{equ:constantminibatch}
    \omega_{f,U} \triangleq 1 - \frac{\mr{Var}\left(L_{f,U}(s)\right)}{\mb{E}_{s\sim P}\left[L_{f,U}(s)^2\right]}\quad \mbox{and}\quad \rho_{f,m,U} \triangleq \left(\sqrt{\frac{1+(m-1)\omega_{f,U}}{m}}+1\right)^2,
\end{equation}
where $L_{f,U}(s)\triangleq L_\infty\cdot \olL_{f,U}(s)$ for each $s\in \mc{S}$.

\subsection{Stability and Almost Sure Convergence}\label{subsec:sippstability}
Before establishing convergence rate guarantees for the model-based ispPPA, we first identify conditions that guarantee the almost sure boundedness of the iterates $\{x_k\}$ generated by (\ref{equ:sippm}). As mentioned in Section \ref{sec:introduction}, the boundedness under mild assumptions is crucial, since our analysis is established without imposing any global Lipschitz condition on the component model functions $f_x(\cdot;s)$ and the regularizer $r$. 

\subsubsection{Weakly Convex Case}
For the weakly convex case, inspired by \cite{gao2024stochastic}, we establish sufficient conditions to ensure the stability of the model-based ispPPA with regularization stepsize strategies. The main results are presented in Theorem \ref{thm:almostsureboundednessadpstep1} and \ref{thm:almostsureboundednessadpstep2}. 

The stability analysis of the model-based ispPPA poses distinct challenges compared with its deterministic counterpart. In the deterministic setting, the objective function sequence $\{\phi(x_k)\}$ is typically nonincreasing, so the boundedness of the iterates follows directly from Assumption \ref{asm:sublevelset}. In contrast, in the stochastic setting, the model-based ispPPA incorporates the stochastic approximation term $\olvphi_{x_k}(\cdot;S_k^{1:m})$ in place of the function $\mb{E}_{s\sim P}[\varphi_{x_k}(\cdot;s)]$, and hence the sequence $\{\phi(x_k)\}$ (or its expectation $\{\mb{E}[\phi(x_k)]\}$) is no longer guaranteed to be monotone. As a result, Assumption \ref{asm:sublevelset} alone is insufficient to guarantee the stability. Meanwhile, when analyzing the almost sure convergence, it is essential to choose a suitable measure of stationarity. Following the framework of \cite{davis2019stochastic}, we adopt the gradient of the Moreau envelope as a continuous measure of stationarity in the weakly convex setting. We first introduce some necessary assumptions for the key stability result in Theorem~\ref{thm:almostsureboundednessadpstep1}.
\begin{assumption}\label{asm:sublevelset}
    The objective function $\phi$ is bounded from below and has bounded sublevel sets, i.e., $\phi^{*}>-\infty$ and 
    the sublevel set $\mc{L}_t\triangleq \{x\in \mb{R}^d\ \vert\ \phi(x)\leq t\}$ is bounded for any $t\in \mb{R}$.
\end{assumption}
\begin{assumption}\label{asm:generalizedlip2}
    For the open convex set $V$ specified in Assumption \ref{asm:stomodel}, there exists a nondecreasing function $\glip\colon \mb{R}_{+}\to \mb{R}_{+}$ and a measurable function $\overline{L}_{f}\colon \mc{S}\to \mb{R}_{+}$, satisfying $\sqrt{\mb{E}_{s\sim P}[\overline{L}_{f}(s)^2]}\leq \overline{L}_F$ for some constant $\overline{L}_F\in \mb{R}_{+}$, such that for $P$-almost all $s\in \mc{S}$, 
    \begin{align*}
        \left\vert f_x(x;s) - f_x(y;s)\right\vert\leq \overline{L}_{f}(s)\glip\left(\Norm{x}_2\right)\Vert x - y\Vert_2
    \end{align*} 
    for all $x,y\in V$.
\end{assumption}
\begin{theorem}\label{thm:almostsureboundednessadpstep1}
    Let Assumption \ref{asm:stomodel}, \ref{asm:mk}, \ref{asm:sppaall}, \ref{asm:sublevelset} and \ref{asm:generalizedlip2} hold, and let the iterates $\{x_k\}$ be generated by the model-based ispPPA (Algorithm \ref{algo:ispppa}) with stepsizes $\{\alpha_k\}$ and parameters $\{\epsilon_k\}$. Denote by $\eta \triangleq L_\infty^2\cdot\oleta$ and $\tau \triangleq L_\infty^2\cdot\oltau$. Set $B_\rho \triangleq \sup_{k\in \mb{Z}_{+}}\rho_k^2$. Suppose that there exists $\gamma\in \mb{R}_{+}$ and $\olrho > B_\rho \cdot(\eta + \tau)$ such that 
    \begin{equation}\label{asm:epsilonkalphakglip}
        \begin{split}
            &\alpha_k = \frac{\olalpha_k}{\max\lrbrackets{1,\mr{G}_{\mr{lip}}(\Vert x_k\Vert_2)}}\quad \mbox{and}\quad \epsilon_k = \gamma \alpha_k^2\quad \mbox{for all}\ k\in \mb{Z}_{+},\\\mbox{where}\quad &\sup_{k\in \mb{Z}_{+}}\olalpha_k < \left(\olrho\max\lrbrackets{2,B_\rho}+\eta\right)^{-1}\quad \mbox{and}\quad \sum_{k=1}^\infty \olalpha_k^2 < \infty.
        \end{split}
    \end{equation}
    Then, $\sup_{k\in \mb{Z}_{+}}\Norm{x_k}_2 < \infty$ holds with probability $1$. 
    Furthermore, if $\sum_{k=1}^\infty \olalpha_k = \infty$, we have 
    \begin{equation}\label{equ:infasmegrad}
        \lim\inf_{k\to \infty} \Norm{\nabla \mephi{k}\left(x_k\right)}_2 = 0\quad \mbox{with probability}\ 1,\quad \mbox{where}\ \mephi{k}(\cdot) \triangleq e^{M_{k-1}}_{(\olrho)^{-1}\phi}(\cdot).
    \end{equation}
\end{theorem}
Assumption \ref{asm:generalizedlip2} is inspired by the condition \cite[C1]{gao2024stochastic}. Under this assumption, unlike \cite[Theorem~4.1]{gao2024stochastic}, which requires the regularizer $r$ to be globally Lipschitz continuous on its domain, Theorem \ref{thm:almostsureboundednessadpstep1} guarantees the almost sure boundedness of the iterates without imposing any Lipschitz condition on $r$. 
Theorem \ref{thm:almostsureboundednessadpstep1} also extends the stability result of \cite[Theorem~4.1]{gao2024stochastic} for the standard stochastic model-based algorithm to its preconditioned and inexact variants. Although this theoretical guarantee bypasses the requirement for the global Lipschitz continuity, the local Lipschitz condition in Assumption \ref{asm:generalizedlip2} depends only on the current iterate $x$, which may limit its applicability in practice. To further address this issue, we propose a mild variant of the local Lipschitz condition, and the corresponding convergence guarantee is established in Theorem \ref{thm:almostsureboundednessadpstep2}.
\begin{assumption}\label{asm:generalizedlip3}
    For the open convex set $V$ specified in Assumption \ref{asm:stomodel}, there exists a continuous nondecreasing function $\glip\colon \mb{R}_{+}\to \mb{R}_{+}$ and a measurable function $\overline{L}_{\varphi}\colon \mc{S}\to \mb{R}_{+}$, satisfying $\sqrt{\mb{E}_{s\sim P}[\overline{L}_{\varphi}(s)^2]}\leq \overline{L}_\phi$ for some constant $\overline{L}_\phi\in \mb{R}_{+}$, such that for $P$-almost all $s\in \mc{S}$, 
    \begin{align*}
        \left\vert \varphi_x(y;s) - \varphi_x(z;s)\right\vert\leq \overline{L}_{\varphi}(s)\glip\left(\max\lrbrackets{\Norm{x}_2,\Norm{y}_2}\right)\Vert y - z\Vert_2
    \end{align*} 
    for all $x,y\in V$ and all $z$ in a neighborhood of $y$, where $\varphi_x(\cdot;s) \triangleq f_x(\cdot;s)+r(\cdot)$.
\end{assumption}
\begin{theorem}\label{thm:almostsureboundednessadpstep2}
    Let Assumption \ref{asm:stomodel}, \ref{asm:mk}, \ref{asm:sppaall}, \ref{asm:sublevelset} and \ref{asm:generalizedlip3} hold, and let the iterates $\{x_k\}$ be generated by the model-based ispPPA (Algorithm \ref{algo:ispppa}) with stepsizes $\{\alpha_k\}$ and parameters $\{\epsilon_k\}$. Denote by $\eta \triangleq L_\infty^2\cdot\oleta$ and $\tau \triangleq L_\infty^2\cdot\oltau$. Set $B_\rho \triangleq \sup_{k\in \mb{Z}_{+}}\rho_k^2$. 
    Suppose that there exists $\gamma\in \mb{R}_{+}$ and $\olrho > B_\rho \cdot(\eta + 2\tau)$ such that 
    \begin{equation}\label{asm:epsilonkalphakglip2}
        \begin{split}
            &\alpha_k = \frac{\olalpha_k}{\max\lrbrackets{1,\mr{G}_{\mr{lip}}(\Vert x_k\Vert_2)}}\quad \mbox{and}\quad \epsilon_k = \gamma \alpha_k^2\quad \mbox{for all}\ k\in \mb{Z}_{+},\\\mbox{where}\quad &\sup_{k\in \mb{Z}_{+}}\olalpha_k < \left(2\olrho\right)^{-1}\quad \mbox{and}\quad \sum_{k=1}^\infty \olalpha_k^2 < \infty.
        \end{split}
    \end{equation}
    Then, all conclusions of Theorem \ref{thm:almostsureboundednessadpstep1} remain valid.
\end{theorem}
\begin{remark}\label{remark:alphaklbubwcvx}
    Under the setting of Theorem \ref{thm:almostsureboundednessadpstep1} or \ref{thm:almostsureboundednessadpstep2}, when the stepsizes are chosen as $\alpha_k = \frac{\olalpha_k}{\max\{1,\glip(\Vert x_k\Vert_2)\}}$ with $\olalpha_k = \olalpha_0k^{-\beta}$ for some $\beta\in (\frac{1}{2},1]$ and $\olalpha_0$ sufficiently small (see (\ref{asm:epsilonkalphakglip}) or (\ref{asm:epsilonkalphakglip2})), with probability $1$, we have $\sup_{k\in \mb{Z}_{+}}\Vert x_k\Vert_2 < \infty$ and there exists $\alphaol,\alphaou\in (0,\olalpha_0]$ such that $\alphaol k^{-\beta}\leq \alpha_k \leq \alphaou k^{-\beta}$ holds for all $k\in \mb{Z}_{+}$. 
\end{remark}
The motivation for introducing Assumption \ref{asm:generalizedlip3} lies in its  ability to broaden the admissible class of model functions for a wide range of stochastic optimization problems. Although Assumption \ref{asm:generalizedlip3} imposes an additional continuity requirement on the growth function $\glip$ compared with Assumption \ref{asm:generalizedlip2}, this condition is mild as it applies to many common choices of $\glip$, such as polynomial functions. For instance, consider $f(\cdot;s)\triangleq h(c(\cdot;s);s)$ and $r\equiv 0$, where $h\colon \mb{R}^m\to \mb{R}$ is convex and $L_h$-Lipschitz continuous, and $c(\cdot;s)\colon \mb{R}^d\to \mathbb{R}^m$ is continuously differentiable with a Jacobian $\nabla c$ that is $L_c$-Lipschitz continuous. Commonly used model functions in this setting include:
\begin{itemize}
    \item {Stochastic subgradient}: $f_x(y;s) = f(x;s) + \lrangle{\nabla c(x;s)^\top v(x;s) ,y-x}$ with $v(x;s)\in \partial h(c(x;s);s)$;
    \item {Stochastic prox-linear}: $f_x(y;s) = h(c(x;s) + \nabla c(x;s)(y-x);s)$;
    \item {Stochastic proximal point}: $f_x(y;s) = h(c(y;s))$.
\end{itemize}
Combining the convexity and Lipschitz continuity of $h$ with the Lipschitz continuity of $\nabla c$ implies that both the stochastic subgradient and prox-linear model functions satisfy Assumption \ref{asm:generalizedlip2} and \ref{asm:generalizedlip3}, with a continuous growth function $\glip(t) \triangleq t$. Furthermore, by leveraging the weak convexity of $f(\cdot;s)$, one can deduce from \cite[Theorem~2.4]{clarke1998nonsmooth} that the stochastic proximal point model functions also satisfy Assumption \ref{asm:generalizedlip3} with a continuous growth function $\glip(t) \triangleq t$. However, when the operator norm $\Norm{\nabla c(\cdot;s)}_2$ is not uniformly bounded on $\mb{R}^d$ (e.g., $c(x;s) = \lrangle{a,x}^2-b$ in the robust phase retrieval problem), the corresponding stochastic proximal point model functions may violate Assumption \ref{asm:generalizedlip2}.

\paragraph*{Proof of Theorem \ref{thm:almostsureboundednessadpstep1} and \ref{thm:almostsureboundednessadpstep2}.} 
To prove Theorem \ref{thm:almostsureboundednessadpstep1} and \ref{thm:almostsureboundednessadpstep2}, we first recall the so-called supermartingale convergence lemma, the proof of which can be found in \cite[Theorem~1]{robbins1971convergence}. 
\begin{lema}\label{lema:supermartingaleconvergence}
    Let $Z_k$, $\beta_k$, $\xi_k$ and $\zeta_k$ be nonnegative random variables adapted to the filtration $\mc{F}_k$ and satisfying 
    \begin{align*}
        \mb{E}\left[Z_{k+1}\ \vert\ \mc{F}_k\right] \leq \left(1+\beta_k\right) Z_k + \xi_k - \zeta_k
    \end{align*} 
    for all $k\in \mb{Z}_{+}$. Then on the event that $\sum_{k=1}^\infty \beta_k < \infty$ and $\sum_{k=1}^\infty \xi_k < \infty$, with probability $1$, we have $\sum_{k=1}^\infty \zeta_k < \infty$ and there exists a random variable $Z_\infty$ such that $\lim_{k\to \infty}Z_k = Z_{\infty}$. 
\end{lema}
The proofs of Theorem \ref{thm:almostsureboundednessadpstep1} and \ref{thm:almostsureboundednessadpstep2} rely on establishing the validity of the following recursive relation, with its proof deferred to Appendix \ref{apx:proofclaimsippinequhatx1tvalueglip}.
\begin{claim}\label{claim:sippinequhatx1tvalueglip}
    Consider the setting of Theorem \ref{thm:almostsureboundednessadpstep1} or \ref{thm:almostsureboundednessadpstep2}. Denote by $\olx_k \triangleq \prox_{(\olrho)^{-1}\phi}^{M_{k-1}}(x_k)$ and $\mephi{k}(\cdot) \triangleq e_{(\olrho)^{-1}\phi}^{M_{k-1}}(\cdot)$ for each $k\in \mb{Z}_{+}$. Then, $\olx_k$ is well-defined, $\mephi{k}(x_k) - \phi^{*}\geq 0$, and there exist constants $\theta_1,\theta_2>0$ such that 
    \begin{equation}\label{equ:medecreaseglip}
        \begin{split}
            &\mb{E}_k\left[\nu_{k+1}^{-2}\cdot\left(\mephi{k+1}\left(x_{k+1}\right) - \phi^{*}\right)\right]\\
            \leq\ &\nu_k^{-2}\cdot\left(\mephi{k}\left(x_k\right) - \phi^{*}\right) - \theta_1 \alpha_k \cdot \nu_k^{-2}\Norm{\olrho \left(x_k - \olx_k\right)}_{M_{k-1}}^2 + \theta_2 \olalpha_k^2
        \end{split}
    \end{equation}
    for all $k\in \mb{Z}_{+}$.
\end{claim}
In view of Claim \ref{claim:sippinequhatx1tvalueglip}, the proof arguments of these two theorems are essentially identical. For brevity, we present only the proof of Theorem \ref{thm:almostsureboundednessadpstep1}.

\begin{proof}[Proof of Theorem \ref{thm:almostsureboundednessadpstep1}]
Using Assumption \ref{asm:mk}, we have
\begin{align*}
    B_\rho\triangleq \sup_{k\in \mb{Z}_{+}}\rho_k^2 \leq \nu_\infty^2 < \infty.
\end{align*}
Fix any $k\in \mb{Z}_{+}$. Claim \ref{claim:sippinequhatx1tvalueglip} implies that $\olx_k \triangleq \prox_{(\olrho)^{-1}\phi}^{M_{k-1}}(x_k)$ is well-defined, 
\begin{equation}\label{equ:lbme}
    \mephi{k}\left(x_{k}\right) - \phi^{*} \geq 0,
\end{equation}
and the estimate (\ref{equ:medecreaseglip}) holds for all $k\in \mb{Z}_{+}$. 

\underline{Proof of the stability.} 
The estimate (\ref{equ:lbme}) combined with $\inf_{k\in \mb{Z}_{+}}\nu_k \geq 1$ yields that both $\{\nu_k^{-2}\cdot(\mephi{k}(x_k) - \phi^{*})\}$ and $\{\alpha_k\nu_k^{-2}\Norm{\olrho (x_k - \olx_k)}_{M_{k-1}}^2\}$ are nonnegative random variables adapted to the filtration $\mc{F}_k$. Note that $\sum_{k=1}^\infty \theta_2 \olalpha_k^2 < \infty$ and $\theta_1 > 0$. By utilizing (\ref{equ:medecreaseglip}) with $Z_k = \nu_k^{-2}\cdot(\mephi{k}(x_k) - \phi^{*})$, $\beta_k = 0$, $\xi_k = \theta_2 \olalpha_k^2$, and $\zeta_k = \theta_1\alpha_k\nu_k^{-2}\Vert \olrho(x_k-\olx_k)\Vert_{M_{k-1}}^2$ in Lemma \ref{lema:supermartingaleconvergence}, one can conclude that 
\begin{subequations}
    \begin{align}
        &\mbox{with probability}\ 1,\quad \sum_{k=1}^\infty \alpha_k \nu_k^{-2}\Norm{\olrho \left(x_k - \olx_k\right)}_{M_{k-1}}^2 < \infty\quad \mbox{and}\label{equ:assumme}\\
        &\hspace{2em}\mbox{sequence}\ \lrbrackets{\nu_k^{-2}\cdot\left(\mephi{k}\left(x_k\right) - \phi^{*}\right)}\ \mbox{converges to some finite value}.\label{equ:convergemephi} 
    \end{align}
\end{subequations}
On the event that (\ref{equ:convergemephi}) holds, using Assumption \ref{asm:mk} and the definition of $\phi^{*}$, we have 
\begin{align*}
        \sup_{k\in \mb{Z}_{+}}\lrbrackets{e_{\left(\olrho\right)^{-1}\phi}(x_k)-\phi^{*}} 
        &\leq \sup_{k\in \mb{Z}_{+}}\lrbrackets{\inf_{x\in \mb{Z}_{+}}\lrbrackets{\phi\left(x\right)+\frac{\olrho L_{k-1}^2}{2}\Norm{x-x_k}_{M_{k-1}}^2 - \phi^{*}}}\\
        &\leq \sup_{k\in \mb{Z}_{+}} \lrbrackets{L_{k-1}^2\cdot \inf_{x\in \mb{Z}_{+}}\lrbrackets{\phi\left(x\right)+\frac{\olrho }{2}\Norm{x-x_k}_{M_{k-1}}^2 - \phi^{*}}}\\
        &= \sup_{k\in \mb{Z}_{+}} \lrbrackets{L_{k-1}^2\cdot(\mephi{k}(x_k) - \phi^{*})}\\
        &\leq \left(\nu_\infty L_\infty\right)^2\cdot \sup_{k\in \mb{Z}_{+}}\lrbrackets{\nu_k^{-2}\cdot\left(\mephi{k}(x_k) - \phi^{*}\right)} < \infty,
\end{align*}
which implies
\begin{equation}\label{equ:sublevelme}
    \sup_{k\in \mb{Z}_{+}} e_{\left(\olrho\right)^{-1}\phi}(x_k) < \infty\quad \mbox{on the event that}\ (\ref{equ:convergemephi})\ \mbox{holds}.
\end{equation}
Since $\phi$ has bounded sublevel sets by Assumption \ref{asm:sublevelset}, according to \cite[Lemma~2]{li2024revisiting}, the Moreau enevlope $e_{\left(\olrho\right)^{-1}\phi}$ has bounded sublevel sets. It then follows immediately from (\ref{equ:sublevelme}) that 
\begin{align*}
        \mb{P}\lrbrackets{\sup_{k\in \mb{Z}_{+}}\Norm{x_k}_2 < \infty} \geq \mb{P}\lrbrackets{\mbox{the sequence}\ \lrbrackets{\nu_k^{-2}\cdot\left(\mephi{k}(x_k) - \phi^{*}\right)}\ \mbox{converges}},
\end{align*}
which combined with (\ref{equ:convergemephi}) confirms that the sequence $\{x_k\}$ is almost surely bounded. 

\underline{Proof of the almost sure convergence.} Assume further that $\sum_{k=1}^\infty \olalpha_k = \infty$. On the event that $\sup_{k\in \mb{Z}_{+}}\Norm{x_k}_2 < \infty$, we have 
\begin{align*}
    0\leq B_{\glip}\triangleq \sup_{k\in \mb{Z}_{+}}\glip\left(\Norm{x_k}_2\right) \leq \glip\left(\sup_{k\in \mb{Z}_{+}}\Norm{x_k}_2\right) < \infty
\end{align*}
as the function $\glip\colon \mb{R}_{+}\to \mb{R}_{+}$ is nondecreasing. It then follows from (\ref{asm:epsilonkalphakglip}) that
\begin{equation}\label{equ:lbalphakwcvx}
    \alpha_k \geq \frac{\olalpha_k}{1+\glip\left(\Norm{x_k}_2\right)} \geq \frac{\olalpha_k}{1+B_{\glip}}\quad \mbox{for all}\ k\in \mb{Z}_{+},
\end{equation}
which combined with condition $\sum_{k=1}^\infty\olalpha_k = \infty$ yields 
\begin{equation}\label{equ:notsummablealphak}
    \sum_{k=1}^\infty \alpha_k = \infty\quad \mbox{with probability}\ 1.
\end{equation}
Denote by 
\begin{equation*}
    \begin{split}
        S_1 &\triangleq \lrbrackets{w\ \vert\ \lim\inf_{k\to \infty}\nu_k^{-1}\Norm{\olrho \left(x_k - \olx_k\right)}_{M_{k-1}} > 0},\\
        S_2 &\triangleq \lrbrackets{w\ \vert\ \sum_{k=1}^\infty \alpha_k \nu_k^{-2}\Norm{\olrho \left(x_k - \olx_k\right)}_{M_{k-1}}^2 < \infty},\quad \mbox{and}\quad 
        S_3 \triangleq \lrbrackets{w\ \vert\ \sum_{k=1}^\infty \alpha_k = \infty}.
    \end{split}
\end{equation*}
Consider an arbitrary sample $w\in S_1\cap S_2\cap S_3$ with associated realizations $\{x_k\}$. The condition $w\in S_1$ implies that 
\begin{align*}
    \theta_0 \triangleq \sup_{m\geq 1}\inf_{k\geq m} \nu_k^{-1}\Vert\olrho \left(x_k - \olx_k\right)\Vert_{M_{k-1}} = \lim\inf_{k\to \infty}\nu_k^{-1}\Vert\olrho \left(x_k - \olx_k\right)\Vert_{M_{k-1}} > 0.
\end{align*}
Using the definition of the supremum, we can choose $K\in \mb{Z}_{+}$ such that 
\begin{align*}
        \inf_{k\geq K} \nu_k^{-1}\Vert\olrho \left(x_k - \olx_k\right)\Vert_{M_{k-1}} > \frac{\theta_0}{2}\quad \Rightarrow\quad &\nu_k^{-1}\Vert\olrho \left(x_k - \olx_k\right)\Vert_{M_{k-1}} > \frac{\theta_0}{2}\quad \mbox{for all}\ k\geq K.
\end{align*}
Since $\theta_0 > 0$, we can thereby deduce from the preceding result with the condition $w\in S_2$ that
\begin{align*}
        \sum_{k=1}^\infty \alpha_k 
        &\leq \sum_{k=1}^{K-1}\alpha_k + \left(\frac{2}{\theta_0}\right)^2\cdot\sum_{k= K}^\infty\alpha_k \cdot \nu_k^{-2}\Vert\olrho \left(x_k - \olx_k\right)\Vert_{M_{k-1}}^2 < \infty,\\ 
\end{align*}
contradicting to the condition $w\in S_3$. Since the sample $w$ was chosen arbitrarily, we have $\mb{P}(S_1\cap S_2\cap S_3) = 0$. We claim that $\mb{P}(S_1)=0$. Suppose on contrary that $\mb{P}(S_1) > 0$. It follows from (\ref{equ:assumme}) and (\ref{equ:notsummablealphak}) that $\mb{P}(S_2^c) = \mb{P}(S_3^c) = 0$ and thus 
\begin{align*}
        \mb{P}\lrrbrackets{S_1 \cap S_2\cap S_3} &= 1 - \mb{P}\lrrbrackets{S_1^2\cup S_2^c\cup S_3^c} \geq 1- \mb{P}\lrrbrackets{S_1^c} - \mb{P}\lrrbrackets{S_2^c} - \mb{P}\lrrbrackets{S_3^c}\\
        &= 1- \mb{P}\lrrbrackets{S_1^c} = \mb{P}\lrrbrackets{S_1} > 0,
\end{align*} 
which is a contradiction. Thus, we have 
\begin{equation}\label{equ:meinfas}
    \lim\inf_{k\to \infty}\nu_k^{-1}\Norm{\olrho \left(x_k - \olx_k\right)}_{M_{k-1}} = 0\quad \mbox{with probability}\ 1.
\end{equation}
Fix any $k\in \mb{Z}_{+}$. Since $M_k$ is self-adjoint and positive-definite, we have 
\begin{align*}
    \frac{\Norm{x}_2^2}{\Norm{x}_{M_k^{-1}}^2} &= \frac{\lrangle{x,x}}{\lrangle{x,M_k^{-1}x}} 
        = \frac{\Norm{M_k^{-\frac{1}{2}}x}_{M_k}^2}{\Norm{M_k^{-\frac{1}{2}}x}_2^2}\quad \mbox{for all}\ 0\neq x\in \mb{R}^d,
\end{align*}
which combined with Assumption \ref{asm:mk} confirms
\begin{equation}\label{equ:normmkm1}
    L_k^{-1}\Norm{x}_{M_k^{-1}} \leq \Norm{x}_2 \leq \mu_k^{-1}\Norm{x}_{M_k^{-1}}\quad \mbox{for all}\ x\in \mb{R}^d\ \mbox{and}\ k\in \mb{Z}_{+}.
\end{equation}
Using the definition of $\olx_k$ and applying Lemma \ref{lema:welldefinedolxk} with $x= x_k$, we can conclude that 
\begin{equation}\label{equ:gradnormmkmeolxk}
    \Norm{\nabla \mephi{k}\left(x_k\right)}_{M_{k-1}^{-1}} = \Norm{\olrho \left(x_k - \olx_k\right)}_{M_{k-1}}.
\end{equation}
Combining (\ref{equ:gradnormmkmeolxk}) with (\ref{equ:normmkm1}) and condition $\mu_{k-1} \geq \mu_\infty > 0$ from Assumption \ref{asm:mk} yields 
\begin{equation}\label{equ:normgradme}
    \begin{split}
        \Norm{\nabla \mephi{k}\left(x_k\right)}_2 &
        \leq \mu_\infty^{-1}\cdot \Norm{\olrho \left(x_k - \olx_k\right)}_{M_{k-1}} \quad \mbox{for all}\ k\in \mb{Z}_{+}.
    \end{split}
\end{equation}
Since $\lim_{k\to \infty}\nu_k = \nu_\infty$ by Assumption \ref{asm:mk}, then (\ref{equ:infasmegrad}) follows immediately from the combination of (\ref{equ:meinfas}) with (\ref{equ:normgradme}). This completes the proof of Theorem \ref{thm:almostsureboundednessadpstep1}. 
\end{proof}

\subsubsection{Convex Case with Full Proximal Model Functions}
When considering the full proximal model functions in the convex setting (i.e., $\oleta = \oltau = 0$), we refer to Algorithm \ref{algo:ispppa} as the \emph{preconditioned isPPA}, which is a preconditioned variant of the isPPA analyzed in \cite{zhu2025tight}. In this case, by imposing the following condition as given in \cite[Assumption~2]{zhu2025tight}, we can establish the almost sure boundedness of the iterates generated by the preconditioned isPPA with square summable stepsizes. 
\begin{assumption}\label{asm:sppbdoptset}
    The optimal solution set $\mc{X}^{*}$ is nonempty. Denote by $\varphi(\cdot;s)\triangleq f(\cdot;s)+r(\cdot)$ for each $s\in \mc{S}$. There exists a scalar $\sigma_\phi\in \mb{R}_{++}$ such that for all $x^{*}\in \mc{X}^{*}$, it holds that 
    \begin{align*}
        \mb{E}_{s\sim P}\left[\Norm{\varphi^\prime(x^{*};s)}_2^2\right] \leq \sigma_\phi^2\quad \mbox{for all measurable selections}\  \varphi^\prime(x^{*};s)\in \partial \varphi(x^{*};s).
    \end{align*}
    For any fixed $x^{*}\in \mc{X}^{*}$ and $\phi^\prime(x^{*})\in \partial \phi(x^{*})$, there exists a measurable mapping $\varphi^\prime(x^{*};\cdot)$ such that 
    \begin{align*}
        \mb{E}_{s\sim P}\left[\varphi^\prime\left(x^{*};s\right)\right] = \phi^\prime(x^{*})\quad \mbox{with}\ \varphi^\prime\left(x^{*};s\right)\in \partial \varphi(x^{*};s)\ \mbox{for}\ P\mbox{-almost}\ s\in \mc{S},
    \end{align*}
    where $\mb{E}_{s\sim P}\left[\varphi^\prime\left(x^{*};s\right)\right]$ is defined as the integral $\int_{\mc{S}} \varphi^\prime(x^{*};s)\mr{d}\, P(s)$.
\end{assumption}
This assumption concerns the behavior of the subdifferential on the optimal solution set. As noted in \cite[Remark~2.1]{zhu2025tight}, this condition is mild and can be satisfied by typical finite-sum regularized regression problems involving real-valued nonnegative convex component functions $f(\cdot;s)$ and a real-valued convex coercive regularizer $r$. By utilizing the stability guarantee and further assuming the local Lipschitz condition in Assumption \ref{asm:sppa1}, we can prove the almost sure convergence for the preconditioned isPPA with stepsizes that are square summable but not summable. The main result is summarized in Theorem \ref{thm:sippstabilityconvex}. 
\begin{theorem}\label{thm:sippstabilityconvex}
    Let Assumption \ref{asm:stomodel}, \ref{asm:mk}, \ref{asm:sppaall} and \ref{asm:sppbdoptset} hold, and let the iterates $\{x_k\}$ be generated by the model-based ispPPA (Algorithm \ref{algo:ispppa}) with stepsizes $\{\alpha_k\}$ and parameters $\{\epsilon_k\}$. Suppose that $\oleta = \oltau = 0$ and there exists $\gamma\in \mb{R}_{+}$ such that 
    \begin{equation}\label{asm:epsilonkalphakcvx}
        \epsilon_k = \gamma \alpha_k^2\quad \mbox{for all}\ k\in \mb{Z}_{+},\quad \mbox{and}\quad \sum_{k=1}^\infty \alpha_k^2<\infty.
    \end{equation}
    Then, $\sup_{k\in \mb{Z}_{+}}\Norm{x_k}_2 < \infty$ holds with probability $1$. 
    If, in addition, Assumption \ref{asm:sppa1} holds and  $\sum_{k=1}^\infty \alpha_k = \infty$, then with probability $1$, there exists $x^{*}\in \mc{X}^{*}$ such that $\lim_{k\to \infty} \Norm{x_k - x^{*}}_2 = 0$.
\end{theorem}

\paragraph*{Proof of Theorem \ref{thm:sippstabilityconvex}.} 
To prove Theorem \ref{thm:sippstabilityconvex}, we introduce the following recursive relation, the proof of which is provided in Appendix \ref{apx:proofclaimsippstablebound}.
\begin{claim}\label{claim:sippstablebound}
    Consider the setting of Theorem \ref{thm:sippstabilityconvex}. Denote by $\tlx_{k+1} \triangleq \prox_{\alpha_k \olvphi_{x_k}(\cdot;S_k^{1:m})}^{M_k}(x_k)$ for each $k\in \mb{Z}_{+}$. Fix any $k\in \mb{Z}_{+}$. Then, $\tlx_{k+1}$ is well-defined and 
    \begin{equation}\label{equ:sippstablebound}
        \begin{split}
            &\mb{E}_k \left[\nu_{k+1}^{-2}\Norm{x_{k+1}-\olx}_{M_k}^2\right]\\ \leq\  &\left(1+\alpha_k^2\right)\nu_k^{-2}\Norm{x_k - \olx}_{M_{k-1}}^2 +   \left(1+\alpha_k^2\right)\alpha_k^2\left(\frac{L_k^2 \sigma_\phi^2}{m}+\gamma^2\right)\quad \mbox{for all}\ \olx\in \mc{X}^{*}.
        \end{split}
    \end{equation}
    If Assumption \ref{asm:sppa1} holds and there exists a bounded open convex subset of $V$, denoted by $U$, such that both $\{x_k\}$ and $\{\tlx_k\}$ are contained in $U$, we deduce
    \begin{equation}\label{equ:sippinequhatx1tvalue}
        \begin{split}
            &\mb{E}_k\left[\nu_{k+1}^{-2}\Norm{x_{k+1} - x}_{M_k}^2\right]\\ \leq\ 
            &\left(1+\alpha_k^2\right)\cdot\left(\nu_k^{-2}\Vert x_k - x\Vert_{M_{k-1}}^2 - 2\nu_{k+1}^{-2}\alpha_k\cdot \mb{E}_k\left[\phi\left(\tilde{x}_{k+1}\right) - \phi\left(x\right)\right]\right.\\
            &\hspace{7.5em}+\ \left.\nu_{k+1}^{-2}\left(\rho_{f,m,U}L_F(U)^2 + \gamma^2\right)\cdot\alpha_k^2\right)\quad \mbox{for all}\ x\in \mb{R}^d,
        \end{split}
    \end{equation}
    where $\rho_{f,m,U} \triangleq (\sqrt{\frac{1+(m-1)\omega_{f,U}}{m}}+1)^2$ for $\omega_{f,U} \triangleq 1 - \frac{\mr{Var}(L_{f,U}(s))}{\mb{E}_{s\sim P}[L_{f,U}(s)^2]}$\footnote{Without loss of generality, we may assume that $\mathbb{E}_{s\sim P}[L_{f,U}(s)^2]>0$.} and $L_{f,U}(s) \triangleq L_\infty\cdot \olL_{f,U}(s)$, and $L_F(U) \triangleq L_\infty \cdot\olL_F(U)$.
\end{claim}
\begin{proof}[Proof of Theorem \ref{thm:sippstabilityconvex}]
Use the same notations as in Claim \ref{claim:sippstablebound}.

\underline{Proof of the stability.} 
Fix any $\olx\in \mc{X}^{*}$. The assumption that $\sum_{k=1}^\infty \alpha_k^2 < \infty$ implies $\alpha_k \leq 1$ for $k$ sufficiently large, which combined with the fact that $\sup_{k\in \mb{Z}_{+}}L_k\leq L_\infty$ yields
\begin{align*}
    \sum_{k=1}^\infty \left(1+\alpha_k^2\right)\alpha_k^2 \left(\frac{L_k^2\sigma_\phi^2}{m} + \gamma^2\right) \leq \sum_{k=1}^\infty \left(1+\alpha_k^2\right)\alpha_k^2 \left(\frac{L_\infty^2\sigma_\phi^2}{m} + \gamma^2\right) < \infty.
\end{align*}
Applying Lemma \ref{lema:supermartingaleconvergence} with $Z_k = \nu_k^{-2}\Norm{x_k - \olx}_{M_{k-1}}^2$, $\beta_k = \alpha_k^2$, $\xi_k = (1+\alpha_k^2)\alpha_k^2 (\frac{L_k^2\sigma_\phi^2}{m} + \gamma^2)$ and $\zeta_k = 0$, we can deduce from Claim \ref{claim:sippstablebound} that 
\begin{equation}\label{equ:boundcvxsequence}
    \begin{split}
        &\mbox{with probability}\ 1,\\
        &\hspace{2em}\mbox{sequence}\ \lrbrackets{\nu_k^{-1}\Norm{x_k - \olx}_{M_{k-1}}}\ \mbox{converges to some finite value}.
    \end{split}
\end{equation}
Recall that $M_0 \triangleq I_d$ and $L_0\triangleq 1$. On the event that $\{\nu_k^{-1}\Norm{x_k - \olx}_{M_{k-1}}\}$ converges, using Assumption \ref{asm:mk}, one can conclude that 
\begin{align*}
    \sup_{k\in \mb{Z}_{+}}\Norm{x_k - \olx}_2 
        &\leq \nu_\infty L_\infty \cdot \sup_{k\in \mb{Z}_{+}} \left\{\nu_k^{-1}\Norm{x_k - \olx}_{M_{k-1}}\right\} < \infty,
\end{align*}
which combined with (\ref{equ:boundcvxsequence}) yields that the sequence $\{x_k\}$ is almost surely bounded. 

\underline{Proof of the almost sure convergence.} Assume further that Assumption \ref{asm:sppa1} holds and $\sum_{k=1}^\infty \alpha_k = \infty$. We first recall the following elementary inequality:
\begin{equation}\label{inequ:exp}
    1+t\leq \mr{exp}(t)\quad \mbox{for all}\ t\in \mb{R}.
\end{equation}
Based on the stability established above, it is sufficient to show that on the event that 
\begin{equation}\label{equ:boundxknsumsqsumalphak}
    \sup_{k\in \mb{Z}_{+}}\ \Norm{x_k}_2<\infty,\quad \sum_{k=1}^\infty \alpha_k = \infty\quad \mbox{and}\quad \sum_{k=1}^\infty \alpha_k^2 < \infty,
\end{equation}
the sequence $\{x_k\}$ converges to some $x^{*}\in \mc{X}^{*}$ with probability $1$. On the event that (\ref{equ:boundxknsumsqsumalphak}) holds, criterion (\ref{equ:criteriaa}) along with the fact that $\Norm{\cdot}_2\leq L_\infty\Norm{\cdot}_{M_k}$ implied by Assumption \ref{asm:mk} and the condition $\epsilon_k = \gamma \alpha_k^2$ by (\ref{asm:epsilonkalphakcvx}) implies 
\begin{equation}\label{equ:boundtlxkp1}
    \begin{split}
        \sup_{k\in \mb{Z}_{+}} \Norm{\tilde{x}_{k+1}}_2 &\leq \sup_{k\in \mb{Z}_{+}} \Norm{x_{k+1}}_2 + \sup_{k\in \mb{Z}_{+}} \Norm{\tlx_{k+1} - x_{k+1}}_2\\
        &\leq \sup_{k\in \mb{Z}_{+}} \Norm{x_{k+1}}_2 +  L_\infty\cdot \sup_{k\in \mb{Z}_{+}} \left\{\Norm{\tlx_{k+1} - x_{k+1}}_{M_k}\right\}\\ 
        &\leq \sup_{k\in \mb{Z}_{+}} \Norm{x_{k+1}}_2 + \gamma L_\infty \cdot \sup_{k\in \mb{Z}_{+}}\alpha_k^2 < \infty.
    \end{split}
\end{equation}
Thus, we denote by $U$ a bounded open convex subset of $V$ such that both $\{x_k\}$ and $\{\tlx_k\}$ are contained in $U$. It then follows from Claim \ref{claim:sippstablebound} that (\ref{equ:sippinequhatx1tvalue}) holds. 

\underline{\textit{*Step 1.}} Since $\inf_{k\in \mb{Z}_{+}}\nu_k \geq 1$, it follows from (\ref{equ:sippinequhatx1tvalue}) that 
\begin{equation}\label{equ:sippinequhatx1tvalueforproofzhutight}
    \begin{split}
        &\mb{E}_k\left[\nu_{k+1}^{-2}\Norm{x_{k+1} - x}_{M_k}^2\right]\\ \leq\ 
        &\left(1+\alpha_k^2\right)\cdot\left(\nu_k^{-2}\Vert x_k - x\Vert_{M_{k-1}}^2 - 2\alpha_k\cdot \mb{E}_k\left[\nu_{k+1}^{-2}\left(\phi\left(\tilde{x}_{k+1}\right) - \phi\left(x\right)\right)\right]\right)\\
            &\hspace{12em}+\ \left(\rho_{f,m,U}L_F(U)^2 + \gamma^2\right)\cdot\left(1+\alpha_k^2\right)\alpha_k^2
    \end{split}
\end{equation}
holds for all $k\in \mb{Z}_{+}$ and $x\in \mb{R}^d$. By applying (\ref{equ:sippinequhatx1tvalueforproofzhutight}) and following an argument nearly identical to that in the proof of \cite[Theorem~2.2]{zhu2025tight}, we can conclude that on the event that  (\ref{equ:boundxknsumsqsumalphak}) holds, 
\begin{equation}\label{equ:infphixkphistarnukp1}
    \liminf_{k\to \infty} \nu_{k+1}^{-2}\left(\phi\left(\tilde{x}_{k+1}\right) - \phi^{*}\right) = 0\quad \mbox{with probability}\ 1.
\end{equation}
Since $\inf_{k\in \mb{Z}_{+}}\nu_{k+1} \geq 1$, the estimate (\ref{equ:infphixkphistarnukp1}) implies 
\begin{equation}\label{equ:infphixkphistar}
    \liminf_{k\to \infty} \phi\left(\tilde{x}_{k+1}\right) = \phi^{*}\quad \mbox{with probability}\ 1.
\end{equation}

\underline{\textit{*Step 2.}} Denote by $C_{f,m,U,\gamma} \triangleq \rho_{f,m,U}L_F(U)^2 + \gamma^2$. Due to (\ref{equ:infphixkphistar}), it suffices to prove that on the event that (\ref{equ:boundxknsumsqsumalphak}) and (\ref{equ:infphixkphistar}) hold, the sequence $\{x_k\}$ converges to some $x^{*}\in \mc{X}^{*}$ almost surely. On the event that (\ref{equ:boundxknsumsqsumalphak}) and (\ref{equ:infphixkphistar}) hold, under Assumption \ref{asm:sppaall} and \ref{asm:sppbdoptset}, we have $\phi^{*}\in \mb{R}$ and thus $\liminf_{k\to \infty} \phi(\tilde{x}_k) = \phi^{*} \in \mb{R}$. Using the definition of the limit inferior, one can easily deduce that there exists a convergent subsequence $\{\phi(\tilde{x}_k)\}_{k\in \mc{K}}$ of $\{\phi(\tilde{x}_k)\}$ such that 
\begin{equation}\label{equ:existsubsequencevalue}
    \lim_{k\to \infty, k\in \mc{K}} \phi\left(\tilde{x}_k\right) = \liminf_{k\to \infty} \phi\left(\tilde{x}_k\right) = \phi^{*}.
\end{equation}
Fix any $\overline{x}\in \mc{X}^{*}$. Since $\inf_{k\in \mb{Z}_{+}}\nu_k \geq 1$ and $\phi(\olx) = \phi^{*}$, we can deduce from (\ref{equ:sippinequhatx1tvalue}) with $x = \overline{x}$ that 
\begin{align*}
    &\mb{E}_k\left[\nu_{k+1}^{-2}\Norm{ x_{k+1} - \overline{x}}_{M_k}^2\right]\
        \leq \left(1+\alpha_k^2\right)\cdot \nu_k^{-2}\Norm{ x_k - \overline{x}}_{M_{k-1}}^2 + C_{f,m,U,\gamma}\cdot\alpha_k^2\left(1+\alpha_k^2\right)
\end{align*}
holds for all $k\in \mb{Z}_{+}$, which combined with Lemma \ref{lema:supermartingaleconvergence}, $\sum_{k=1}^\infty \alpha_k^2 < \infty$ and $\nu_{k+1}\geq 1$ confirms that 
the sequence $\{\nu_k^{-1}\Vert x_k - \overline{x}\Vert_{M_{k-1}}\}$ converges to some finite value with probability $1$. 
By Assumption \ref{asm:mk}, we have $\lim_{k\to \infty}\nu_k = \nu_\infty$ and thus for each $\overline{x}\in \mc{X}^{*}$, 
\begin{align*}
    \mbox{the sequence}\ \{\Norm{x_k - \overline{x}}_{M_{k-1}}\}\ \mbox{converges with probability}\ 1.
\end{align*}
Set $\tlx_1 \triangleq x_1$. Combining the previous result with the fact that 
\begin{equation}\label{equ:asumalmostconvergence}
    \begin{split}
        \lim_{k\to \infty} \alpha_k = 0\quad \mbox{and}\quad \Vert x_{k+1} - \tilde{x}_{k+1}\Vert_{M_k} \leq \epsilon_{k} = \gamma \alpha_k^2\quad \mbox{for all}\ k\in \mb{Z}_{+}
    \end{split}
\end{equation}
yields
\begin{equation}\label{equ:sippconvergedistanceas1}
    \mbox{the sequence}\ \{\Norm{ \tilde{x}_k - \overline{x}}_{M_{k-1}}\}\ \mbox{converges with probability}\ 1.
\end{equation}
Note that Assumption \ref{asm:mk} implies that $\Vert \cdot\Vert_{M_{k-1}}\leq \mu_\infty^{-1} \Vert \cdot\Vert_2$ holds for all $k\in \mb{Z}_{+}$. 
Then, on the event that (\ref{equ:boundxknsumsqsumalphak}) and (\ref{equ:infphixkphistar}) hold, by following an argument nearly identical to that in the proof of \cite[Theorem~2.2]{zhu2025tight}, we can deduce from (\ref{equ:existsubsequencevalue}) together with (\ref{equ:sippconvergedistanceas1}) that there exists some $x^{*}\in \mc{X}^{*}$ such that 
\begin{align*}
    \lim_{k\to \infty}\Vert \tlx_k - x^{*}\Vert_{M_{k-1}} = 0\ \mbox{with probability}\ 1,
\end{align*} 
which combined with  (\ref{equ:asumalmostconvergence}) yields 
\begin{equation}\label{equ:alconvxkmknorm}
    \lim_{k\to \infty}\Vert x_k - x^{*}\Vert_{M_{k-1}} = 0\ \mbox{with probability}\ 1.
\end{equation}
Since Assumption \ref{asm:mk} implies $\Vert \cdot \Vert_2 \leq L_\infty \Vert \cdot\Vert_{M_{k-1}}$ holds for all $k\in \mb{Z}_{+}$, it then follows from (\ref{equ:alconvxkmknorm}) that $\lim_{k\to \infty}\Vert x_k - x^{*}\Vert_2 = 0$ holds with probability $1$. 
This demonstrates the almost sure convergence of the iterate sequence $\{x_k\}$, and thus completes the proof of Theorem \ref{thm:sippstabilityconvex}. 
\end{proof}

In the convex setting, Theorem \ref{thm:sippstabilityconvex} extends the convergence guarantee derived in \cite[Theorem~2.1~and~2.2]{zhu2025tight} for the isPPA to its preconditioned variants. It is worth noting that, for the convex case, the stability of the preconditioned isPPA does not require any additional coercivity assumptions (e.g., Assumption \ref{asm:sublevelset}), mirroring the convergence properties of its deterministic counterpart. 

\subsection{Convergence Rate for the Weakly Convex Case}\label{subsec:sipprateingradme}
Based on the established almost sure boundedness of the iterates generated by Algorithm \ref{algo:ispppa}, we now proceed to derive convergence rates of the model-based ispPPA under a local Lipschitz condition on $f_x(\cdot;s)$ in the weakly convex setting. Unless stated otherwise, for the set $V$ specified in Assumption \ref{asm:stomodel}, if there exists a bounded open convex subset of $V$ that contains both $\lrbrackets{x_k}$ and $\lrbrackets{\tilde{x}_k}$, we denote this subset by $U$. Recall the function $\varsigma_\beta$ defined in (\ref{equ:varphi}) and the constant $\rho_{f,m,U}$ introduced in (\ref{equ:constantminibatch}). The main convergence results are presented in Theorem \ref{thm:sippmdiminishing} and \ref{thm:sippmdiminishing2}. Compared with the results in \cite{davis2019stochastic,gao2024stochastic}, the convergence analysis developed herein allows each subproblem to be solved inexactly and extends to the preconditioned variants of the stochastic model-based algorithms. It is worth noting that, as demonstrated in Remark \ref{remark:alphaklbubwcvx}, the requirement that the stepsizes $\{\alpha_k\}$ satisfy $\alphaol k^{-\beta} \leq \alpha_k \leq \alphaou k^{-\beta}$ for some positive scalars $\alphaol,\alphaou$ is practically reasonable. 
Since $\nabla e_{(\olrho)^{-1}\phi}^{M_{k-1}}(x_k) = \olrho M_{k-1}(x_k - \prox_{(\olrho)^{-1}\phi}^{M_{k-1}}(x_k))$, we also refer to $\Vert \nabla e_{(\olrho)^{-1}\phi}^{M_{k-1}}(x_k)\Vert_2$ as the Karush–Kuhn–Tucker (KKT) residual. 
\begin{theorem}\label{thm:sippmdiminishing}
    Let Assumption \ref{asm:stomodel}-\ref{asm:sppa1} hold, and assume that $\phi^{*}>-\infty$. Let $\{x_k\}$ denote the iterates generated by the model-based ispPPA (Algorithm \ref{algo:ispppa}) with stepsizes $\{\alpha_k\}$ and parameters $\{\epsilon_k\}$. Denote by $\eta \triangleq L_\infty^2\cdot\oleta$ and $\tau \triangleq L_\infty^2\cdot\oltau$. Set $B_\rho \triangleq \sup_{k\in \mb{Z}_{+}}\rho_k^2$. Suppose that there exists $\gamma\in \mb{R}_{+}$ and $\olrho > B_\rho \cdot(\eta + \tau)$ such that 
    \begin{equation}\label{asm:epsilonkalphakdiminishingwcvx}
        \begin{split}
            &\alphaol k^{-\beta}\leq \alpha_k \leq \alphaou k^{-\beta}\quad \mbox{and}\quad \epsilon_k = \gamma \alpha_k^2\quad \mbox{for all}\ k\in \mb{Z}_{+},\\
            \mbox{where}\quad &0<\alphaol\leq\alphaou < \left(\olrho\max\lrbrackets{2,B_\rho}\right)^{-1}\quad \mbox{and}\quad \beta\in \left(\frac{1}{2},1\right).
        \end{split}
    \end{equation}
    Then, on the event that $\sup_{k\in \mb{Z}_{+}}\Norm{x_k}_2 < \infty$, we have 
    \begin{equation}\label{equ:bdminnorm}
        \frac{1}{k}\sum_{i=1}^k \mb{E}\left[\Norm{\nabla \mephi{i}\left(x_i\right)}_2^2\right] \leq \frac{\nu_\infty^2}{\theta_1\mu_\infty^2} \left(\frac{B_\phi}{\alphaol}k^{-\left(1-\beta\right)} +  \frac{2\theta_2 \alphaou}{1-\beta}k^{-\beta} - \frac{2^\beta\theta_2\alphaou}{1-\beta}k^{-1}\right)  
    \end{equation}
    holds for all $k\in \mb{Z}_{+}$, where $\mephi{i}(\cdot) \triangleq e^{M_{i-1}}_{(\olrho)^{-1}\phi}(\cdot)$,
    \begin{equation}\label{def:theta1theta2}
        \begin{split}
            \theta_1 &= \frac{\left(\olrho B_\rho^{-1}-\eta-\tau\right)\left(1-2\olrho\alphaou\right)^2}{2\olrho},\\
            \theta_2 &= \frac{\olrho\left[1+4\left(\olrho - \eta\right)\left(\olrho-\eta-\tau\right)\alphaou^2\right]}{2}\\
            &\hspace{2em}\cdot \left[\frac{\rho_{f,m,U}L_F\left(U\right)^2}{\left(1-\olrho \alphaou\right)\left(1-\olrho B_\rho\alphaou\right)} + \frac{\gamma^2}{4\left(\olrho B_\rho^{-1} - \eta\right)\left(\olrho B_\rho^{-1}-\eta-\tau\right)}\right]\quad \mbox{with}\\
            B_\phi &= e_{\left(\olrho\right)^{-1}\phi}(x_1) - \phi^{*} + \frac{4^\beta\theta_2\alphaou^2}{2\beta-1}\quad \mbox{and}\quad L_F\left(U\right)\triangleq L_\infty\cdot\olL_F\left(U\right).
        \end{split}
    \end{equation}
\end{theorem}
\begin{theorem}\label{thm:sippmdiminishing2}
    Let Assumption \ref{asm:stomodel}-\ref{asm:sppa1} hold and assume that $\phi^{*}>-\infty$. Let $\{x_k\}$ denote the iterates generated by the model-based ispPPA (Algorithm \ref{algo:ispppa}) with stepsizes $\{\alpha_k\}$ and parameters $\{\epsilon_k\}$. Denote by $\eta \triangleq L_\infty^2\cdot\oleta$ and $\tau \triangleq L_\infty^2\cdot\oltau$. Set $B_\rho \triangleq \sup_{k\in \mb{Z}_{+}}\rho_k^2$. Suppose that there exists $\gamma\in \mb{R}_{+}$ and $\olrho > B_\rho \cdot(\eta + \tau)$ such that 
    \begin{equation}\label{asm:epsilonkalphakdiminishingwcvx2}
        \begin{split}
            &\alphaol k^{-\beta}\leq \alpha_k \leq \alphaou k^{-\beta}\quad \mbox{and}\quad \epsilon_k = \gamma \alpha_k^2\quad \mbox{for all}\ k\in \mb{Z}_{+},\\
            \mbox{where}\quad &0<\alphaol\leq \alphaou < \left(\olrho\max\lrbrackets{2,B_\rho}\right)^{-1}\quad \mbox{and}\quad \beta\in \left(0,1\right].
        \end{split}
    \end{equation}
    Assume further that the output variable $x = x_{i_{*}}$ is sampled from the iterates $\{x_i\}_{i\in [k]}$ with probability $\mb{P}\{i_{*} = i\} = \frac{\alpha_i}{\sum_{j=1}^k \alpha_j}$. Then, on the event that $\sup_{k\in \mb{Z}_{+}}\Norm{x_k}_2 < \infty$, we have
    \begin{equation}\label{equ:bdminnormstar}
        \mb{E}\left[\Norm{\nabla \mephi{i_{*}}\left(x_{i_{*}}\right)}_2^2\right] \leq \frac{\nu_\infty^2} {\theta_1\mu_\infty^2} \cdot\frac{e_{\left(\olrho\right)^{-1}\phi}(x_1)-\phi^{*} + \theta_2\sum_{i=1}^k \alpha_i^2}{\sum_{i=1}^k \alpha_i}  
    \end{equation}
    holds for all $k\in \mb{Z}_{+}$, where $\mephi{i_{*}}(\cdot) \triangleq e^{M_{i_{*}-1}}_{(\olrho)^{-1}\phi}(\cdot)$, and $\theta_1$ and $\theta_2$ are positive scalars defined in (\ref{def:theta1theta2}). 
\end{theorem}
Focusing on the dominant terms in Theorem \ref{thm:sippmdiminishing2}, we can immediately obtain the asymptotic convergence rates of the model-based ispPPA, summarized as follows. A detailed proof of Corollary \ref{coro:sipprate} can be found in Appendix \ref{apx:proofcorosipprate}.
\begin{coro}\label{coro:sipprate}
    Consider the setting of Theorem \ref{thm:sippmdiminishing2}. On the event that $\sup_{k\in \mb{Z}_{+}}\Norm{x_k}_2 < \infty$, the following assertion holds:
    \begin{equation}\label{equ:asymptoticrate}
        \begin{split}
            \mb{E}\left[\Norm{\nabla \mephi{i_{*}}\left(x_{i_{*}}\right)}_2^2\right] &\leq \frac{\nu_\infty^2} {\theta_1\mu_\infty^2} \cdot\frac{e_{\left(\olrho\right)^{-1}\phi}(x_1)-\phi^{*} + 4^\beta\theta_2\alphaou^2\cdot\varsigma_{1-2\beta}\left(k+1\right)}{\alphaol\cdot \varsigma_{1-\beta}\left(k+1\right)}\\
            &\leq \begin{cases}
                \mc{O}\left(k^{-\beta}\right)\quad &\mbox{if}\ \beta\in \left(0,\frac{1}{2}\right),\\
                \mc{O}\left(k^{-\frac{1}{2}}\ln k\right)\quad &\mbox{if}\ \beta = \frac{1}{2},\\
                \mc{O}\left(k^{-\left(1-\beta\right)}\right)\quad &\mbox{if}\ \beta\in \left(\frac{1}{2},1\right),\\
                \mc{O}\left(\left(\ln k\right)^{-1}\right)\quad &\mbox{if}\ \beta=1.
            \end{cases}  
        \end{split}  
    \end{equation}
\end{coro}

\paragraph*{Proof of Theorem \ref{thm:sippmdiminishing} and \ref{thm:sippmdiminishing2}.} 
To establish Theorem \ref{thm:sippmdiminishing} and \ref{thm:sippmdiminishing2}, we present the following recursive relation, with its proof given in Appendix \ref{apx:proofclaimmephiinequ1further}.
\begin{claim}\label{claim:mephiinequ1further}
    Consider the setting of Theorem \ref{thm:sippmdiminishing} or \ref{thm:sippmdiminishing2}. Denote by $\olx_k \triangleq \prox_{(\olrho)^{-1}\phi}^{M_{k-1}}(x_k)$, $\mephi{k}(\cdot) \triangleq e_{(\olrho)^{-1}\phi}^{M_{k-1}}(\cdot)$ and $\tlx_{k+1}\triangleq \mr{prox}_{\alpha_k \olvphi_{x_k}(\cdot;S_k^{1:m})}^{M_k}(x_k)$ for each $k\in \mb{Z}_{+}$. Then, both $\olx_k$ and $\tlx_{k+1}$ are well-defined, and $\mephi{k}(x_k) - \phi^{*} \geq 0$ for all $k\in \mb{Z}_{+}$. Moreover, on the event that $\sup_{k\in \mb{Z}_{+}}\Norm{x_k}_2 < \infty$, there exists a bounded open convex subset of $V$, denoted by $U$, such that both $\{\tlx_k\}$ and $\{x_k\}$ are contained in $U$, and the following assertion holds for all $k\in \mb{Z}_{+}$:
    \begin{equation}\label{equ:medecrease}
        \begin{split}
            &\mb{E}_k\left[\nu_{k+1}^{-2}\cdot\left(\mephi{k+1}\left(x_{k+1}\right) - \phi^{*}\right)\right]\\
            \leq\ &\nu_k^{-2}\cdot\left(\mephi{k}\left(x_k\right) - \phi^{*}\right) - \theta_1 \alpha_k \cdot \nu_k^{-2}\Norm{\olrho \left(x_k - \olx_k\right)}_{M_{k-1}}^2 + \theta_2 \alpha_k^2
        \end{split}
    \end{equation}
    where $\theta_1$ and $\theta_2$ are positive scalars defined in (\ref{def:theta1theta2}).
\end{claim}
Use the same notations as in Claim \ref{claim:mephiinequ1further}. 
\begin{proof}[Proof of Theorem \ref{thm:sippmdiminishing}]  
Define $Y_k \triangleq \nu_k^{-2}\cdot(\mephi{k}(x_k) - \phi^{*}) + \theta_2\sum_{j=k}^\infty \alpha_j^2$ for each $k\in \mb{Z}_{+}$. Note that $\theta_1>0$ by (\ref{def:theta1theta2}). Fix any $k\in \mb{Z}_{+}$. Using the definition of $\{Y_k\}$ and appealing to the estimate (\ref{equ:medecrease}), we have 
\begin{equation}\label{equ:supermartingaleyk}
    \begin{split}
        &\mb{E}_k \left[Y_{k+1}\right] 
        \leq \nu_k^{-2}\cdot\left(\mephi{k}\left(x_k\right) - \phi^{*}\right) + \theta_2\sum_{j=k}^\infty \alpha_j^2 = Y_k.
    \end{split}
\end{equation}
Since condition (\ref{asm:epsilonkalphakdiminishingwcvx}) implies $\alpha_k \leq \alphaou k^{-\beta}$ and $\beta\in (\frac{1}{2},1)$, we have $\sum_{j=1}^\infty \alpha_j^2 \leq \frac{4^\beta \alphaou^2}{2\beta - 1}$\footnote{Since $\beta > \frac{1}{2}$, it follows from Lemma \ref{lema:inequtech} in Appendix \ref{secapx:proofofsipprategradme} that $\sum_{j=1}^k j^{-2\beta} \leq \frac{4^\beta}{2\beta - 1}$ holds for all $k\in \mb{Z}_{+}$.}. 
Recall that $\nu_1 = 1$, $M_0 = I_d$ and $B_\phi = e_{(\olrho)^{-1}\phi}(x_1) - \phi^{*} + \frac{4^\beta\theta_2\alphaou^2}{2\beta-1}$. Then, using  the estimate $\mephi{k}\left(x_k\right) -\phi^{*}\geq 0$ by Claim \ref{claim:mephiinequ1further}, we deduce from (\ref{equ:supermartingaleyk}) that 
\begin{equation}\label{equ:boundf}
    \begin{split}
        0\leq \mb{E}\left[Y_k\right] &\leq Y_1 = e_{\left(\olrho\right)^{-1}\phi}(x_1) - \phi^{*} + \theta_2\sum_{j=1}^\infty \alpha_j^2 \leq B_\phi < \infty\quad \mbox{for all}\ k\in \mb{Z}_{+}.
    \end{split}
\end{equation}
Taking expectation on both sides of (\ref{equ:medecrease}) and summing over $i = 1,\cdots,k$ yields 
\begin{align*}
    &\theta_1 \sum_{i=1}^k \nu_i^{-2} \mb{E}\left[\Norm{\olrho \left(x_i - \olx_i\right)}_{M_{i-1}}^2\right]\\
        \leq\ 
        &\frac{1}{\alpha_1}\nu_1^{-2}\cdot\left(\mephi{1}\left(x_1\right) - \phi^{*}\right) + \sum_{i=2}^k\left(\frac{1}{\alpha_i} - \frac{1}{\alpha_{i-1}}\right)\mb{E}\left[\nu_i^{-2}\cdot\left(\mephi{i}\left(x_i\right) - \phi^{*}\right)\right]\\
        &\hspace{0em} - \frac{1}{\alpha_k} \mb{E}\left[\nu_{k+1}^{-2}\cdot\left(\mephi{k+1}\left(x_{k+1}\right) - \phi^{*}\right)\right] + \theta_2 \sum_{i=1}^k \alpha_i\\
        \leq\ &\frac{B_\phi}{\alpha_1} + B_\phi \sum_{i=2}^k\left(\frac{1}{\alpha_i} - \frac{1}{\alpha_{i-1}}\right) + \theta_2 \sum_{i=1}^k \alpha_i = \frac{B_\phi}{\alpha_k} + \theta_2 \sum_{i=1}^k \alpha_i\\
        \leq\ &
        \frac{B_\phi}{\alphaol}k^\beta + \theta_2 \alphaou 2^\beta \frac{\left(k+1\right)^{1-\beta}-1}{1-\beta}\\
        \leq\ &
        \frac{B_\phi}{\alphaol}k^\beta +  \frac{2\theta_2 \alphaou}{1-\beta}k^{1-\beta} - \frac{2^\beta\theta_2\alphaou}{1-\beta},
\end{align*}
where the second inequality uses (\ref{equ:boundf}) together with the definition of $Y_k$, the third holds because $\alphaol k^{-\beta} \leq \alpha_k \leq \alphaou k^{-\beta}$ and $\sum_{i=1}^k i^{-\beta} \leq 2^\beta \frac{\left(k+1\right)^{1-\beta}-1}{1-\beta}$\footnote{Since $\frac{1}{2} < \beta < 1$, it follows from Lemma \ref{lema:inequtech} in Appendix \ref{secapx:proofofsipprategradme} that $\sum_{j=1}^k j^{-\beta} \leq 2^\beta\frac{\left(k+1\right)^{1-\beta}-1}{1-\beta}$ holds for all $k\in \mb{Z}_{+}$.}, and the last follows from $k+1\leq 2k$ as $k\geq 1$. Combining the preceding inequality with (\ref{equ:gradnormmkmeolxk}), $\sup_{k\in \mb{Z}_{+}}\nu_k \leq \nu_\infty$ and $\theta_1 > 0$, we deduce 
\begin{equation}\label{equ:upratemk}
    \begin{split}
        &\frac{1}{k}\sum_{i=1}^k \mb{E}\left[\Norm{\nabla \mephi{i}\left(x_i\right)}_{M_{i-1}^{-1}}^2\right] = \frac{1}{k}\sum_{i=1}^k \mb{E}\left[\Norm{\olrho \left(x_i - \olx_i\right)}_{M_{i-1}}^2\right]\\
        \leq\ &\frac{1}{k}\sum_{i=1}^k \nu_\infty^2\cdot\nu_i^{-2}\mb{E}\left[\Norm{\olrho \left(x_i - \olx_i\right)}_{M_{i-1}}^2\right]\\
        \leq\ &\frac{\nu_\infty^2}{\theta_1} \left(\frac{B_\phi}{\alphaol}k^{-\left(1-\beta\right)} +  \frac{2\theta_2 \alphaou}{1-\beta}k^{-\beta} - \frac{2^\beta\theta_2\alphaou}{1-\beta}k^{-1}\right).
    \end{split}
\end{equation}
Using the estimate (\ref{equ:normmkm1}) and condition $\inf_{k\in \mb{Z}_{+}} \mu_k \geq \mu_\infty$ by Assumption \ref{asm:mk}, we have 
\begin{align*}
    \Norm{\nabla \mephi{i}\left(x_i\right)}_2 \leq \mu_{i-1}^{-1}\cdot\Norm{\nabla \mephi{i}\left(x_i\right)}_{M_{i-1}^{-1}} \leq \mu_\infty^{-1}\cdot\Norm{\nabla \mephi{i}\left(x_i\right)}_{M_{i-1}^{-1}}
\end{align*}
and thus 
\begin{equation}\label{equ:mknormreg}
    \mb{E}\left[\Norm{\nabla \mephi{i}\left(x_i\right)}_2 ^2\right] \leq \mu_\infty^{-2}\cdot \mb{E}\left[\Norm{\nabla \mephi{i}\left(x_i\right)}_{M_{i-1}^{-1}}^2\right]\quad \mbox{for all}\ i\in [k]. 
\end{equation}
Combining (\ref{equ:upratemk}) together with (\ref{equ:mknormreg}) implies that the estimate (\ref{equ:bdminnorm}) holds. This completes the proof of Theorem \ref{thm:sippmdiminishing}.
\end{proof}

\begin{proof}[Proof of Theorem \ref{thm:sippmdiminishing2}] 
Recall that the estimate (\ref{equ:medecrease}) from Claim \ref{claim:mephiinequ1further} holds, i.e., 
\begin{equation}\label{equ:medecreaserate2}
    \begin{split}
        &\theta_1 \alpha_k \cdot \nu_k^{-2}\Norm{\olrho \left(x_k - \olx_k\right)}_{M_{k-1}}^2\\
        \leq\ &\nu_k^{-2}\cdot\left(\mephi{k}\left(x_k\right) - \phi^{*}\right) - \mb{E}_k\left[\nu_{i+1}^{-2}\cdot\left(\mephi{k+1}\left(x_{k+1}\right) - \phi^{*}\right)\right] + \theta_2 \alpha_k^2
    \end{split}
\end{equation}
holds for all $k\in \mb{Z}_{+}$, where $\theta_1$ and $\theta_2$ are positive scalars defined in (\ref{def:theta1theta2}). 
Fix any $k\in \mb{Z}_{+}$. Taking expectation on both sides of (\ref{equ:medecreaserate2}) and summing over $i=1,\cdots,k$ yields 
\begin{align}
        &\theta_1\sum_{i=1}^k \alpha_i\mb{E}\left[\nu_i^{-2}\Norm{\olrho \left(x_i - \olx_i\right)}_{M_{i-1}}^2\right]\nonumber\\
        \leq\ 
        &\sum_{i=1}^k \lrbrackets{\mb{E}\left[\nu_i^{-2}\cdot\left(\mephi{i}\left(x_i\right) - \phi^{*}\right)\right] - \mb{E}\left[\nu_{i+1}^{-2}\cdot\left(\mephi{i+1}\left(x_{i+1}\right) - \phi^{*}\right)\right]} + \theta_2\sum_{i=1}^k \alpha_i^2\nonumber\\
        =\ &e_{\left(\olrho\right)^{-1}\phi}\left(x_1\right)-\phi^{*} - \mb{E}\left[\nu_{k+1}^{-2}\cdot\left(\mephi{k+1}\left(x_{k+1}\right) - \phi^{*}\right)\right] + \theta_2\sum_{i=1}^k \alpha_i^2\nonumber\\
        \leq\ &e_{\left(\olrho\right)^{-1}\phi}\left(x_1\right)-\phi^{*} + \theta_2\sum_{i=1}^k \alpha_i^2,\label{equ:medecreaserate2a}
\end{align}
where the equality is due to $\nu_1 = \rho_0 = 1$ and $M_0 = I_d$ by Assumption \ref{asm:mk}, and the last inequality holds because $\mephi{k+1}\left(x_{k+1}\right) - \phi^{*}\geq 0$ by Claim \ref{claim:mephiinequ1further}. Combining the estimate (\ref{equ:medecreaserate2a}) with (\ref{equ:mknormreg}), (\ref{equ:gradnormmkmeolxk}) and $\sup_{k\in \mb{Z}_{+}}\nu_k \leq \nu_\infty$ suggests that 
\begin{equation}\label{equ:medecreaserate2b}
    \begin{split}
        &\sum_{i=1}^k \alpha_i\mb{E}\left[\Norm{\nabla \mephi{i}\left(x_i\right)}_2^2\right] \leq \mu_\infty^{-2}\cdot\sum_{i=1}^k \alpha_i\mb{E}\left[\Norm{\nabla \mephi{i}\left(x_i\right)}_{M_{i-1}^{-1}}^2\right]\\
        \leq\ &\frac{\nu_\infty^2}{\mu_\infty^2}\cdot \sum_{i=1}^k \alpha_i\mb{E}\left[\nu_2^{-2}\Norm{\olrho\left(x_i - \olx_i\right)}_{M_{i-1}}^2\right] 
        \leq \frac{\nu_\infty^2}{\theta_1\mu_\infty^2}\left(e_{\left(\olrho\right)^{-1}\phi}\left(x_1\right)-\phi^{*} + \theta_2\sum_{i=1}^k \alpha_i^2\right).
    \end{split}
\end{equation}
Since the output variable $x = x_{i_{*}}$ is sampled from the iterates $\{x_i\}_{i\in [k]}$ with probability $\mb{P}\{i_{*} = i\} = \frac{\alpha_i}{\sum_{j=1}^k \alpha_j}$, we have 
\begin{equation}\label{equ:medecreaserate2c}
    \mb{E}\left[\Norm{\nabla \mephi{i_{*}}\left(x_{i_{*}}\right)}_2^2\right] = \left(\sum_{i=1}^k \alpha_i\right)^{-1}\cdot \sum_{i=1}^k \alpha_i\mb{E}\left[\Norm{\nabla \mephi{i}\left(x_i\right)}_2^2\right].
\end{equation}
The estimate (\ref{equ:bdminnormstar}) follows immediately from the combination of (\ref{equ:medecreaserate2b}) and (\ref{equ:medecreaserate2c}). This completes the proof of Theorem \ref{thm:sippmdiminishing2}.

\end{proof}

\subsection{Convergence Rate Under the Quadratic Growth Condition}\label{subsec:sipprateinsqdistopt}
In this section, we derive nonasymptotic and asymptotic convergence rates for the model-based ispPPA under an additional quadratic growth condition on the objective function. Let Assumption \ref{asm:stomodel}-\ref{asm:sppa2} hold. Denote by $c_1\triangleq \mu_\infty^2\cdot\overline{c}_1$ and assume that $c_1 > \frac{\tau+\eta}{2}$. Let $s\in (0,2)$ be an arbitrary scalar and define 
\begin{equation}\label{equ:sipprateboundconstant1}
    \begin{split}
        &\tilde{C}_{f,m,U,\tau,\eta,c_1,\gamma}(\alpha)\\ \triangleq\ &\rho_{f,m,U} L_F\left(U\right)^2 + \frac{\left(1+\left(2c_1-\tau\right)\alpha\right)\left(2\left(1+ \eta\alpha\right) + s\left(2c_1-\tau-\eta\right)\alpha\right)}{(2-s)\left(2c_1-\tau-\eta\right)\left(1+\eta\alpha\right)}\gamma^2\\
    \end{split}
\end{equation}
for all $\alpha\in \mb{R}_{+}$, where $\tau = L_\infty^2\cdot \oltau$, $\eta = L_\infty^2\cdot\oleta$ and $L_F(U)\triangleq L_\infty\cdot \olL_F(U)$. By using the stability and imposing additional constraints on the constants involved in Assumption \ref{asm:stomodel}-\ref{asm:sppa2}, we establish nonasymptotic convergence retes measured by the distance to the optimal solution set. The main results are summarized in Theorem \ref{thm:sippmconstantquadratic} and \ref{thm:sippmdiminishingquadratic}. 
\begin{theorem}\label{thm:sippmconstantquadratic}
    Let Assumption \ref{asm:stomodel}-\ref{asm:sppa2} hold, and assume that $\mc{X}^{*}$ is nonempty. Let $\{x_k\}$ denote the iterates generated by the model-based ispPPA (Algorithm \ref{algo:ispppa}) with stepsizes $\{\alpha_k\}$ and parameters $\{\epsilon_k\}$. Denote by $\eta \triangleq L_\infty^2\cdot\oleta$ and $\tau \triangleq L_\infty^2\cdot\oltau$. Suppose that there exists $\gamma\in \mb{R}_{+}$ and $\olrho > \eta + \tau$ such that 
    \begin{align*}
        \alpha_k = \alpha_0\quad \mbox{and}\quad \epsilon_k = \gamma \alpha_0^{\frac{3}{2}}\quad \mbox{for all}\ k\in \mb{Z}_{+},\quad \mbox{where}\quad \alpha_0\in \left(0,\left(\olrho\right)^{-1}\right).
    \end{align*}
    Set $c_{\tau,\eta} \triangleq c_1 - \frac{\tau + \eta}{2}$, $\tilde{c}_1 \triangleq \frac{\olrho c_{\tau,\eta}}{\olrho + \eta}$ and $\delta_1\triangleq \mr{dist}(x_1,\mc{X}^{*})^2$. Let $s\in (0,2)$ be an arbitrary scalar. If $c_1>\frac{\tau+\eta}{2}$, then on the event that $\sup_{k\in \mb{Z}_{+}}\Norm{x_k}_2 < \infty$, we have 
    \begin{equation}\label{equ:sippconstantbound1}
        \begin{split}
            \mb{E}\left[\mr{dist}(x_k,\mc{X}^{*})^2\right] &\leq \left(\nu_\infty L_\infty\right)^2\cdot\left[\left(\frac{1}{1+s\tilde{c}_1\alpha_0}\right)^{k-1} \delta_1 + \frac{\tilde{C}_{f,m,U,\tau,\eta,c_1,\gamma}(\alpha_0)}{s\tilde{c}_1}\alpha_0\right]
        \end{split}
    \end{equation}
    for all $k\in \mb{Z}_{+}$, where $\tilde{C}_{f,m,U,\tau,\eta,c_1,\gamma}(\cdot)$ is defined in (\ref{equ:sipprateboundconstant1}).  
\end{theorem}
\begin{theorem}\label{thm:sippmdiminishingquadratic}
    Let Assumption \ref{asm:stomodel}-\ref{asm:sppa2} hold, and assume that $\mc{X}^{*}$ is nonempty. Let $\{x_k\}$ denote the iterates generated by the model-based ispPPA (Algorithm \ref{algo:ispppa}) with stepsizes $\{\alpha_k\}$ and parameters $\{\epsilon_k\}$. Denote by $\eta \triangleq L_\infty^2\cdot\oleta$ and $\tau \triangleq L_\infty^2\cdot\oltau$. Suppose that there exists $\gamma\in \mb{R}_{+}$ and $\olrho > \eta + \tau$ such that 
    \begin{equation}\label{asm:epsilonkalphakdiminishingquadratic}
        \begin{split}
            &\alphaol k^{-\beta} \leq \alpha_k \leq \alphaou k^{-\beta}\quad \mbox{and}\quad \epsilon_k = \gamma \alpha_k^{\frac{3}{2}}\quad \mbox{for all}\ k\in \mb{Z}_{+},\\ \mbox{where}\quad &\alphaol,\alphaou\in \left(0,\left(\olrho\right)^{-1}\right)\quad \mbox{and}\quad \beta\in (0,1].
        \end{split}
    \end{equation}
    Set $c_{\tau,\eta} \triangleq c_1 - \frac{\tau + \eta}{2}$, $\tilde{c}_1 \triangleq \frac{\olrho c_{\tau,\eta}}{\olrho + \eta}$ and $\delta_1\triangleq \mr{dist}(x_1,\mc{X}^{*})^2$. Let $s\in (0,2)$ be an arbitrary scalar. If $c_1>\frac{\tau+\eta}{2}$, then on the event that $\sup_{k\in \mb{Z}_{+}}\Norm{x_k}_2 < \infty$, the following assertions hold:
    \begin{enumerate}
        \item If $\beta\in (0,1)$, then 
        \begin{equation}\label{equ:sippdiminishingbound1}
            \begin{split}
                &\mb{E}\left[\mr{dist}(x_k,\mc{X}^{*})^2\right]
                \leq \left(\nu_\infty L_\infty\right)^2\cdot\left[\exp\left(-\frac{s\tilde{c}_1\alphaol}{1+s\tilde{c}_1\alphaou}\varsigma_{1-\beta}(k)\right)\delta_1\right.\\ &\hspace{0em}+\ \tilde{C}_{f,m,U,\tau,\eta,c_1,\gamma}(\alphaou)\cdot4^{\beta}\frac{\alphaou^2}{1-\mr{exp}(-\frac{s\tilde{c}_1\alphaou}{1+s\tilde{c}_1\alphaou})}\cdot \frac{1}{k^\beta}\\
                &\left.\hspace{0em}+\ \tilde{C}_{f,m,U,\tau,\eta,c_1,\gamma}(\alphaou)\cdot4^{\beta}\alphaou^2\cdot\mr{exp}\left(-\frac{s\tilde{c}_1\alphaol}{2\left(1+s\tilde{c}_1\alphaou\right)}k^{1-\beta}\right)\cdot \varsigma_{1-2\beta}\left(\frac{k}{2}\right)\right]
            \end{split}
        \end{equation}
        holds for all $k\in \mb{Z}_{+}$, where $\tilde{C}_{f,m,U,\tau,\eta,c_1,\gamma}(\cdot)$ is defined in (\ref{equ:sipprateboundconstant1}).
        \item If $\beta = 1$, then 
        \begin{equation}\label{equ:sippdiminishingbound3}
            \begin{split}
                &\mb{E}\left[\mr{dist}(x_k,\mc{X}^{*})^2\right] \leq \left(\nu_\infty L_\infty\right)^2\cdot\\ &\hspace{2em}\begin{cases}
                    &\left[\frac{4\delta_1 }{1+s\tilde{c}_1\alphaol}+ \frac{4\tilde{C}_{f,m,U,\tau,\eta,c_1,\gamma}(\alphaou)(2+s\tilde{c}_1\alphaol)\alphaou^2}{2-s\tilde{c}_1\alphaol}\right]\cdot \left(\frac{1}{k}\right)^{\frac{2s\tilde{c}_1\alphaol}{2+s\tilde{c}_1\alphaol}}\\ &\hspace{15em}\mbox{if}\ c_{\tau,\eta} \alphaol < \frac{2}{s}\left(1+\frac{\eta}{\olrho}\right),\\
                    &\frac{4\delta_1}{3}\cdot \frac{1}{k} + 4\tilde{C}_{f,m,U,\tau,\eta,c_1,\gamma}(\alphaou)\alphaou^2\cdot \frac{\ln(k)}{k}\\ &\hspace{15em}\mbox{if}\ c_{\tau,\eta} \alphaol = \frac{2}{s}\left(1+\frac{\eta}{\olrho}\right),\\
                    &\frac{4\delta_1 }{1+s\tilde{c}_1\alphaol}\cdot \left(\frac{1}{k}\right)^{\frac{2s\tilde{c}_1\alphaol}{2+s\tilde{c}_1\alphaol}} + \frac{4\tilde{C}_{f,m,U,\tau,\eta,c_1,\gamma}(\alphaou)(2+s\tilde{c}_1\alphaol)\alphaou^2}{s\tilde{c}_1\alphaol - 2}\cdot\frac{1}{k}\\ &\hspace{15em}\mbox{if}\ c_{\tau,\eta} \alphaol > \frac{2}{s}\left(1+\frac{\eta}{\olrho}\right),
                \end{cases}
            \end{split}
        \end{equation}
        holds for all $k\in \mb{Z}_{+}$, where $\tilde{C}_{f,m,U,\tau,\eta,c_1,\gamma}(\cdot)$ is defined in (\ref{equ:sipprateboundconstant1}).
    \end{enumerate}
\end{theorem}
By focusing on the dominant terms in Theorem \ref{thm:sippmdiminishingquadratic} and following the proof arguments of \cite[Corollary~2.5]{zhu2025tight}, we further derive asymptotic convergence rates for the model-based ispPPA, summarized in Corollary \ref{coro:sippratequadratic} below.
\begin{coro}\label{coro:sippratequadratic}
    Consider the setting of Theorem \ref{thm:sippmdiminishingquadratic} and let $s\in (0,2)$ be an arbitrary scalar. Denote by $c_{\tau,\eta} \triangleq c_1 - \frac{\tau + \eta}{2}$ and $\tilde{c}_1 \triangleq \frac{\olrho c_{\tau,\eta}}{\olrho + \eta}$. If $c_1>\frac{\tau+\eta}{2}$, then on the event that $\sup_{k\in \mb{Z}_{+}}\Norm{x_k}_2 < \infty$, the following assertions hold: 
    \begin{enumerate}
        \item If $\beta\in (0,1)$, then 
        $\mb{E}\left[\mr{dist}(x_k,\mc{X}^{*})^2\right] \leq \mc{O}\left(k^{-\beta}\right)$ 
        holds for all $k\in \mb{Z}_{+}$.
        \item If $\beta = 1$, then 
        \begin{align*}
            \mb{E}\left[\mr{dist}(x_k,\mc{X}^{*})^2\right] \leq \begin{cases}
                \mc{O}\left(k^{-\frac{2s\tilde{c}_1\alphaol}{2+s\tilde{c}_1\alphaol}}\right)\quad &\mbox{if}\ c_{\tau,\eta} \alphaol < \frac{2}{s}\left(1+\frac{\eta}{\olrho}\right),\\
                    \mc{O}\left(k^{-1}\ln(k)\right)\quad &\mbox{if}\ c_{\tau,\eta} \alphaol = \frac{2}{s}\left(1+\frac{\eta}{\olrho}\right),\\
                    \mc{O}\left(k^{-1}\right)\quad &\mbox{if}\ c_{\tau,\eta} \alphaol > \frac{2}{s}\left(1+\frac{\eta}{\olrho}\right),
            \end{cases}
        \end{align*}
        holds for all $k\in \mb{Z}_{+}$.
    \end{enumerate}
\end{coro}
Under Assumption \ref{asm:stomodel} and \ref{asm:sppaall}, the objective function $\phi$ is $(\oleta+\oltau)$-weakly convex\footnote{See Lemma \ref{lema:lipoff}(\romannumeral1) in Appendix \ref{secapx:proofofsippstability}.}. For any $\hrho>\oleta+\oltau$, it follows from \cite[Proposition~3.1~and~3.3]{hoheisel2010proximal} that $0\in \partial \phi(\olx)$ if and only if $\olx = \prox_{(\hrho)^{-1} \phi}(\olx)$, and the proximal mapping $\prox_{(\hrho)^{-1} \phi}(\cdot)$ is $\frac{\hrho}{\hrho-\oleta-\oltau}$-Lipschitz continuous. Hence, by applying the proof techniques of \cite[Lemma~2.6]{zhu2025tight}, we obtain the following estimate:
\begin{align*}
    \Norm{x-\prox_{\left(\hrho\right)^{-1} \phi}(x)}_2 \leq \left(1+\frac{\hrho}{\hrho-\oleta-\oltau}\right)\dist\left(x,\mc{X}^{*}\right)\quad \mbox{for all}\ x\in \mb{R}^d,
\end{align*}
which combined with Corollary \ref{coro:sippratequadratic} yields asymptotic convergence rates in terms of the KKT residual as follows.
\begin{coro}\label{coro:sippratequadratickkt}
    Consider the setting of Theorem \ref{thm:sippmdiminishingquadratic} and let $s\in (0,2)$ be an arbitrary scalar. Denote by $\eta \triangleq L_\infty^2\cdot\oleta$ and $\tau \triangleq L_\infty^2\cdot \oltau$. Set $c_{\tau,\eta} \triangleq c_1 - \frac{\tau+\eta}{2}$ and $\tilde{c}_1 \triangleq \frac{\olrho c_{\tau,\eta}}{\olrho + \eta}$. If $c_1>\frac{\tau+\eta}{2}$, then on the event that $\sup_{k\in \mb{Z}_{+}}\Norm{x_k}_2 < \infty$, the following assertions hold for any fixed $\hrho>\oleta+\oltau$: 
    \begin{enumerate}
        \item If $\beta\in (0,1)$, then 
        $\mb{E}\left[\Norm{x_k - \prox_{(\hrho)^{-1}\phi}\left(x_k\right)}_2\right] \leq \mc{O}\left(k^{-\frac{\beta}{2}}\right)$ 
        holds for all $k\in \mb{Z}_{+}$.
        \item If $\beta = 1$, then 
        \begin{align*}
            \mb{E}\left[\Norm{x_k - \prox_{\left(\hrho\right)^{-1}\phi}\left(x_k\right)}_2\right] \leq \begin{cases}
                \mc{O}\left(k^{-\frac{s\tilde{c}_1\alphaol}{2+s\tilde{c}_1\alphaol}}\right)\quad &\mbox{if}\ c_{\tau,\eta} \alphaol < \frac{2}{s}\left(1+\frac{\eta}{\olrho}\right),\\
                    \mc{O}\left(k^{-{\frac{1}{2}}}\left(\ln(k)\right)^{\frac{1}{2}}\right)\quad &\mbox{if}\ c_{\tau,\eta} \alphaol = \frac{2}{s}\left(1+\frac{\eta}{\olrho}\right),\\
                    \mc{O}\left(k^{-\frac{1}{2}}\right)\quad &\mbox{if}\ c_{\tau,\eta} \alphaol > \frac{2}{s}\left(1+\frac{\eta}{\olrho}\right),
            \end{cases}
        \end{align*}
        holds for all $k\in \mb{Z}_{+}$.
    \end{enumerate}
\end{coro}
\begin{remark}
    For the convex case, when applying the preconditioned isPPA (i.e., $\oleta = \oltau = 0$) with diminishing stepsizes $\alpha_k = \alpha_0 k^{-\beta}$ and parameters $\epsilon_k = \gamma \alpha_k^2$, where $\alpha_0\in \mb{R}_{++}$, $\beta\in (\tfrac{1}{2},1]$ and $\gamma\in \mb{R}_{+}$, Theorem \ref{thm:sippstabilityconvex} ensures that under Assumption \ref{asm:sppbdoptset}, with probability $1$, the iterate sequence $\{x_k\}$ is bounded and converges to some $x^{*}\in \mc{X}^{*}$. In this setting, as shown in \cite{zhu2025tight}, the iterates $\{x_k\}$ also satisfy criterion (\ref{equ:criteriaa}) with accuracy parameter $\tilde{\epsilon}_k \triangleq\tilde{\gamma}\alpha_k^{3/2}$, where $\tilde{\gamma}\triangleq \gamma \alpha_0^{1/2}$. Consequently, even without assuming $\sup_{k\in \mb{Z}_{+}}\Vert x_k\Vert_2 < \infty$, the convergence rates derived in Theorem \ref{thm:sippmdiminishingquadratic}, Corollary \ref{coro:sippratequadratic} and \ref{coro:sippratequadratickkt} hold with probability $1$. Furthermore, since $\{\dist(x_k,\mc{X}^{*})\}$ converges to zero almost surely, Assumption \ref{asm:sppa2} in Theorem \ref{thm:sippmdiminishing} can be replaced with its localized version:
    \vspace{0.5em}
    \begin{adjustwidth}{1em}{1em}
        \emph{The objective function $\phi$ satisfies the quadratic growth condition on $\mc{X}^{*}$ locally, that is, there exists $\overline{c}_1>0$ and $\delta>0$ such that 
        \begin{align*}
            \phi(x)\geq \phi^{*} + \overline{c}_1 \mr{dist}(x,\mc{X}^{*})^2\quad \mbox{for all}\ x\in \mc{U}\left(\mc{X}^{*},\delta\right),
        \end{align*}
        where $\mc{U}(\mc{X}^{*},\delta)$ denotes the set $\{x\in\mb{R}^d\ \vert\ \mr{dist}(x,\mc{X}^{*})\leq \delta\}$.
        }
    \end{adjustwidth}
    \vspace{0.5em}
    Under this localized quadratic growth condition, analogous convergence results remain valid for sufficiently large $k$. 
\end{remark}

\paragraph*{Proof of Theorem \ref{thm:sippmconstantquadratic} and \ref{thm:sippmdiminishingquadratic}.} 
By following an argument nearly identical to that in the proof of \cite[Lemma~C.1-C.3]{zhu2025tight}, we can prove that the following three technical lemmas remain valid.
\begin{lema}\label{lema:constantrecusion}
    Consider the scalar sequence $\{\delta_k\}$ satisfying the recursive relation:
    \begin{align*}
        \delta_{k+1}\leq a_1 \delta_k + a_2\quad \mbox{for all}\ k\in \mb{Z}_{+}
    \end{align*}
    for some constants $a_1 \in (0,1)$ and $a_2 \in \mathbb{R}_{+}$. Then it holds that:
    \begin{align*}
        \delta_k \leq a_1^{k-1}\delta_1 + a_2\left(1-a_1\right)^{-1}\quad \mbox{for all}\ k\in \mb{Z}_{+}.
    \end{align*}
\end{lema}
\begin{lema}\label{lema:diminishingrecusion}
    Consider the scalar sequence $\{\delta_k\}$ satisfying the recursive relation:
    \begin{align*}
        \delta_{k+1}\leq \mr{exp}(-a_1\alpha_k) \delta_k + a_2\alpha_k^2\quad \mbox{with}\quad \alphaol k^{-\beta}\leq \alpha_k \leq \alphaou k^{-\beta}\quad \mbox{for all}\ k\in \mb{Z}_{+},
    \end{align*}
    where $\beta\in (0,1)$, $\alphaol,\alphaou\in \mb{R}_{++}$, $a_1\in \mb{R}_{++}$ and $a_2 \in \mathbb{R}_{+}$. Then it holds that 
    \begin{align*}
        \delta_k \leq\ &\mr{exp}\left(-a_1\alphaol\varsigma_{1-\beta}(k)\right)\delta_1\\ &\hspace{0em} + \frac{4^\beta a_2\alphaou^2}{1-\exp\left(-a_1\alphaou\right)}\cdot \frac{1}{k^\beta} + 4^\beta a_2 \alphaou^2 \mr{exp}\left(-\frac{a_1\alphaol}{2}k^{1-\beta}\right)\cdot \varsigma_{1-2\beta}\left(\frac{k}{2}\right)
    \end{align*}
    for all $k\in \mb{Z}_{+}$.
\end{lema}
\begin{lema}\label{lema:diminishingrecusionbeta1}
    Consider the scalar sequence $\{\delta_k\}$ satisfying the recursive relation:
    \begin{align*}
        \delta_{k+1}\leq \frac{1}{1+a_1\alpha_k} \delta_k + a_2\frac{\alpha_k^2}{1+a_1\alpha_k}\quad \mbox{with}\quad \alphaol k^{-1}\leq \alpha_k \leq \alphaou k^{-1}\quad \mbox{for all}\ k\in \mb{Z}_{+},
    \end{align*}
    where $\alphaol,\alphaou\in \mb{R}_{++}$, $a_1\in \mb{R}_{++}$ and $a_2 \in \mathbb{R}_{+}$. Then it holds that 
    \begin{align*}
        \delta_k\leq \begin{cases}
            \left[\frac{4\delta_1}{1+a_1\alphaol} + \frac{4a_2(2+a_1\alphaol)\alphaou^2}{2-a_1\alphaol}\right]\cdot \left(\frac{1}{k}\right)^{\frac{2a_1\alphaol}{2+a_1\alphaol}}\quad &\mbox{if}\ a_1\alphaol < 2,\\
            \frac{4\delta_1}{3}\cdot \frac{1}{k} + 4a_2\alphaou^2\cdot \frac{\ln(k)}{k}\quad &\mbox{if}\ a_1\alphaol = 2,\\
            \frac{4\delta_1}{1+a_1\alphaol}\cdot \left(\frac{1}{k}\right)^{\frac{2a_1\alphaol}{2+a_1\alphaol}} + \frac{4a_2(2+a_1\alphaol)\alphaou^2}{a_1\alphaol - 2}\cdot \frac{1}{k}\quad &\mbox{if}\ a_1\alphaol > 2.
        \end{cases}
    \end{align*}
\end{lema}
Furthermore, on the event that $\sup_{k\in \mb{Z}_{+}}\Norm{x_k}_2 < \infty$, and under the setting of Theorem \ref{thm:sippmconstantquadratic} or \ref{thm:sippmdiminishingquadratic}, by following an argument nearly identical to that in the proof of Theorem \ref{thm:sippstabilityconvex}, we can show that (\ref{equ:boundtlxkp1}) holds, i.e., $\sup_{k\in \mb{Z}_{+}}\Vert \tlx_{k+1}\Vert_2 < \infty$. 
Denote by $U$ a bounded open convex subset of $V$ such that both $\{x_k\}$ and $\{\tlx_k\}$ are contained in $U$. Let $\delta_{k}\triangleq \mb{E}[\nu_k^{-2}\mr{dist}_{M_{k-1}}(x_k,\mc{X}^{*})^2]$ for each $k\in \mb{Z}_{+}$. Recall that $M_0 \triangleq I_d$ and $\rho_0 \triangleq 1$. Then, $\nu_1 = \rho_0 = 1$ and 
\begin{equation}\label{equ:delta1}
    \delta_1 = \mb{E}\left[\nu_1^{-2}\dist_{M_0}\left(x_1,\mc{X}^{*}\right)^2\right] = \mb{E}\left[\dist\left(x_1,\mc{X}^{*}\right)^2\right]=\dist\left(x_1,\mc{X}^{*}\right)^2,
\end{equation}
that is, the definition of $\delta_1$ coincides with the one given in the statement of Theorem \ref{thm:sippmconstantquadratic} and \ref{thm:sippmdiminishingquadratic}. Fix any $k\in \mb{Z}_{+}$. Since Assumption \ref{asm:mk} implies $\Vert \cdot \Vert_2 \leq L_\infty \Vert \cdot\Vert_{M_{k-1}}$ holds for all $k\in \mb{Z}_{+}$, we can deduce from condition $\sup_{k\in \mb{Z}_{+}}\nu_k \leq \nu_\infty$ that
\begin{align*}
    \nu_k^{-1}\inf_{x^{*}\in \mc{X}^{*}} \Norm{x-x^{*}}_{M_{k-1}} &
        \geq \nu_\infty^{-1}\inf_{x^{*}\in \mc{X}^{*}} \Norm{x-x^{*}}_{M_{k-1}}
        \geq \nu_\infty^{-1}\inf_{x^{*}\in \mc{X}^{*}} L_\infty^{-1}\Norm{x-x^{*}}_2
\end{align*}
and thus 
$\mb{E}[\dist(x_k,\mc{X}^{*})^2] \leq (\nu_\infty L_\infty)^2\cdot \mb{E}[\nu_k^{-2}\dist_{M_{k-1}}(x_k,\mc{X}^{*})^2]$. 
It then follows from the definition of $\delta_k$ that 
\begin{equation}\label{equ:distkdeltak}
    \mb{E}\left[\dist\left(x_k,\mc{X}^{*}\right)^2\right] \leq \left(\nu_\infty L_\infty\right)^2\cdot \delta_k\quad \mbox{for all}\ k\in \mb{Z}_{+}.
\end{equation}
The proofs of Theorem \ref{thm:sippmconstantquadratic} and \ref{thm:sippmdiminishingquadratic} proceed by establishing the validity of the following deterministic recursive relation, with its proof deferred to Appendix \ref{subsubsecapx:proofrecursive}.
\begin{claim}\label{claim:sippinequ1imply}
    Let Assumption \ref{asm:stomodel}-\ref{asm:sppa2} hold, and let the iterates $\{x_k\}$ be generated by the model-based ispPPA (Algorithm \ref{algo:ispppa}) with stepsizes $\{\alpha_k\}$ and parameters $\{\epsilon_k\}$. Denote by $\eta \triangleq L_\infty^2\cdot\oleta$ and $\tau \triangleq L_\infty^2\cdot\oltau$. Suppose that there exists $\gamma\in \mb{R}_{+}$ and $\olrho > \eta + \tau$ such that $\sup_{k\in \mb{Z}_{+}}\alpha_k \leq \alphaou$ for some $\alphaou\in (0,(\olrho)^{-1})$, and $\epsilon_k = \gamma \alpha_k^{\frac{3}{2}}$ for each $k\in \mb{Z}_{+}$. Set $c_{\tau,\eta} \triangleq c_1 - \frac{\tau+\eta}{2}$ and $\tilde{c}_1 \triangleq \frac{\olrho c_{\tau,\eta}}{\olrho+\eta}$. Define $\delta_{k}\triangleq \mb{E}[\nu_k^{-2}\mr{dist}_{M_{k-1}}(x_k,\mc{X}^{*})^2]$ for each $k\in \mb{Z}_{+}$. Let $s\in (0,2)$ be an arbitrary scalar. Then on the event that $\sup_{k\in \mb{Z}_{+}}\Norm{x_k}_2 < \infty$, we have   
    \begin{align*}
        \delta_{k+1} \leq \frac{1}{1+s\tilde{c}_1\alpha_k} \delta_k + \tilde{C}_{f,m,U,\tau,\eta,c_1,\gamma}\left(\alpha_k\right)\cdot\frac{\alpha_k^2}{1+s\tilde{c}_1\alpha_k}\quad \mbox{for all}\ k\in \mb{Z}_{+}
    \end{align*}
    for all $k\in \mb{Z}_{+}$.
\end{claim}
\begin{proof}[Proof of Theorem \ref{thm:sippmconstantquadratic} and \ref{thm:sippmdiminishingquadratic}]
Since $c_1>\frac{\tau+\eta}{2}$ by assumption, it follows that $\tilde{c}_1 > 0$ and $\tilde{C}_{f,m,\tau,\eta,c_1,\gamma}(\cdot)$ is nondecreasing on $\mb{R}_{+}$. Then, under the setting of Theorem \ref{thm:sippmconstantquadratic}, applying Claim \ref{claim:sippinequ1imply} with $\alphaou = \alpha_0$ yields 
\begin{align}\label{equ:recursiveconstantqg}
    \delta_{k+1}\leq \frac{1}{1+s\tilde{c}_1\alpha_0} \delta_k + \tilde{C}_{f,m,U,\tau,\eta,c_1,\gamma}\left(\alpha_0\right)\cdot\frac{\alpha_0^2}{1+s\tilde{c}_1\alpha_0}
\end{align}
for all $k\in \mb{Z}_{+}$. Similarly, under the setting of Theorem \ref{thm:sippmdiminishingquadratic}, it follows from Claim \ref{claim:sippinequ1imply} combined with condition $\alpha_k\leq \alphaou$ and inequality (\ref{inequ:exp}) that 
\begin{subequations}
    \begin{align}
        \delta_{k+1} &\leq \frac{1}{1+s\tilde{c}_1\alpha_k} \delta_k + \tilde{C}_{f,m,U,\tau,\eta,c_1,\gamma}\left(\alphaou\right)\cdot\frac{\alpha_k^2}{1+s\tilde{c}_1\alpha_k}\label{inequ:recurbeta1}\\
        &\leq \mr{exp}\left(-\frac{s\tilde{c}_1\alpha_k}{1+s\tilde{c}_1\alphaou}\right) \delta_k + \tilde{C}_{f,m,U,\tau,\eta,c_1,\gamma}\left(\alphaou\right)\cdot\alpha_k^2\label{inequ:recurbetal1}
    \end{align}
\end{subequations}
for all $k\in \mb{Z}_{+}$. By applying Lemma \ref{lema:constantrecusion}, \ref{lema:diminishingrecusion} and \ref{lema:diminishingrecusionbeta1}, and following an argument nearly identical to that in the proof of \cite[Theorem~2.3~and~2.4]{zhu2025tight}, we can deduce from (\ref{equ:recursiveconstantqg}), (\ref{inequ:recurbetal1}) and (\ref{inequ:recurbeta1}) that Theorem \ref{thm:sippmconstantquadratic} and \ref{thm:sippmdiminishingquadratic} remain valid.
\end{proof}

\section{Numerical Experiments}\label{sec:sippnumericalexperiments}
In this section, we present numerical results to validate the theoretical convergence rate guarantees derived in Section \ref{sec:convergenceispp}. We generate the synthetic data by setting
\begin{align*}
    b_i = \begin{cases}
        1\quad &\mbox{if}\ (Ax^{*}+\sigma \xi)_i \geq 0,\\
        -1\quad &\mbox{otherwise},
    \end{cases}
\end{align*}
for the logistic regression model, and $b = Ax^{*} + \sigma \xi$ for the linear regression model. Here, $A$ is sampled from the standard normal distribution, $x^{*}$ denotes the predefined true solution, $\xi$ is a Gaussian noise vector, and $\sigma$ specifies the noise level. The number of nonzero elements in $x^{*}$ is set to $\lfloor \rho_s d\rfloor$ with $\rho_s = 0.1$. The results are averaged over $3$ trials with consistent parameters and initialization across each trial of the model-based ispPPA. 

\subsection{$\ell_1$ Regularized Logistic Regression Model}
\begin{figure}[htbp]
    \centering
    \subfigure[Noiseless setting]{
    \label{fig:synnoiseless}
    \begin{minipage}[b]{0.23\textwidth}
        \centering
        \includegraphics[width=1\linewidth]{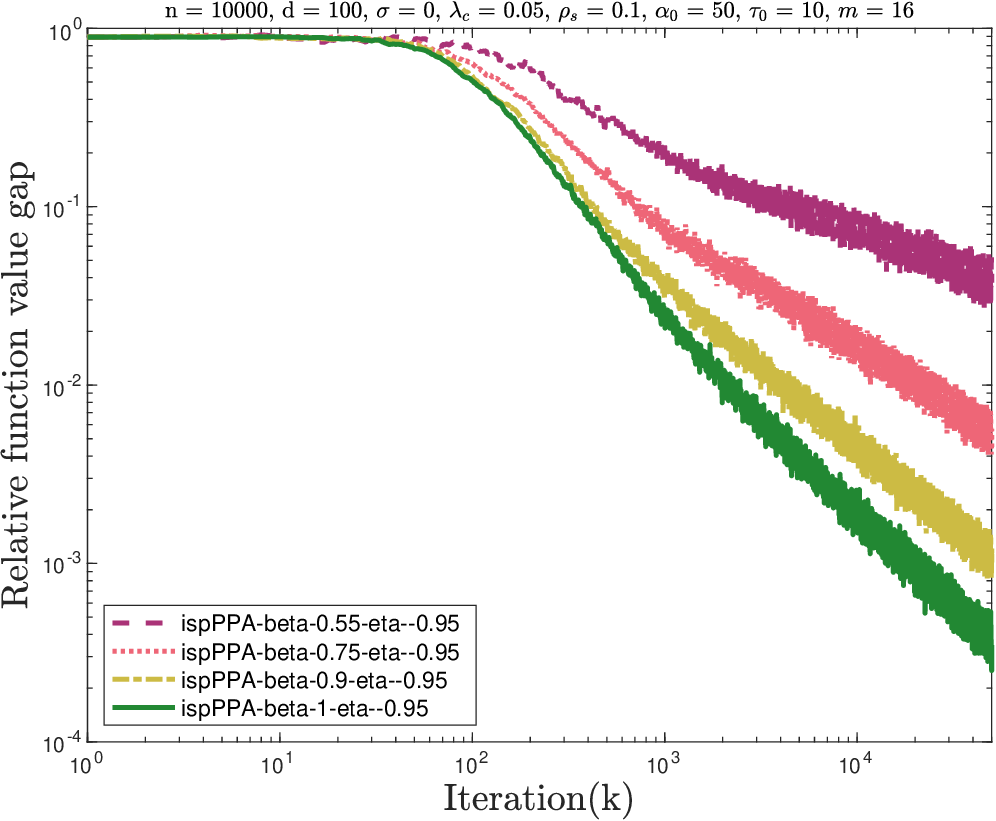}
    \end{minipage}
    \begin{minipage}[b]{0.23\textwidth}
        \centering
        \includegraphics[width=1\linewidth]{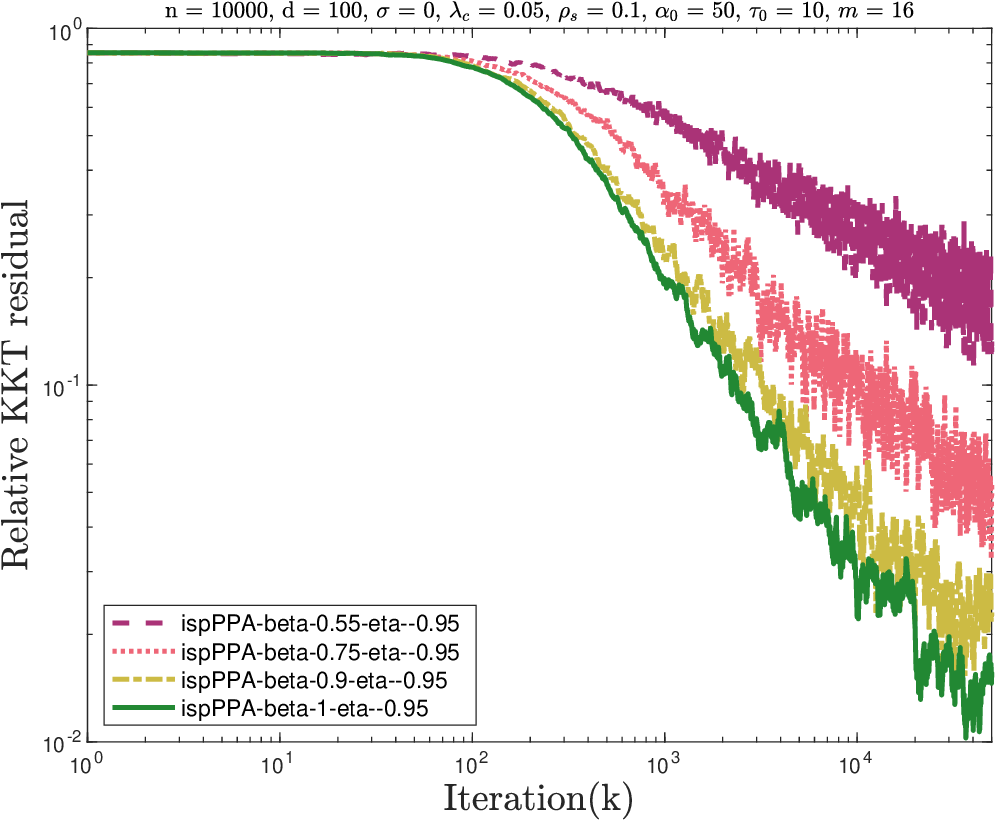}
    \end{minipage}
    }
    \subfigure[Noisy setting]{
    \label{fig:synnoisy}
    \begin{minipage}[b]{0.23\textwidth}
        \centering
        \includegraphics[width=1\linewidth]{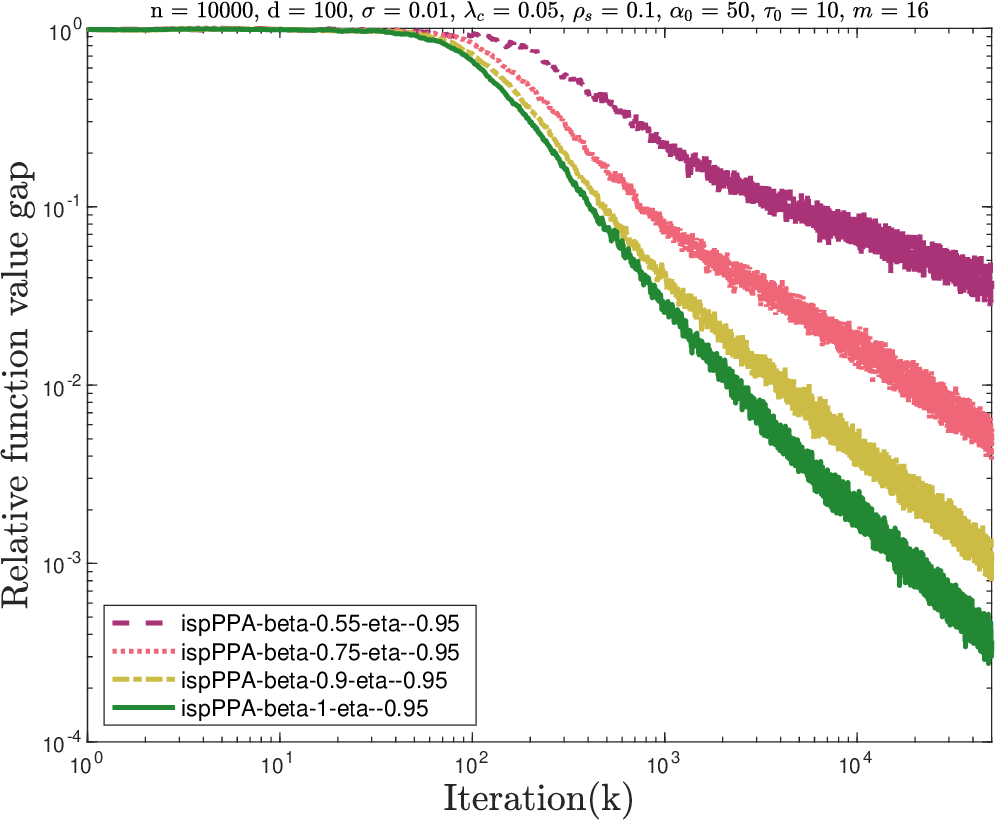}
    \end{minipage}
    \begin{minipage}[b]{0.23\textwidth}
        \centering
        \includegraphics[width=1\linewidth]{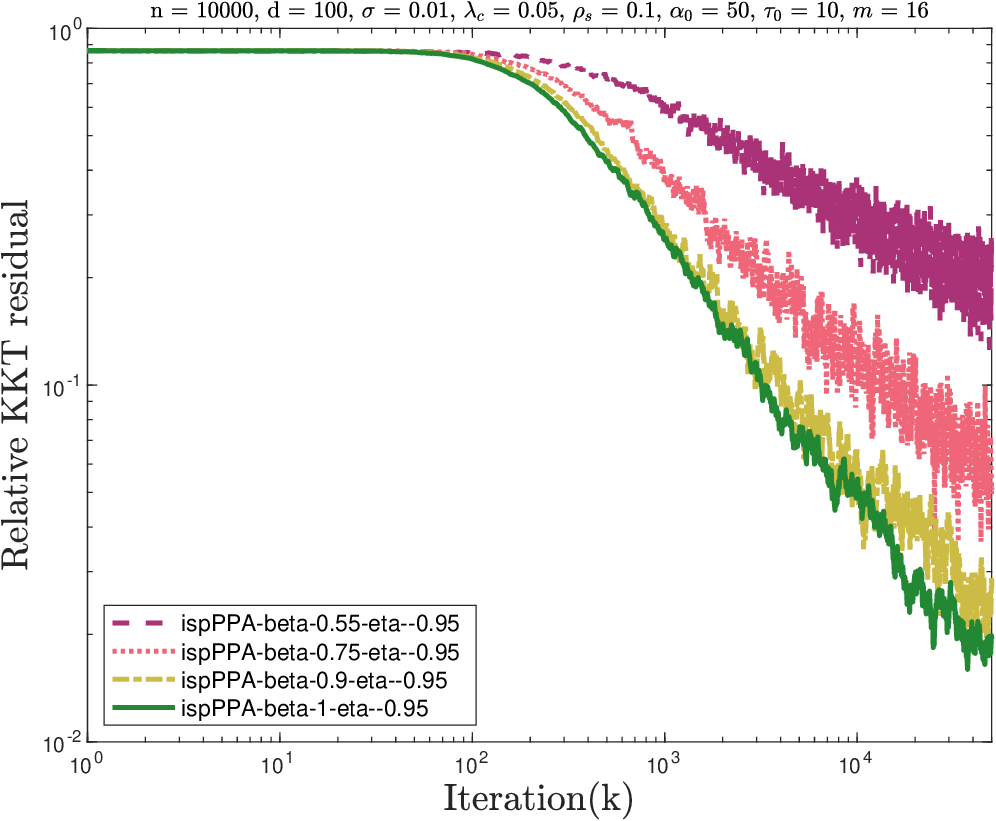}
    \end{minipage}
    }
    
    \caption{Performance of the preconditioned isPPA in solving the $\ell_1$-regularized logistic regression model for four values of stepsize exponents $\beta = 0.55$, $0.75$, $0.90$ and $1$. The legend ``ispPPA-beta-$\beta$-eta-$\eta$'' denotes the preconditioned isPPA with preconditioners $\{M_k\}$ specified in (\ref{equ:precontioner}), where stepsizes $\alpha_k = \alpha_0 k^{-\beta}$ and $\tau_k = \tau_0 k^\eta$. On the y-axis, the function value gap represents $\psi_{\ell_1}(x_k)-\psi_{\ell_1}(\hat{x}^{*})$, and the KKT residual denotes $\Vert x_k - \prox_{r_{\ell_1}}(x_k - \nabla \losslog(x))\Vert_2$, where $r_{\ell_1}(\cdot) \triangleq \lambda_1 \Vert \cdot\Vert_1$ and $\losslog(\cdot) \triangleq \sum_{i=1}^n \log\left(1 + \exp\left(-b_i \langle a_i, \cdot\rangle\right)\right)$. (a) The noiseless setting with $\sigma = 0$. (b) The noisy setting with $\sigma = 10^{-2}$.}
    \label{fig:synlassodiminishing}
\end{figure}

Consider the logistic regression model with the $\ell_1$ regularizer:
\begin{align}\label{pro:loglasso}
    \min_{x\in \mb{R}^d}\ \psi_{\ell_1}(x) \triangleq \sum_{i=1}^n \log\left(1 + \exp\left(-b_i \langle a_i,x\rangle\right)\right) + \lambda_1 \Norm{x}_1,
\end{align}
where $A = \begin{pmatrix}
    a_1 &\cdots &a_n
\end{pmatrix}^\top \in \mb{R}^{n\times d}$ and $b\in \{-1,1\}^n$ are given data, and $\lambda_1\in \mb{R}_{++}$ is the regularization parameter. We apply the \emph{preconditioned isPPA}, i.e., Algorithm \ref{algo:ispppa} with model functions $f_x(\cdot;s) = f(\cdot;s)$, to solve problem (\ref{pro:loglasso}). Under this setting, Assumption \ref{asm:stomodel} and \ref{asm:sppaall} hold with $\oleta = \oltau = 0$. The preconditioners are defined as 
\begin{equation}\label{equ:precontioner}
    \begin{split}
        M_k \triangleq I_d + \alpha_k \tau_k A_{S_k^{1:m}}^\top A_{S_k^{1:m}}\quad \mbox{with}\quad A_{S_k^{1:m}} \triangleq 
        \begin{pmatrix}
            a_i 
        \end{pmatrix}_{i\in S_k^{1:m}}^\top\in \mb{R}^{m\times d}
    \end{split}
\end{equation}
for each $k\in \mb{Z}_{+}$, where $\{\alpha_k\}$ denotes the stepsize sequence and $\{\tau_k\}$ is a sequence of positive parameters specified later. With this choice of preconditioners, each inner-loop subproblem (\ref{equ:sippm}) is solved inexactly by the efficient semismooth Newton (SSN) method \cite[Algorithm~SSN]{li2018highly} such that the iterate sequence $\lrbrackets{x_k}$ satisfies criterion (\ref{equ:criteriaa}); see Claim \ref{claim:ippaequialmnew} in Appendix \ref{apx:subsecinnerproblem} for details. In all experiments, we set $n = 10000$ and $d = 100$, using $\sigma=0$ for the noiseless case and $\sigma=10^{-2}$ otherwise. The regularization parameter is chosen as $\lambda_1 = \lambda_c \Vert A^\top b\Vert_\infty$ with $\lambda_c = 10^{-2}$. The accuracy parameter is set to $\epsilon_k= \gamma \alpha_k^2$ with $\gamma = 10^{-2}$, and the minibatch size is fixed at $m=16$. Fig. \ref{fig:synlassodiminishing} illustrates the convergence behavior of the preconditioned isPPA with diminishing stepsizes $\alpha_k = \alpha_0 k^{-\beta}$ for various stepsize exponents $\beta\in \{0.55,0.75,0.90,1\}$, where the initial stepsize is $\alpha_0 = 50$. The parameter $\tau_k$ is updated as $\tau_0 k^\eta$ with $\tau_0 = 10$ and $\eta = -0.95$. By Claim \ref{claim:choicemk} in Appendix \ref{apx:subsecpreconditioners}, the preconditioner sequence $\{M_k\}$ defined in (\ref{equ:precontioner}) satisfies Assumption \ref{asm:mk}. The numerical results demonstrate that the preconditioned isPPA with a stepsize exponent of $\beta = 1$ exhibits the best performance. Specifically, the asymptotic convergence rate measured by the squared distance to the optimal solution set is $\mc{O}(k^{-\beta})$, viewing $\psi_{\ell_1}(x_k)-\psi_{\ell_1}(\hat{x}^{*})$ as an upper bound of $\mc{O}(\dist(x_k,\mc{X}^{*})^2)$ under the quadratic growth condition. Here, $\hat{x}^{*}$ denotes an approximate solution of (\ref{pro:loglasso}) obtained by the semismooth Newton augmented Lagrangian (SSNAL) method \cite{li2018highly}. Moreover, the asymptotic convergence rate in terms of the KKT residual $\Vert x_k - \prox_{\psi_{\ell_1}}(x_k)\Vert_2$ is $\mc{O}(k^{-\beta/2})$\footnote{Claim \ref{claim:kktresidualcompute} in Appendix \ref{apx:kktresidualcompute} implies $\Vert x - \prox_{\psi_{\ell_1}}(x)\Vert_2 \leq \mc{O}(\Vert x - \prox_{r_{\ell_1}}(x - \nabla \losslog(x))\Vert_2)$.}. As all assumptions required in Corollary \ref{coro:sippratequadratic} and \ref{coro:sippratequadratickkt} are satisfied, these empirical observations confirm the theoretical results established therein. 

\subsection{MCP Regularized Linear Regression Model}
\begin{figure}[!htbp]
    \centering
    \subfigure[Noiseless setting]{
    \label{fig:synmcpnoiseless}
    \begin{minipage}[b]{0.23\textwidth}
        \centering
        \includegraphics[width=1\linewidth]{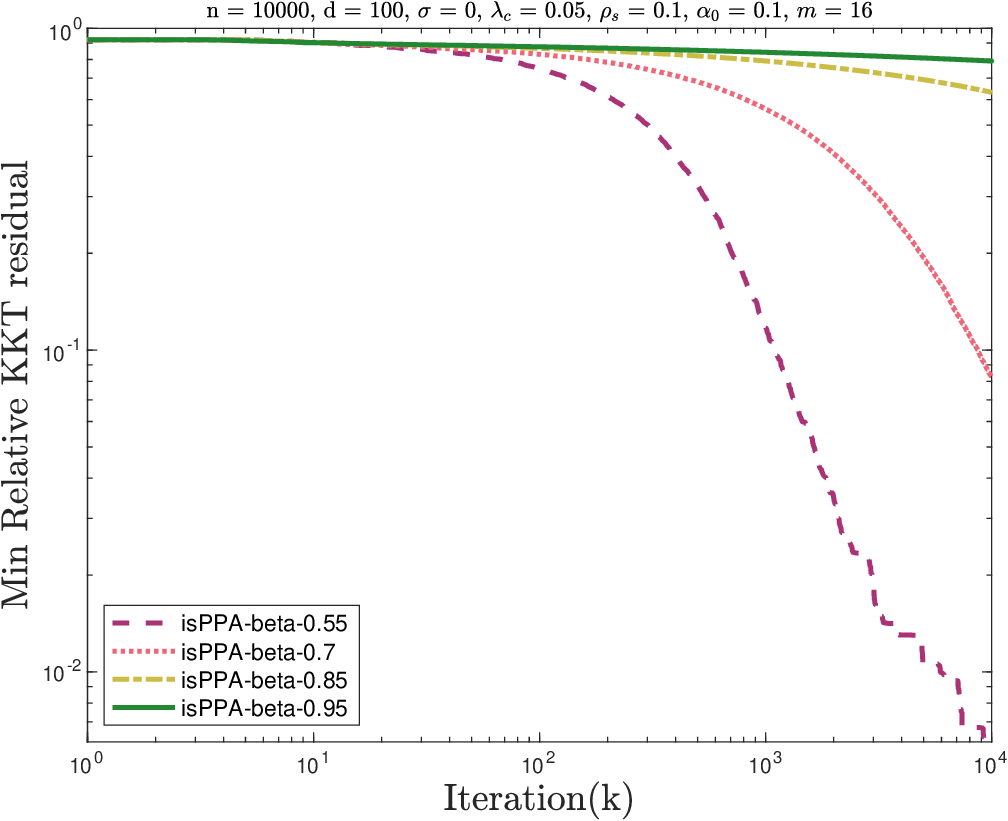}\\
        \includegraphics[width=1\linewidth]{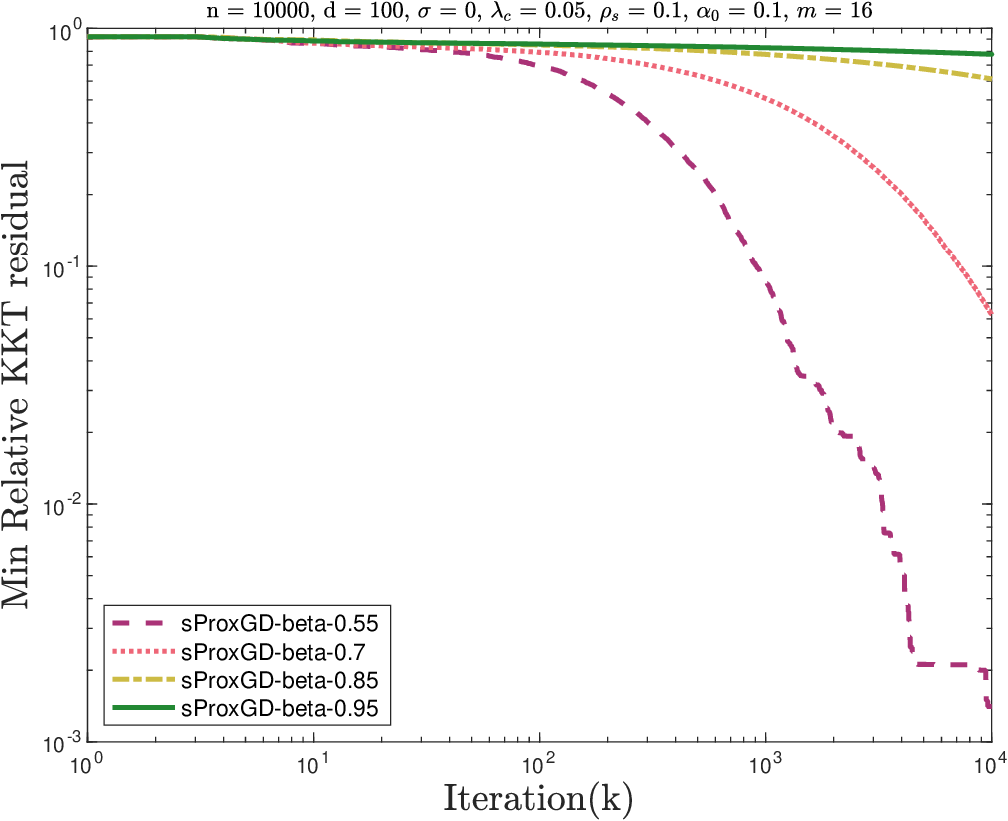}
    \end{minipage}
    \begin{minipage}[b]{0.23\textwidth}
        \centering
        \includegraphics[width=1\linewidth]{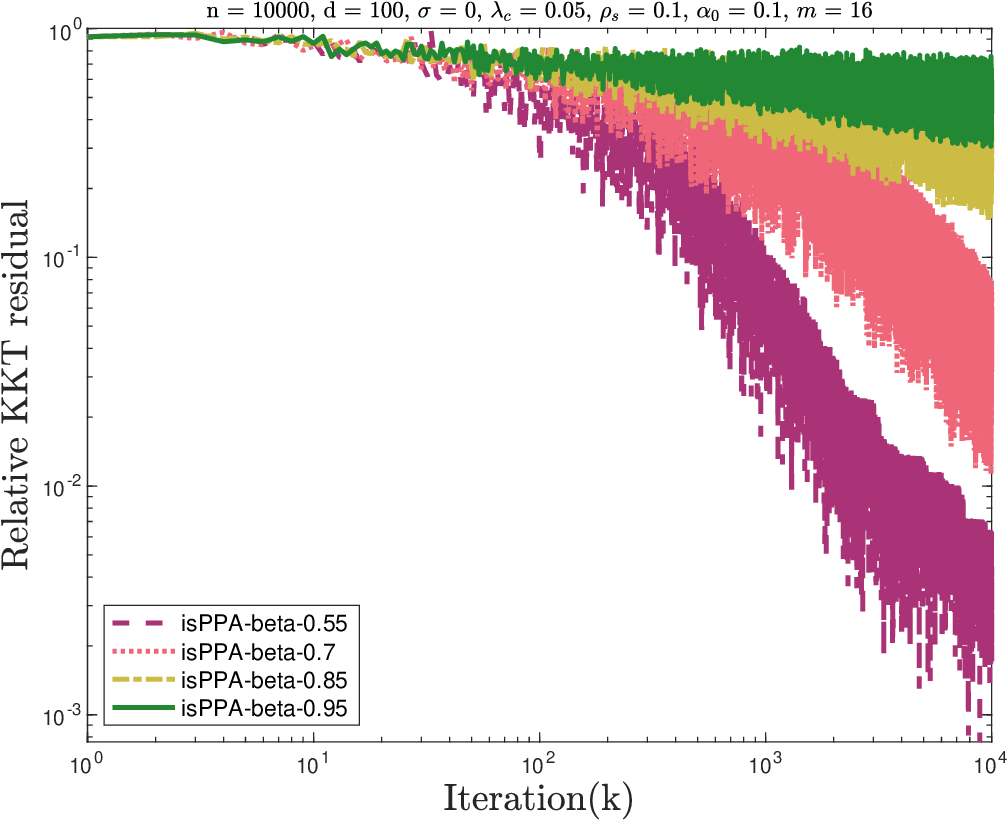}\\
        \includegraphics[width=1\linewidth]{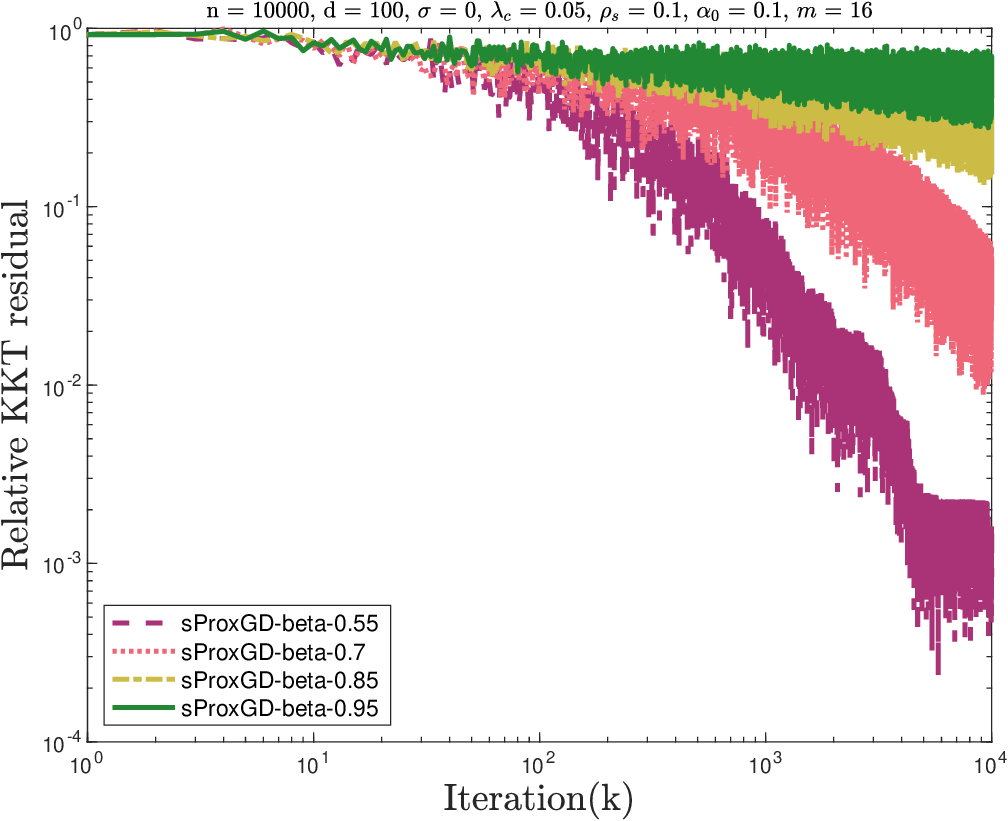}
    \end{minipage}
    }
    \subfigure[Noisy setting]{
    \label{fig:synmcpnoisy}
    \begin{minipage}[b]{0.23\textwidth}
        \centering
        \includegraphics[width=1\linewidth]{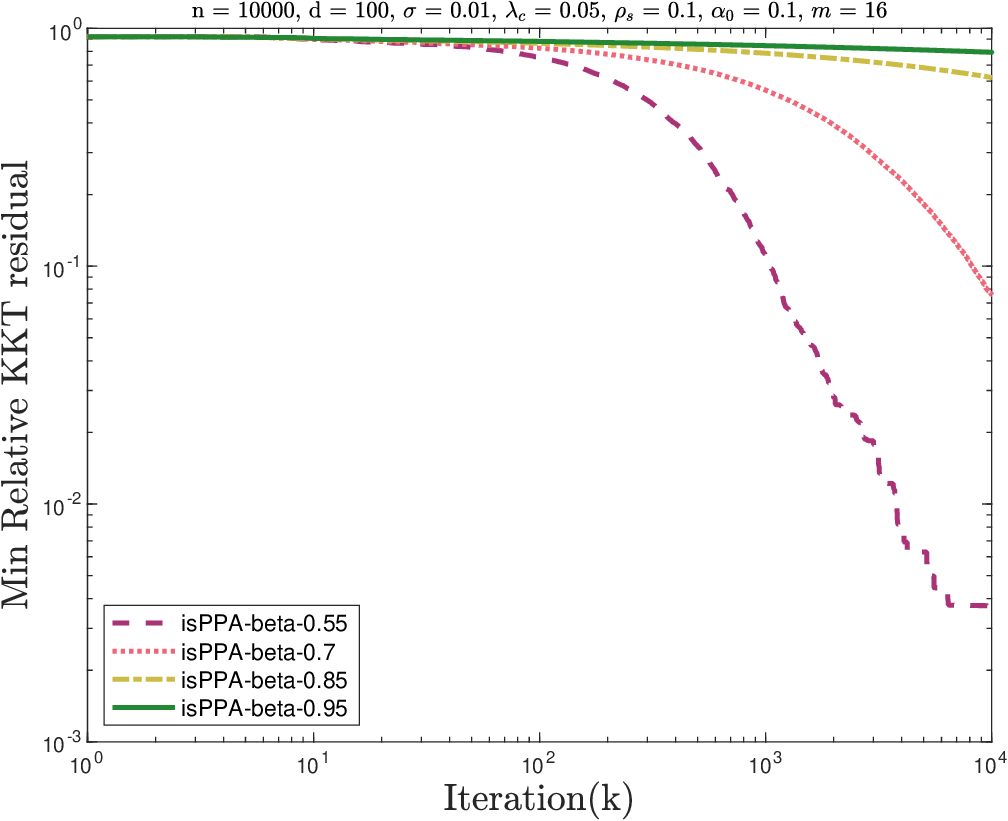}\\
        \includegraphics[width=1\linewidth]{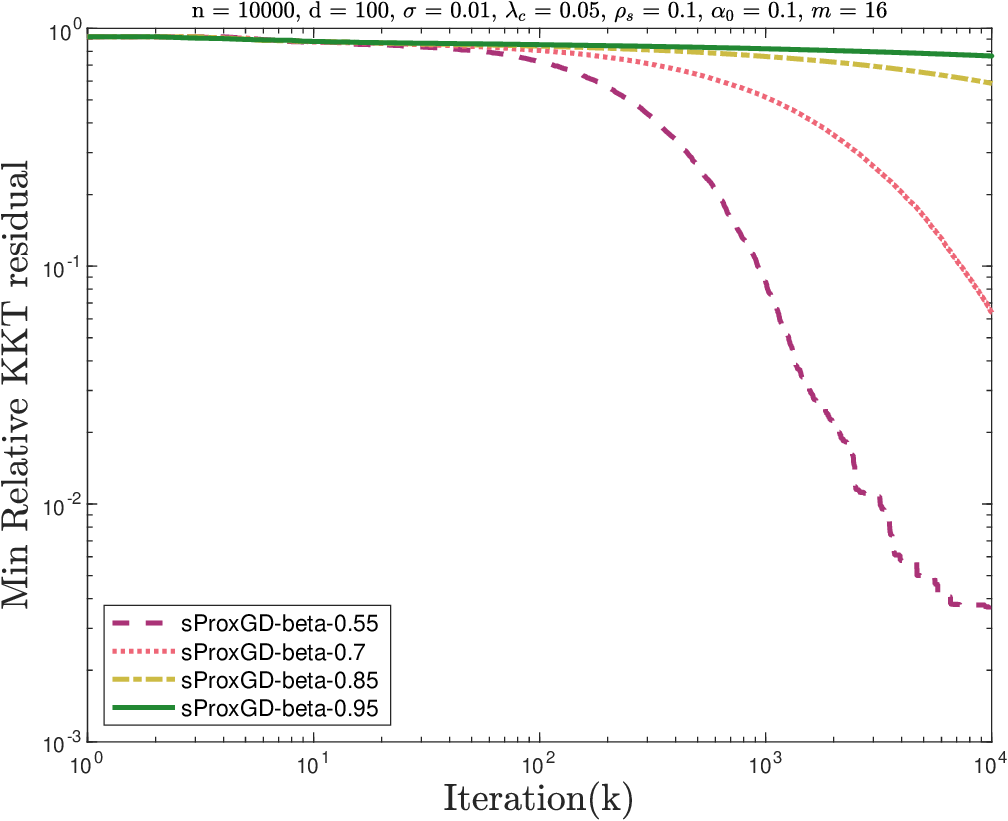}
    \end{minipage}
    \begin{minipage}[b]{0.23\textwidth}
        \centering
        \includegraphics[width=1\linewidth]{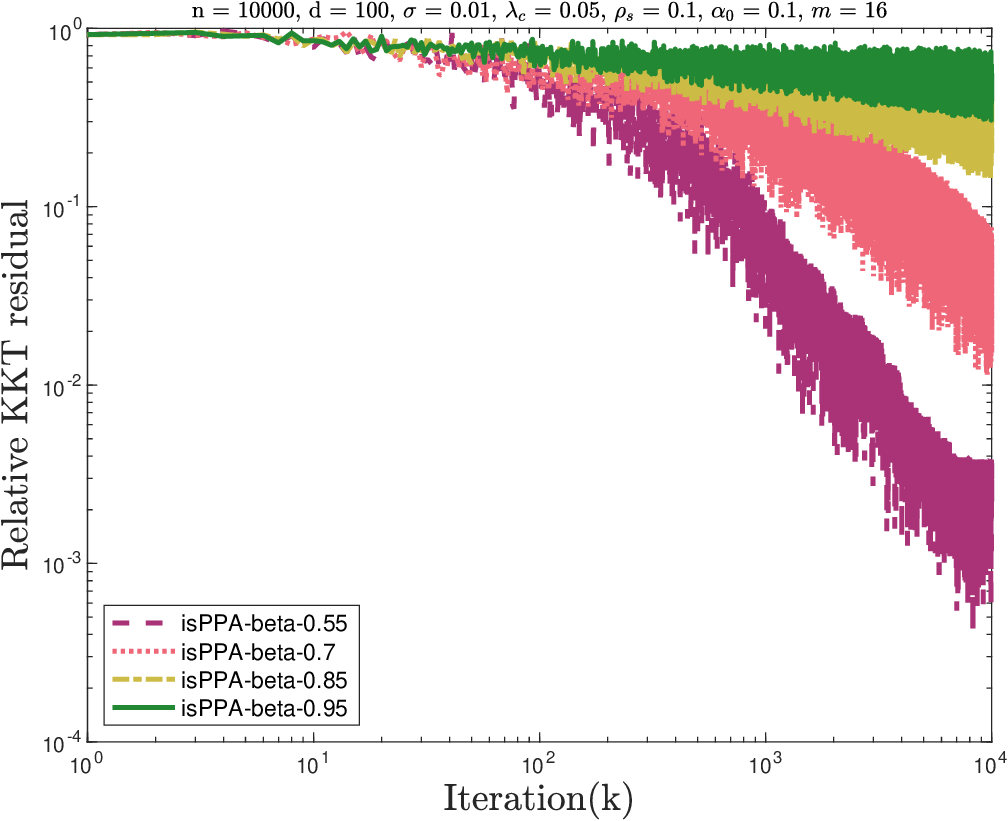}
        \includegraphics[width=1\linewidth]{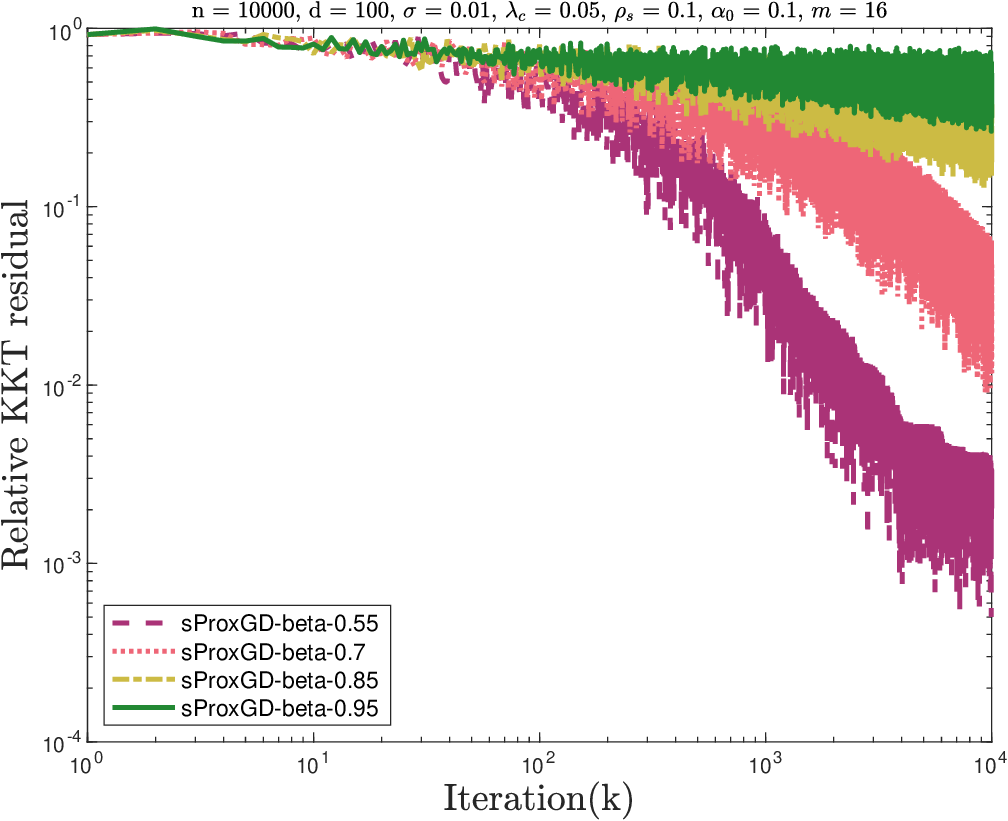}
    \end{minipage}
    }
    
    \caption{Performance of the model-based ispPPA in solving the MCP-regularized linear regression model (isPPA - top; stochastic proximal gradient - bottom) for four values of stepsize exponents $\beta = 0.55$, $0.70$, $0.85$ and $0.95$. The legend ``isPPA-beta-$\beta$'' denotes the isPPA with stepsizes $\alpha_k = \frac{\alpha_0 k^{-\beta}}{\max\{1,\Vert x_k\Vert_2\}}$. On the y-axis, the ``min KKT residual'' represents $\min_{1\leq i\leq k} \Vert x_i - \prox_{\tlrho r_{\mr{MCP}}}(x_i - \tlrho\nabla \losslin(x_i))\Vert_2$, and the ``KKT residual'' corresponds to $\Vert x_{i_{*}} - \prox_{\tlrho r_{\mr{MCP}}}(x_{i_{*}} - \tlrho\nabla \losslin(x_{i_{*}}))\Vert_2$, where $\tlrho \triangleq (\olrho)^{-1}$ with $\olrho > 2\lambda_2^{-1}$ for the isPPA and $\olrho > 2\lambda_2^{-1}+ n\max_{1\leq i\leq n}\Vert a_i\Vert_2^2$ for the stochastic proximal gradient algorithm. Here $\losslin(\cdot) \triangleq \frac{1}{2}\Vert A\cdot - b\Vert_2^2$, and $x_{i_{*}}$ denotes the iterate sampled from $\{x_i\}_{i\in [k]}$ as described in Theorem \ref{thm:sippmdiminishing2}.  (a) The noiseless setting with $\sigma = 0$. (b) The noisy setting with $\sigma = 10^{-2}$.}
    \label{fig:synmcpdiminishing}
\end{figure}

The linear regression model with the minimax concave penalty (MCP) regularizer \cite{zhang2010nearly} is defined by
\begin{align}\label{pro:logmcp}
    \min_{x\in \mb{R}^d}\ \psi_{\mr{MCP}}(x) \triangleq \frac{1}{2}\Vert Ax - b\Vert_2^2 + r_{\mr{MCP}}\left(x\right),
\end{align}
where $A = \begin{pmatrix}
    a_1 &\cdots &a_n
\end{pmatrix}^\top\in \mb{R}^{n\times d}$ and $b\in \mb{R}^n$ are given data, $\lambda_1,\lambda_2\in \mb{R}_{++}$ are regularization parameters, and 
\begin{align*}
    r_{\mr{MCP}}\left(x\right) = \sum_{i=1}^d r_{\lambda_1,\lambda_2}\left(x_i\right)\quad \mbox{with}\quad r_{\lambda_1,\lambda_2}\left(t\right) = \begin{cases}
        \lambda_1\norm{t} - \frac{t^2}{2\lambda_2}\quad &\mbox{if}\ \norm{t}\leq \lambda_1\lambda_2\\
        \frac{\lambda_2\lambda_1^2}{2}\quad &\mbox{otherwise}.
    \end{cases}
\end{align*}
The regularizer $r_{\mr{MCP}}$ is $\lambda_2^{-1}$-weakly convex, and its proximal mapping is known as the firm-threshold operator \cite{gao1997waveshrink}. 
We evaluate two algorithmic variants of the model-based ispPPA for solving this problem:
\begin{itemize}
    \item The \emph{isPPA} refers to Algorithm \ref{algo:ispppa} with model functions $f_x(\cdot;s)=f(\cdot;s)$ and trivial preconditioners $M_k \equiv I_d$. In this case, Assumption \ref{asm:stomodel} and \ref{asm:sppaall} hold with $\oleta = 0$ and $\oltau = \lambda_2^{-1}$. Each inner-loop subproblem (\ref{equ:sippm}) is solved inexactly via the SSN method, ensuring that the iterates $\lrbrackets{x_k}$ satisfy criterion (\ref{equ:criteriaa}); see Remark \ref{remark:predid} in Appendix \ref{apx:subsecinnerproblem} for details.  
    \item The \emph{stochastic proximal gradient} algorithm refers to Algorithm \ref{algo:ispppa} with model functions $f_x(\cdot;s) = f(x;s) + \langle \nabla f(x;s),\cdot-x\rangle$ and trivial preconditioners $M_k \equiv I_d$. Here, Assumption \ref{asm:stomodel} and \ref{asm:sppaall} hold with $\oleta = n\max_{1\leq i\leq n}\Vert a_i\Vert_2^2$ and $\oltau = \lambda_2^{-1}$, since each component function $f(\cdot;s) \triangleq\frac{n}{2}( a_i^\top x - b_i)^2$ is $n\Vert a_i\Vert_2^2$-Lipschitz smooth. 
\end{itemize}
All experiments use $n = 10000$ and $d = 100$, with $\sigma=0$ in the noiseless setting and $\sigma=10^{-2}$ otherwise. The regularization parameters are set as $(\lambda_1,\lambda_2)=(\lambda_{c_1} \Norm{A^\top b}_\infty,\lambda_{c_2} \Norm{A^\top b}_\infty)$ with $(\lambda_{c_1},\lambda_{c_2}) = (5\times 10^{-2},5\times 10^{-3})$. We fix the minibatch size $m = 16$ and define $\epsilon_k =\gamma \alpha_k^2$ with $\gamma = 10^{-2}$. Fig. \ref{fig:synmcpdiminishing} displays the convergence curves of the isPPA and the stochastic proximal gradient algorithm with diminishing stepsizes $\alpha_k = \alpha_0 k^{-\beta}$ for various stepsize exponents $\beta\in \{0.55,0.70,0.85,0.95\}$. For the isPPA, since $M_k \equiv I_d$, we have $\sup_{k\in \mb{Z}_{+}}\rho_k^2 = 1$, and choosing $\alpha_0 < (2\olrho)^{-1}$ for some $\olrho > 2\lambda_2^{-1}$ ensures that $\{\alpha_k\}$ satisfies the requirements of Theorem \ref{thm:sippmdiminishing} and \ref{thm:sippmdiminishing2}. Similarly, for the stochastic proximal gradient algorithm, it holds that $\sup_{k\in \mb{Z}_{+}}\rho_k^2 = 1$, and setting $\alpha_0 < (2\olrho)^{-1}$ for some $\olrho > 2\lambda_2^{-1}+ n\max_{1\leq i\leq n}\Vert a_i\Vert_2^2$ provides the same guarantee. As shown in Fig. \ref{fig:synmcpdiminishing}, both the isPPA and the stochastic proximal gradient algorithm exhibit the best performance when the stepsize exponent $\beta$ is close to $0.5$. The asymptotic convergence rates measured by the KKT residual, in both $\min_{1\leq i\leq k} \Vert x_i - \prox_{(\olrho)^{-1}\psi_{\mr{MCP}}}(x_i)\Vert_2$ and $\Vert x_{i_{*}} - \prox_{(\olrho)^{-1}\psi_{\mr{MCP}}}(x_{i_{*}})\Vert_2$, are $\mc{O}(k^{-(1-\beta)/2})$\footnote{Claim \ref{claim:kktresidualcompute} in Appendix \ref{apx:kktresidualcompute} implies $\Vert x - \prox_{\tlrho\psi_{\mr{MCP}}}(x)\Vert_2\leq \mc{O}(\Vert x - \prox_{\tlrho r_{\mr{MCP}}}(x - \tlrho\nabla \losslin(x))\Vert_2)$ with $\tlrho \triangleq (\olrho)^{-1}$.}, consistent with the theoretical findings of Theorem \ref{thm:sippmdiminishing} and \ref{thm:sippmdiminishing2}. 

\section{Conclusion}
In this paper, we consider the stochastic model-based algorithm combined with preconditioning techniques for solving stochastic composite optimization problems, where a family of variable preconditioners satisfying suitable relations is admissible, and each subproblem is allowed to be solved inexactly under implementable inner-loop stopping criteria. Without assuming global Lipschitz continuity of either the component model functions or the regularizer, we establish several sufficient conditions ensuring the stability in both weakly convex and convex settings. Under a local Lipschitz condition on the component model functions, we derive convergence rate bounds of $\mc{O}(k^{-\frac{1}{4}})$ measured by the gradient of the Moreau envelope for the weakly convex case. By further imposing a quadratic growth condition on the objective function and appropriate relations among the quadratic growth constant $c_1$, the model accuracy constant $\oleta$, and the weak convexity constant $\oltau$, we obtain convergence rate bounds of $\mc{O}(k^{-\frac{1}{2}})$ in terms of both the distance to the optimal solution set and the KKT residual. Numerical experiments corroborate the theoretical results.

\bibliography{reference}{}
\bibliographystyle{abbrvnat}

\clearpage
\appendix
\section{Comparison of Convergence Results for Stochastic Model-Based Algorithms}
Table \ref{tab:asumstability} summarizes representative sufficient conditions from existing work that guarantee the almost sure boundedness of the iterates generated by stochastic optimization algorithms for 
\begin{align*}
    \min_{x\in \mb{R}^d}\ \phi(x)\triangleq F(x) + r(x)\quad \mbox{where}\quad F(x) \triangleq \mb{E}_{s\sim P}\left[f(x;s)\right]
\end{align*}
under square summable stepsizes. Abbreviations used in the table are as follows:
\begin{itemize}
    \item \textbf{ME}: Moreau envelope of $\phi$; \textbf{subpb}: subproblem; \textbf{PreC}: preconditioner.
    \item \textbf{cvx}/\textbf{w.cvx}/\textbf{s.cvx}: convex/weakly convex/strongly convex.
    \item \textbf{Lip.sm}: Lipschitz smooth; \textbf{Lip}: globally Lipschitz continuous; \textbf{bd.subg*}: bounded subgradients over the optimal solution set.
    \item \textbf{loc.Lip*$\boldsymbol{\nu}$}, \textbf{loc.Lip*1}, \textbf{loc.Lip*2}: locally Lipschitz continuous functions defined in (\ref{equ:generalizedloclipnu}), Assumption \ref{asm:generalizedlip2}, Assumption \ref{asm:generalizedlip3}.
    \item \textbf{$\boldsymbol{\sigma}$-reg.coer}: $(\lambda,\sigma$)-regular coercivity defined in (\ref{equ:regcoercive}); \textbf{bd.sublv}: bounded sublevel sets; \textbf{blw}: bounded from below.
    \item $\boldsymbol{\mf{1}_{\mc{X}}}$: indicator function of a closed convex set $\mc{X}\subset \mb{R}^d$.
\end{itemize}
In Table \ref{tab:asumstability}, ``\textit{exact}'' and ``\textit{inexact}'' under ``\textbf{subpb}'' indicate whether each subproblem is solved exactly or inexactly, while ``$I_d$'' and ``$M_k$'' under ``\textbf{PreC}'' denote fixed-metric and preconditioned variants, respectively. 
\begin{sidewaystable}[!htbp]
    \centering
    \caption{Comparison of assumptions required for the stability of stochastic optimization algorithms.}
    \label{tab:asumstability}
    \vspace{0.3em}
    \setlength{\tabcolsep}{4pt}
    \renewcommand{\arraystretch}{1.3}
    \resizebox{\textwidth}{!}{%
    \begin{tabular}{@{}ccccccccc@{}}
            \toprule[1.1pt]
            \multirow{2}{*}{Method} & \multirow{2}{*}{Literature} & \multicolumn{7}{c}{Assumptions} \\
            \cmidrule(lr){3-9}
            &  & Loss & comp.Loss & Reg & $\phi$ & ME & subpb & PreC \\ 
            \midrule[1.05pt]
            SGD & \cite{nguyen2018sgd} & \begin{tabular}[c]{@{}c@{}}cvx \& Lip.sm\\ on $f(\cdot;s)$\end{tabular} & \_ & $r\equiv 0$ & s.cvx & \_ & exact & $I_d$ \\ \specialrule{0.2pt}{0pt}{0pt} 
            \specialrule{0.2pt}{1.5pt}{0pt}
            \begin{tabular}[c]{@{}c@{}}Stochastic\\ Proximal\\ Subgradient\end{tabular} & \cite{duchi2018stochastic} & \begin{tabular}[c]{@{}c@{}}w.cvx \& loc.Lip*$\nu$\\ on $f(\cdot;s)$\end{tabular} & \_ & \begin{tabular}[c]{@{}c@{}}real-valued \& cvx\\ $\sigma$-reg.coer\\ $\sigma > \nu + 1$\end{tabular} & \_ & \_ & exact & $I_d$ \\ \specialrule{0.2pt}{0pt}{0pt}  \specialrule{0.2pt}{1.5pt}{0pt}
            \multirow{4}{*}{\begin{tabular}[c]{@{}c@{}}Stochastic\\ PPA\end{tabular}} & \cite{zhu2025tight} & \_ & \begin{tabular}[c]{@{}c@{}}cvx \& bd.subg*\\ on $\varphi(\cdot;s)$\end{tabular} & \_ & \_ & \_ & inexact & $I_d$ \\ \cmidrule(lr){2-2} \cmidrule(lr){3-9}
            & \begin{tabular}[c]{@{}c@{}}Theorem \ref{thm:sippstabilityconvex}\\ \textbf{(our paper)}\end{tabular} & \_ & \begin{tabular}[c]{@{}c@{}}cvx \& bd.subg*\\ on $\varphi(\cdot;s)$\end{tabular} & \_ & \_ & \_ & inexact & $M_k$ \\ \specialrule{0.2pt}{0pt}{0pt}  \specialrule{0.2pt}{1.5pt}{0pt}
            \multirow{8}{*}{\begin{tabular}[c]{@{}c@{}}Model-Based\\ Stochastic\\ Optimization\end{tabular}} & \cite{asi2019stochastic} & \begin{tabular}[c]{@{}c@{}}cvx \& bd.subg*\\ on $f_x(\cdot;s)$\end{tabular} & \_ & $r=\mf{1}_{\mc{X}}$ & \_ & \_ & exact & $I_d$ \\ \cmidrule(lr){2-2} \cmidrule(lr){3-9}
            & \cite{gao2024stochastic} & \begin{tabular}[c]{@{}c@{}}cvx \& loc.Lip*1\\ on $f_x(\cdot;s)$\end{tabular} & \_ & w.cvx \& Lip & \_ & blw \& bd.sublv & exact & $I_d$ \\ \cmidrule(lr){2-2} \cmidrule(lr){3-9}
            & \begin{tabular}[c]{@{}c@{}}Theorem \ref{thm:almostsureboundednessadpstep1}\\ \textbf{(our paper)}\end{tabular} & \begin{tabular}[c]{@{}c@{}}loc.Lip*1\\ on $f_x(\cdot;s)$\end{tabular} & \begin{tabular}[c]{@{}c@{}}w.cvx\\ on $\varphi_x(\cdot;s)$\end{tabular} & \_ & blw \& bd.sublv & \_ & inexact & $M_k$ \\ \cmidrule(lr){2-2} \cmidrule(lr){3-9}
            & \begin{tabular}[c]{@{}c@{}}Theorem \ref{thm:almostsureboundednessadpstep2}\\ \textbf{(our paper)}\end{tabular} & \begin{tabular}[c]{@{}c@{}}loc.Lip*2\\ on $\varphi_x(\cdot;s)$\end{tabular} & \begin{tabular}[c]{@{}c@{}}w.cvx\\ on $\varphi_x(\cdot;s)$\end{tabular} & \_ & blw \& bd.sublv & \_ & inexact & $M_k$ \\ 
            \bottomrule[1.1pt]
    \end{tabular}%
    }
\end{sidewaystable}

Table \ref{tab:asumcounterexample} presents several illustrative examples that violate the assumptions required for the stability analysis in prior work. As discussion in Section \ref{subsec:motivationrelatedworks}, these observations demonstrate that the assumptions adopted in this paper are considerably milder and more practically reasonable than those in previous studies. Furthermore, our proposed framework accounts for the inexactness of subproblem solutions and extends the stability analysis to preconditioned variants. 

\begin{table}[!htbp]
    \centering
    \caption{Illustrative examples violating the required conditions for the stability of stochastic optimization algorithms.}
    \label{tab:asumcounterexample}
    \vspace{0.3em}
    \setlength{\tabcolsep}{4pt}
    \renewcommand{\arraystretch}{1.3}
    \begin{tabular*}{\textwidth}{@{\extracolsep\fill}ccl@{}}
        \toprule[1.1pt]
        Method & Literature & Examples not satisfy required conditions \\
        \midrule[1.05pt]
        SGD & \cite{nguyen2018sgd} & \begin{tabular}[l]{@{}l@{}}$\frac{1}{n}\sum_{i=1}^n \vert \langle a_i, x\rangle^2 - b_i\vert$\\ (\textbf{nonconvex \& nondifferentiable})\end{tabular} \\ \specialrule{0.2pt}{0pt}{0pt}  \specialrule{0.2pt}{1.5pt}{0pt}
        \begin{tabular}[c]{@{}c@{}}Stochastic\\ Proximal\\ Subgradient\end{tabular} & \cite{duchi2018stochastic} & \begin{tabular}[l]{@{}l@{}}$\frac{1}{n}\sum_{i=1}^n \vert \langle a_i, x\rangle^2 - b_i\vert + \lambda \Vert x\Vert_2^2$ \\ (\textbf{$\boldsymbol{\sigma < 2}$ and $\boldsymbol{\nu = 1}$})\end{tabular} \\ \specialrule{0.2pt}{0pt}{0pt}  \specialrule{0.2pt}{1.5pt}{0pt}
        \begin{tabular}[c]{@{}c@{}}Stochastic\\ PPA\end{tabular} & \cite{zhu2025tight} & $\frac{1}{n}\sum_{i=1}^n \vert \langle a_i, x\rangle^2 - b_i\vert$ (\textbf{nonconvex}) \\ 
        \specialrule{0.2pt}{0pt}{0pt}  \specialrule{0.2pt}{1.5pt}{0pt}
        \begin{tabular}[c]{@{}c@{}}Model-Based\\ Stochastic\\ Optimization\end{tabular} & \cite{gao2024stochastic} & \begin{tabular}[l]{@{}l@{}}$\frac{1}{n}\sum_{i=1}^n \vert \langle a_i, x\rangle^2 - b_i\vert$ using \\ model function $f_x(\cdot;s) = f(\cdot;s)$ where \\ $\vert f_x(x;s) - f_x(y;s)\vert \leq L_f \boldsymbol{\max\{\Vert x\Vert_2,\Vert y\Vert_2\}} \Vert x-y\Vert_2$\end{tabular}\\
        \bottomrule[1.1pt]
    \end{tabular*}
\end{table}

\section{Proofs for the Results in Section \ref{subsec:sippstability}}\label{secapx:proofofsippstability}
In this section, we present detailed proofs of Claim \ref{claim:sippinequhatx1tvalueglip} and \ref{claim:sippstablebound} stated in Section \ref{subsec:sippstability}. 
Recall that $\eta \triangleq L_\infty^2\cdot\oleta$, $\tau \triangleq L_\infty^2\cdot \oltau$ and $B_\rho \triangleq \sup_{k\in \mb{Z}_{+}} \rho_k^2$. 
By Remark \ref{remark:lipqgmknorm}, if Assumption \ref{asm:stomodel} and \ref{asm:mk} hold, it follows that for $P$-almost all $s\in \mc{S}$ and any $k\in \mb{Z}_{+}$,
\begin{equation}\label{equ:stomodelmknorm}
    \norm{f_x(y;s) - f(y;s)} \leq \frac{\eta}{2}\Norm{x-y}_{M_k}^2\quad \mbox{for all}\ x,y\in V.
\end{equation}
For simplicity, we denote by
\begin{align*}
    \tlx_{k+1}\triangleq \mr{prox}_{\alpha_k \olvphi_{x_k}\left(\cdot;S_k^{1:m}\right)}^{M_k}\left(x_k\right),\quad \olx_k \triangleq \prox_{\left(\olrho\right)^{-1}\phi}^{M_{k-1}}\left(x_k\right)\quad \mbox{and}\quad \mephi{k}(\cdot) \triangleq e_{\left(\olrho\right)^{-1}\phi}^{M_{k-1}}(\cdot)
\end{align*}
for each $k\in \mb{Z}_{+}$. The subsequent proof makes use of the following technical lemmas, the proof details of which are deferred to Appendix \ref{apx:proofstabilityasconvergencelemmaused}. 
\begin{lema}\label{lema:lipoff}
    The following assertions hold:
    \begin{enumerate}
        \item Let Assumption \ref{asm:stomodel} and \ref{asm:sppaall} hold. Then the function $\phi$ is proper, closed and $(\oleta + \oltau)$-weakly convex.
        \item Let Assumption \ref{asm:stomodel}, \ref{asm:sppaall} and \ref{asm:sppa1} hold. Fix any bounded open convex subset $U\subset V$. Then, we have 
        \begin{align*}
            \vert F(x) - F(y)\vert \leq\olL_{F}\left(U\right)\Norm{x - y}_2\quad \mbox{for all}\ x,y\in U.
        \end{align*} 
        \item Let Assumption \ref{asm:stomodel}, \ref{asm:sppaall} and \ref{asm:generalizedlip2} hold. Then, we have 
        \begin{align*}
            \vert F(x) - F(y)\vert \leq\olL_{F}\glip\left(\Norm{x}_2\right)\Norm{x - y}_2 + \frac{\oleta}{2}\Norm{x-y}_2^2\quad \mbox{for all}\ x,y\in V.
        \end{align*} 
        \item Let Assumption \ref{asm:stomodel}, \ref{asm:sppaall} and \ref{asm:generalizedlip3} hold. Fix any $k\in \mb{Z}_{+}$ and $x_k\in V$. Then, we have 
        \begin{align*}
            \Norm{u}_2 \leq \left(\frac{1}{m}\sum_{i=1}^m \olL_\varphi\left(S_k^i\right)\right) \glip\left(\max\lrbrackets{\Norm{x_k}_2,\Norm{x}_2}\right)\quad \mbox{for all}\ u\in \partial \olvphi_{x_k}\left(x;S_k^{1:m}\right)
        \end{align*} 
        holds for all $x\in V$.
    \end{enumerate}
\end{lema} 
\begin{lema}\label{lema:onestepprogress}
    Let Assumption \ref{asm:mk} and \ref{asm:sppaall} hold, and let the iterates $\{x_k\}$ be generated by the model-based ispPPA (Algorithm \ref{algo:ispppa}) with stepsizes $\{\alpha_k\}$ and parameters $\{\epsilon_k\}$. Denote by $\tau \triangleq L_\infty^2\cdot\oltau$. Suppose that $\sup_{k\in \mb{Z}_{+}}\alpha_k < \tau^{-1}$. Then for all $x\in \mb{R}^d$ and $k\in \mb{Z}_{+}$, we have $\tlx_{k+1}\triangleq \mr{prox}_{\alpha_k \olvphi_{x_k}(\cdot;S_k^{1:m})}^{M_k}(x_k)$ is well-defined, which satisfies
    \begin{equation}\label{equ:cauchyschwarz}
        \begin{split}
            \left\Vert x_{k+1} - x\right\Vert_{M_k}^2 \leq \left(1 + \frac{1}{t}\right)\left\Vert\tilde{x}_{k+1} - x \right\Vert_{M_k}^2 + \left(1+t\right)\epsilon_k^2\quad \mbox{for all}\ t\in \mb{R}_{++},
        \end{split}
    \end{equation}
    and
    \begin{equation}\label{inequ:onestep}
        \begin{split}
            &\olvphi_{x_k}\left(x;S_k^{1:m}\right) + \frac{1}{2\alpha_k}\Norm{x-x_k}_{M_k}^2\\
            \geq\ &\olvphi_{x_k}\left(\tilde{x}_{k+1};S_k^{1:m}\right) + \frac{1}{2\alpha_k}\Norm{\tilde{x}_{k+1} - x_k}_{M_k}^2 + \frac{1-\tau \alpha_k}{2\alpha_k} \Norm{x - \tilde{x}_{k+1}}_{M_k}^2.
        \end{split}
    \end{equation}
\end{lema}
\begin{lema}\label{lema:welldefinedolxk}
    Let Assumption \ref{asm:stomodel}, \ref{asm:mk} and \ref{asm:sppaall} hold, and set $k\in \mb{Z}_{+}$. Then for $\olrho > \eta + \tau$, the following assertions hold:
    \begin{enumerate}
        \item The proximal mapping $\prox_{(\olrho)^{-1}\phi}^{M_{k-1}}(\cdot)$ is single-valued on $\mb{R}^d$, i.e., for each $x\in \mb{R}^d$, there exists a unique point $\olx\in \mb{R}^d$ such that $\olx = \arg\min_{y\in \mb{R}^d}\{ \phi(y) + \frac{\olrho}{2}\Vert y - x\Vert_{M_{k-1}}^2\}$. 
        \item The Moreau envelope $e_{(\olrho)^{-1}\phi}^{M_{k-1}}(\cdot)$ is differentiable on $\mb{R}^d$ with gradient given by $\nabla e_{(\olrho)^{-1}\phi}^{M_{k-1}}\left(x\right) = \olrho M_{k-1}(x - \prox_{(\olrho)^{-1}\phi}^{M_{k-1}}(x))$ for each $x\in \mb{R}^d$, which satisfies
        \begin{equation}\label{equ:gradnormmkme}
            \Norm{\nabla e_{\left(\olrho\right)^{-1}\phi}^{M_{k-1}}\left(x\right)}_{M_{k-1}^{-1}} = \Norm{\olrho \left(x - \prox_{\left(\olrho\right)^{-1}\phi}^{M_{k-1}}\left(x\right)\right)}_{M_{k-1}}.
        \end{equation}
    \end{enumerate}
\end{lema}
\begin{lema}\label{lema:inequalitydistance}
    Let Assumption \ref{asm:stomodel}, \ref{asm:mk} and \ref{asm:sppaall} hold, and let $\{x_k\}$ be generated by the model-based ispPPA (Algorithm \ref{algo:ispppa}) with stepsizes $\{\alpha_k\}$ and parameters $\{\epsilon_k\}$. Denote by $\eta \triangleq L_\infty^2\cdot\oleta$ and $\tau \triangleq L_\infty^2\cdot\oltau$. Set $B_\rho \triangleq \sup_{k\in \mb{Z}_{+}} \rho_k^2$ and let $t_k \triangleq \sqrt{\mb{E}_k[\Vert\tlx_{k+1}-x_k\Vert_{M_k}^2]}$ for all $k\in \mb{Z}_{+}$. Then the following assertions hold:
    \begin{enumerate}
        \item Let Assumption \ref{asm:sppa1} hold and assume that $\sup_{k\in \mb{Z}_{+}}\alpha_k < \tau^{-1}$. Suppose that $\{x_k\}$ and $\{\tlx_k\}$ are contained in some bounded open convex subset $U$ of $V$. Then for all $x\in \mb{R}^d$ and $k\in \mb{Z}_{+}$, we have 
        \begin{subequations}
            \begin{align}
                &\frac{1-\tau \alpha_k}{2\alpha_k}\mb{E}_k\left[\Vert \tilde{x}_{k+1} - x\Vert_{M_k}^2\right] - \frac{1+\eta \alpha_k}{2\alpha_k}\mb{E}_k\left[\Vert x_k - x\Vert_{M_k}^2\right]\nonumber\\ \leq\ & - \mb{E}_k\left[\phi\left(\tilde{x}_{k+1}\right) - \phi\left(x\right)\right] - \frac{1}{2\alpha_k}t_k^2 + \sqrt{\rho_{f,m,U}}L_F\left(U\right) t_k\label{equ:inequ1}\\
                \leq\ &- \mb{E}_k\left[\phi\left(\tilde{x}_{k+1}\right) - \phi\left(x\right)\right] + \frac{\rho_{f,m,U}L_F\left(U\right)^2}{2}\alpha_k.\label{equ:inequ1max}
            \end{align}
        \end{subequations}
        Denote by $b = 4(\olrho\rho_k^{-2} - \eta)(\olrho\rho_k^{-2}-\eta-\tau)$. Furthermore, if $\sup_{k\in \mb{Z}_{+}}\alpha_k < (\olrho\max\{2,B_\rho\})^{-1}$ for some $\olrho > B_\rho \cdot(\eta + \tau)$ and $\phi$ is bounded from below, it holds that 
        \begin{equation}\label{equ:mephiinequ1}
            \begin{split}
                &\mb{E}_k\left[\nu_{k+1}^{-2}\cdot\left(\mephi{k+1}\left(x_{k+1}\right) - \phi^{*}\right)\right] - \nu_k^{-2}\cdot\left(\mephi{k}\left(x_k\right) - \phi^{*}\right)\\ 
                \leq\ & \frac{\olrho\left(1+b\alpha_k^2\right)}{2}\left[\frac{\rho_{f,m,U}L_F\left(U\right)^2}{\left(1-\tau \alpha_k\right)\left(1-\olrho\rho_k^2\alpha_k\right)} + \frac{\epsilon_k^2}{b\alpha_k^4}\right]\alpha_k^2\\
                &\hspace{0em} - \frac{\left(\olrho\rho_k^{-2}-\eta-\tau\right)\left[1-2\left(\olrho\rho_k^{-2}-\eta\right)\alpha_k\right]^2}{2\olrho\left(1-\tau \alpha_k\right)}\alpha_k \cdot \nu_k^{-2}\Norm{\olrho \left(x_k - \olx_k\right)}_{M_{k-1}}^2,
            \end{split}
        \end{equation}
        where $\rho_{f,m,U} \triangleq (\sqrt{\frac{1+(m-1)\omega_{f,U}}{m}}+1)^2$ for $\omega_{f,U} \triangleq 1 - \frac{\mr{Var}(L_{f,U}(s))}{\mb{E}_{s\sim P}[L_{f,U}(s)^2]}$\footnote{Without loss of generality, we may assume that $\mathbb{E}_{s\sim P}[L_{f,U}(s)^2]>0$.} and $L_{f,U}(s) \triangleq L_\infty\cdot \olL_{f,U}(s)$, and $L_F(U) \triangleq L_\infty \cdot\olL_F(U)$.
        \item Let Assumption \ref{asm:generalizedlip2} hold and assume that $\sup_{k\in \mb{Z}_{+}}\alpha_k < (\olrho\max\{2,B_\rho\}+\eta)^{-1}$ for $\olrho > B_\rho \cdot(\eta + \tau)$. Denote by $b = 4(\olrho\rho_k^{-2} - \eta)(\olrho\rho_k^{-2}-\eta-\tau)$. Then for all $x\in \mb{R}^d$ and $k\in \mb{Z}_{+}$, we have 
        \begin{equation}\label{equ:inequ2}
            \begin{split}
                &\frac{1-\tau \alpha_k}{2\alpha_k}\mb{E}_k\left[\Vert \tilde{x}_{k+1} - x\Vert_{M_k}^2\right] - \frac{1+\eta \alpha_k}{2\alpha_k}\mb{E}_k\left[\Vert x_k - x\Vert_{M_k}^2\right]\\ \leq\ & - \mb{E}_k\left[\phi\left(\tilde{x}_{k+1}\right) - \phi\left(x\right)\right] - \frac{1-\eta\alpha_k}{2\alpha_k}t_k^2 + \sqrt{\rho_{f,m}}L_F\glip\left(\Norm{x_k}_2\right) t_k
            \end{split}
        \end{equation}
        and when assuming further that $\phi$ is bounded from below, it holds that 
        \begin{equation}\label{equ:mephiinequ2}
            \begin{split}
                &\mb{E}_k\left[\nu_{k+1}^{-2}\cdot\left(\mephi{k+1}\left(x_{k+1}\right) - \phi^{*}\right)\right] - \nu_k^{-2}\cdot\left(\mephi{k}\left(x_k\right) - \phi^{*}\right)\\ 
                \leq\ & \frac{\olrho\left(1+b\alpha_k^2\right)}{2}\left[\frac{\rho_{f,m}L_F^2\glip\left(\Norm{x_k}_2\right)^2}{\left(1-\tau \alpha_k\right)\left(1-\olrho\rho_k^2\alpha_k-\eta\alpha_k\right)} + \frac{\epsilon_k^2}{b\alpha_k^4}\right]\alpha_k^2\\
                &\hspace{0em} - \frac{\left(\olrho\rho_k^{-2}-\eta-\tau\right)\left[1-2\left(\olrho\rho_k^{-2}-\eta\right)\alpha_k\right]^2}{2\olrho\left(1-\tau \alpha_k\right)}\alpha_k \cdot \nu_k^{-2}\Norm{\olrho \left(x_k - \olx_k\right)}_{M_{k-1}}^2,
            \end{split}
        \end{equation}
        where $\rho_{f,m} \triangleq (\sqrt{\frac{1+(m-1)\omega_{f}}{m}}+1)^2$ for $\omega_{f} \triangleq 1 - \frac{\mr{Var}(L_{f}(s))}{\mb{E}_{s\sim P}[L_{f}(s)^2]}$ and $L_{f}(s) \triangleq L_\infty\cdot \olL_{f}(s)$, and $L_F \triangleq L_\infty \cdot\olL_F$.
        \item Let Assumption \ref{asm:generalizedlip3} hold and assume that $\sup_{k\in \mb{Z}_{+}}\alpha_k < (2\olrho)^{-1}$ for $\olrho > B_\rho \cdot(\eta + 2\tau)$. Denote by $b = 4(\olrho\rho_k^{-2} - \eta)(\olrho\rho_k^{-2}-\eta-2\tau)$. Then for all $x\in \mb{R}^d$ and $k\in \mb{Z}_{+}$, we have 
        \begin{equation}\label{equ:inequ3}
            \begin{split}
                &\frac{1-2\tau \alpha_k}{2\alpha_k}\mb{E}_k\left[\Vert \tilde{x}_{k+1} - x\Vert_{M_k}^2\right] - \frac{1+\eta \alpha_k}{2\alpha_k}\mb{E}_k\left[\Vert x_k - x\Vert_{M_k}^2\right]\\ \leq\ & - \mb{E}_k\left[\phi\left(x_k\right) - \phi\left(x\right)\right] - \frac{1}{2\alpha_k}t_k^2 + \sqrt{\rho_{\varphi,m}}L_\phi\glip\left(\Norm{x_k}_2\right) t_k
            \end{split}
        \end{equation}
        and when assuming further that $\phi$ is bounded from below, it holds that 
        \begin{equation}\label{equ:mephiinequ3}
            \begin{split}
                &\mb{E}_k\left[\nu_{k+1}^{-2}\cdot\left(\mephi{k+1}\left(x_{k+1}\right) - \phi^{*}\right)\right] -\nu_k^{-2}\cdot\left(\mephi{k}\left(x_k\right) - \phi^{*}\right)  \\ 
                \leq\ & \frac{\olrho\left(1+b\alpha_k^2\right)}{2}\left[\frac{\rho_{\varphi,m}L_\phi^2\glip\left(\Norm{x_k}_2\right)^2}{1-2\tau \alpha_k} + \frac{\epsilon_k^2}{b\alpha_k^4}\right]\alpha_k^2\\
                &\hspace{0em} - \frac{\left(\olrho\rho_k^{-2}-\eta-2\tau\right)\left[1-2\left(\olrho\rho_k^{-2}-\eta\right)\alpha_k\right]^2}{2\olrho\left(1-2\tau \alpha_k\right)}\alpha_k \cdot \nu_k^{-2}\Norm{\olrho \left(x_k - \olx_k\right)}_{M_{k-1}}^2,
            \end{split}
        \end{equation}
        where $\rho_{\varphi,m} \triangleq \frac{1+(m-1)\omega_{\varphi}}{m}$ for $\omega_{\varphi} \triangleq 1 - \frac{\mr{Var}(L_{\varphi}(s))}{\mb{E}_{s\sim P}[L_{\varphi}(s)^2]}$ and $L_{\varphi}(s) \triangleq L_\infty\cdot \olL_{\varphi}(s)$, and $L_\phi \triangleq L_\infty \cdot\olL_\phi$.
    \end{enumerate}
\end{lema}

\subsection{Proof of Claim \ref{claim:sippinequhatx1tvalueglip}}\label{apx:proofclaimsippinequhatx1tvalueglip}
Using Assumption \ref{asm:mk}, one can conclude that 
\begin{equation}\label{equ:rhonurelationapx}
    \nu_k \geq 1\quad\mbox{and}\quad 1\leq \rho_k\leq \rho_k\nu_k = \nu_{k+1} \leq \nu_\infty\quad \mbox{for all}\ k\in \mb{Z}_{+}.
\end{equation}
\subsubsection{Under the Setting of Theorem \ref{thm:almostsureboundednessadpstep1}}
Fix any $k\in \mb{Z}_{+}$. Recall that $B_\rho \triangleq \sup_{k\in \mb{Z}_{+}}\rho_k^2$. Condition $\olrho > B_\rho\cdot(\eta + \tau)$ implies $\olrho > \eta + \tau$ since $B_\rho\geq 1$. It then follows from Lemma \ref{lema:welldefinedolxk} that $\olx_k$ is well-defined. Recall that $\phi^{*}\triangleq \inf_{x\in \mb{R}^d}\phi(x)$. Under Assumption \ref{asm:sppaall} and \ref{asm:sublevelset}, $\phi$ is proper and bounded from below, implying $\phi^{*}\in \mb{R}$. Using the definition of the Moreau envelope, we deduce 
\begin{equation}\label{equ:lbmeapx}
    \mephi{k}\left(x_{k}\right) - \phi^{*} = \phi\left(\olx_k\right) + \frac{\olrho}{2}\Norm{\olx_k - x_k}_{M_{k-1}}^2 - \phi^{*} \geq \phi\left(\olx_k\right) - \phi^{*}\geq 0.
\end{equation} 
Use the same notations as in Lemma \ref{lema:inequalitydistance}(\romannumeral2). It follows from condition (\ref{asm:epsilonkalphakglip}) that
\begin{equation}\label{equ:alphakolalphakglip}
    \alpha_k \glip\left(\Norm{x_k}_2\right)\leq \olalpha_k\quad \mbox{and}\quad \alpha_k \leq \olalpha_k.
\end{equation}
Using (\ref{equ:alphakolalphakglip}) and applying (\ref{equ:mephiinequ2}) in Lemma \ref{lema:inequalitydistance}(\romannumeral2) with $\epsilon_k = \gamma \alpha_k^2$, we deduce
\begin{equation}\label{equ:sippinequhatx1tvalueglipbeta}
    \begin{split}
        &\mb{E}_k\left[\nu_{k+1}^{-2}\cdot\left(\mephi{k+1}\left(x_{k+1}\right) - \phi^{*}\right)\right]\\
        \leq\ &\nu_k^{-2}\cdot\left(\mephi{k}\left(x_k\right) - \phi^{*}\right) - \beta_k^{(1)} \alpha_k \cdot \nu_k^{-2}\Norm{\olrho \left(x_k - \olx_k\right)}_{M_{k-1}}^2 + \beta_k^{(2)} \olalpha_k^2,
    \end{split}
\end{equation}
where 
\begin{align*}
    \beta_k^{(1)} &\triangleq \frac{\left(\olrho\rho_k^{-2}-\eta-\tau\right)\left[1-2\left(\olrho\rho_k^{-2}-\eta\right)\alpha_k\right]^2}{2\olrho\left(1-\tau \alpha_k\right)}\quad \mbox{and}\\
    \beta_k^{(2)} &\triangleq \frac{\olrho\left(1+b\alpha_k^2\right)}{2}\left[\frac{\rho_{f,m}L_F^2}{\left(1-\tau \alpha_k\right)\left(1-\olrho\rho_k^2\alpha_k-\eta\alpha_k\right)} + \frac{\gamma^2}{b}\right]
\end{align*}
with $b \triangleq 4(\olrho\rho_k^{-2} - \eta)(\olrho\rho_k^{-2}-\eta-\tau)$. 
Let $\alpha_0 \triangleq \sup_{k\in \mb{Z}_{+}}\olalpha_k$, and set
\begin{equation}\label{equ:lbbetakglip}
    \begin{split}
        \theta_1 &\triangleq \frac{\left(\olrho B_\rho^{-1}-\eta-\tau\right)\left(1-2\olrho\alpha_0\right)^2}{2\olrho}\quad \mbox{and}\\
        \theta_2 &\triangleq \frac{\olrho\left[1+4\left(\olrho - \eta\right)\left(\olrho-\eta-\tau\right)\alpha_0^2\right]}{2}\\
        &\hspace{1em}\cdot \left[\frac{\rho_{f,m}L_F^2}{\left(1-\olrho \alpha_0\right)\left(1-\left(\olrho B_\rho+\eta\right)\alpha_0\right)} + \frac{\gamma^2}{4\left(\olrho B_\rho^{-1} - \eta\right)\left(\olrho B_\rho^{-1}-\eta-\tau\right)}\right].
    \end{split}
\end{equation}
Since $\olrho > B_\rho\cdot (\eta+\tau)$ and $\sup_{k\in \mb{Z}_{+}}\alpha_k \leq \sup_{k\in \mb{Z}_{+}}\olalpha_k < (\olrho\max\{2,B_\rho\}+\eta)^{-1}$, we can thus deduce from (\ref{equ:rhonurelationapx}) that 
\begin{equation}\label{equ:lbbetakglipbd}
    \begin{split}
        \beta_k^{(1)} &\geq \theta_1 > 0\quad \mbox{and}\quad
        0<\beta_k^{(2)} \leq \theta_2 < \infty.
    \end{split}
\end{equation}
It then follows from (\ref{equ:sippinequhatx1tvalueglipbeta})-(\ref{equ:lbbetakglipbd}) that the estimate (\ref{equ:medecreaseglip}) holds with positive scalars $\theta_1$ and $\theta_2$ defined in (\ref{equ:lbbetakglip}).

\subsubsection{Under the Setting of Theorem \ref{thm:almostsureboundednessadpstep2}}
Fix any $k\in \mb{Z}_{+}$. Condition $\olrho > B_\rho\cdot(\eta + 2\tau)$ implies $\olrho > \eta + 2\tau$ since $B_\rho\geq 1$ by (\ref{equ:rhonurelationapx}). It then follows from Lemma \ref{lema:welldefinedolxk} that $\olx_k$ is well-defined and thus $\mephi{k}(x_k) - \phi^{*}\geq 0$ (see the proof of (\ref{equ:lbmeapx})). Use the same notations as in Lemma \ref{lema:inequalitydistance}(\romannumeral3). By setting $\epsilon_k = \gamma \alpha_k^2$ in (\ref{equ:mephiinequ3}) from Lemma \ref{lema:inequalitydistance}(\romannumeral3), we have 
\begin{equation}\label{equ:sippinequhatx1tvaluenglip2}
    \begin{split}
        &\mb{E}_k\left[\nu_{k+1}^{-2}\cdot\left(\mephi{k+1}\left(x_{k+1}\right) - \phi^{*}\right)\right] - \nu_k^{-2}\cdot\left(\mephi{k}\left(x_k\right) - \phi^{*}\right)\\ 
        \leq\ & \frac{\olrho\left(1+b\alpha_k^2\right)}{2}\left[\frac{\rho_{\varphi,m}L_\phi^2\glip\left(\Norm{x_k}_2\right)^2}{1-2\tau \alpha_k} + \frac{\gamma^2}{b}\right]\alpha_k^2\\
        &\hspace{0em} - \frac{\left(\olrho\rho_k^{-2}-\eta-2\tau\right)\left[1-2\left(\olrho\rho_k^{-2}-\eta\right)\alpha_k\right]^2}{2\olrho\left(1-2\tau \alpha_k\right)}\alpha_k \cdot \nu_k^{-2}\Norm{\olrho \left(x_k - \olx_k\right)}_{M_{k-1}}^2,
    \end{split}
\end{equation}
where $b = 4(\olrho\rho_k^{-2} - \eta)(\olrho\rho_k^{-2}-\eta-2\tau)$. Let $\alpha_0 \triangleq \sup_{k\in \mb{Z}_{+}}\olalpha_k$, and set
\begin{equation}\label{equ:lbbetakglip2}
    \begin{split}
        \theta_1 &\triangleq \frac{\left(\olrho B_\rho^{-1}-\eta-2\tau\right)\left(1-2\olrho\alpha_0\right)^2}{2\olrho}\quad \mbox{and}\\
        \theta_2 &\triangleq \frac{\olrho\left[1+4\left(\olrho - \eta\right)\left(\olrho-\eta-2\tau\right)\alpha_0^2\right]}{2}\\
        &\hspace{2em}\cdot \left[\frac{\rho_{\varphi,m}L_\phi^2}{1-\olrho \alpha_0} + \frac{\gamma^2}{4\left(\olrho B_\rho^{-1} - \eta\right)\left(\olrho B_\rho^{-1}-\eta-2\tau\right)}\right].
    \end{split}
\end{equation}
Since $\olrho > B_\rho\cdot (\eta + 2\tau)$ and $\sup_{k\in \mb{Z}_{+}}\alpha_k \leq \sup_{k\in \mb{Z}_{+}}\olalpha_k < (2\olrho)^{-1}$, by following an argument nearly identical to that in the proof of Claim \ref{claim:sippinequhatx1tvalueglip}, we can deduce from (\ref{equ:sippinequhatx1tvaluenglip2}) that the estimate (\ref{equ:medecreaseglip}) holds with positive scalars $\theta_1$ and $\theta_2$ defined in (\ref{equ:lbbetakglip2}).

\subsection{Proof of Claim \ref{claim:sippstablebound}}\label{apx:proofclaimsippstablebound}
\underline{Proof of the estimate (\ref{equ:sippstablebound}).} Fix any $k\in \mb{Z}_{+}$ and $\olx\in \mc{X}^{*}$. Recall that $\oleta=\oltau=0$. Lemma \ref{lema:lipoff}(\romannumeral1) implies that $\phi$ is convex. It follows that $0\in \partial \phi(\olx)$ and thus, by the second condition in Assumption \ref{asm:sppbdoptset}, we can take $\varphi^\prime(\olx;S_k^i)\in \partial \varphi(\olx;S_k^i)$ such that 
\begin{equation}\label{equ:expksubgvarphi0}
    \mb{E}_k \left[\varphi^\prime(\olx;S_k^i)\right] = 0\quad \mbox{for}\ i = 1,\cdots,m.
\end{equation}
Denote by $\olvphi(\cdot;S_k^{1:m}) \triangleq \olvphi_{x_k}(\cdot;S_k^{1:m})$. Set $\olvphi^\prime(\olx;S_k^{1:m})\triangleq \frac{1}{m}\sum_{i=1}^m \varphi^\prime(\olx;S_k^i)$. Condition $\oleta = 0$ implies $\olvphi^\prime(\olx;S_k^{1:m})\in \partial \olvphi(\olx;S_k^{1:m})$ and $\mb{E}_k \left[\olvphi^\prime(\olx;S_k^{1:m})\right] = 0$. Using Assumption \ref{asm:sppaall} and condition $\oltau = 0$, we have $\olvphi(\cdot;S_k^{1:m})$ is convex and $\tlx_{k+1} \triangleq \prox_{\alpha_k \olvphi_{x_k}(\cdot;S_k^{1:m})}^{M_k}(x_k)$ is well-defined by Lemma \ref{lema:welldefinedolxk}, satisfying
\begin{equation}\label{equ:convexvphixstar}
    \begin{split}
        &\olvphi\left(\olx;S_k^{1:m}\right) - \olvphi\left(\tilde{x}_{k+1};S_k^{1:m}\right) \leq \lrangle{\olvphi^\prime\left(\olx;S_k^{1:m}\right), \olx - \tilde{x}_{k+1}}\\
        =\ &\lrangle{\olvphi^\prime\left(\olx;S_k^{1:m}\right), \olx - x_k} + \lrangle{\olvphi^\prime\left(\olx;S_k^{1:m}\right), x_k - \tilde{x}_{k+1}}\\
        \leq\ &\lrangle{\olvphi^\prime\left(\olx;S_k^{1:m}\right), \olx - x_k} + \frac{\alpha_kL_k^2}{2} \Norm{\olvphi^\prime\left(\olx;S_k^{1:m}\right)}_2^2 + \frac{1}{2\alpha_kL_k^2}\Norm{\tilde{x}_{k+1}-x_k}_2^2\\
        \leq\ &\lrangle{\olvphi^\prime\left(\olx;S_k^{1:m}\right), \olx - x_k} + \frac{\alpha_kL_k^2}{2} \Norm{\olvphi^\prime\left(\olx;S_k^{1:m}\right)}_2^2 + \frac{1}{2\alpha_k}\Norm{\tilde{x}_{k+1}-x_k}_{M_k}^2
    \end{split}
\end{equation}
where the second last inequality uses the Cauchy-Schwarz inequality, and the last is due to Assumption \ref{asm:mk}. Using (\ref{equ:convexvphixstar}) and applying (\ref{inequ:onestep}) in Lemma \ref{lema:onestepprogress} with $\eta =\tau= 0$ and $x = \olx$, we deduce
\begin{equation}\label{equ:onestepprogressxstar}
    \begin{split}
        &\Norm{\tilde{x}_{k+1} - \olx}_{M_k}^2\\ \leq\ &\Norm{x_k - \olx}_{M_k}^2 + 2\alpha_k \left[\lrangle{\olvphi^\prime\left(\olx;S_k^{1:m}\right), \olx - x_k} + \frac{\alpha_kL_k^2}{2} \Norm{\olvphi^\prime\left(\olx;S_k^{1:m}\right)}_2^2\right.\\ &\hspace{8.5em}\left.+ \frac{1}{2\alpha_k}\Norm{\tilde{x}_{k+1}-x_k}_{M_k}^2\right] - \Norm{\tilde{x}_{k+1} - x_k}_{M_k}^2\\
        =\ &\Norm{x_k - \olx}_{M_k}^2 + 2\alpha_k \lrangle{\olvphi^\prime\left(\olx;S_k^{1:m}\right), \olx - x_k} + \alpha_k^2 L_k^2 \Norm{\olvphi^\prime\left(\olx;S_k^{1:m}\right)}_2^2.
    \end{split}
\end{equation}
It follows from the definition of $\olvphi^\prime\left(\olx;S_k^{1:m}\right)$, the independence of the random variables and condition (\ref{equ:expksubgvarphi0}) that 
\begin{equation}\label{equ:linearspeedupbatchsize}
    \begin{split}
        \mb{E}_k \left[\Norm{\olvphi^\prime\left(\olx;S_k^{1:m}\right)}_2^2\right] 
        &= \frac{1}{m^2} \sum_{i=1}^m \mb{E}_k \left[\Norm{\varphi^\prime\left(\olx;S_k^i\right)}_2^2\right].
    \end{split}
\end{equation}
Applying (\ref{equ:cauchyschwarz}) in Lemma \ref{lema:onestepprogress} with $x = \olx$, $t = \frac{1}{\alpha_k^2}$ and $\epsilon_k = \gamma \alpha_k^2$, we have 
\begin{align*}
        \Norm{x_{k+1} - \olx}_{M_k}^2 \leq\ &\left(1+\alpha_k^2\right) \Norm{\tilde{x}_{k+1}-\olx}_{M_k}^2 + \left(1+\alpha_k^2\right)\alpha_k^2 \gamma^2\\
        \leq\ &\left(1+\alpha_k^2\right) \left[\Norm{x_k - \olx}_{M_k}^2 + 2\alpha_k \lrangle{\olvphi^\prime\left(\olx;S_k^{1:m}\right), \olx - x_k}\right.\\ &\hspace{9.5em}\left.+\ \alpha_k^2L_k^2 \Norm{\olvphi^\prime\left(\olx;S_k^{1:m}\right)}_2^2\right] + \left(1+\alpha_k^2\right)\alpha_k^2 \gamma^2,
\end{align*}
where the second inequality follows from (\ref{equ:onestepprogressxstar}). Taking expectations on both sides of the previous inequality with respect to $\mc{F}_{k-1}$ yields that 
\begin{align}
        &\mb{E}_k \left[\Norm{x_{k+1} - \olx}_{M_k}^2\right]\nonumber\\ \leq\ &\left(1+\alpha_k^2\right)\mb{E}_k \left[\Norm{x_k - \olx}_{M_k}^2\right] + \left(1+\alpha_k^2\right)\alpha_k^2 \left(L_k^2\mb{E}_k \left[\Norm{\olvphi^\prime\left(\olx;S_k^{1:m}\right)}_2^2\right] + \gamma^2\right)\nonumber\\
        =\ &\left(1+\alpha_k^2\right)\mb{E}_k \left[\Norm{x_k - \olx}_{M_k}^2\right] + \left(1+\alpha_k^2\right)\alpha_k^2 \left(\frac{L_k^2}{m^2} \sum_{i=1}^m \mb{E}_k \left[\Norm{\varphi^\prime\left(\olx;S_k^i\right)}_2^2\right] + \gamma^2\right)\nonumber\\
        \leq\ &\left(1+\alpha_k^2\right)\rho_k^2\Norm{x_k - \olx}_{M_{k-1}}^2 + \left(1+\alpha_k^2\right)\alpha_k^2 \left(\frac{L_k^2\sigma_\phi^2}{m} + \gamma^2\right),\label{equ:sippstableboundoldrho}
\end{align}
where the first inequality comes from the condition that $\mb{E}_k \left[\olvphi^\prime(\olx;S_k^{1:m})\right] = 0$ and the last is due to Assumption \ref{asm:mk} and \ref{asm:sppbdoptset}, and the equality holds because of (\ref{equ:linearspeedupbatchsize}). By multiplying both sides of (\ref{equ:sippstableboundoldrho}) by $\nu_{k+1}^{-2}$, we deduce from $\nu_{k+1} = \rho_k \nu_k$ and $\nu_{k+1} \geq 1$ that
\begin{align*}
        &\mb{E}_k \left[\nu_{k+1}^{-2}\Norm{x_{k+1} - \olx}_{M_k}^2\right]\\ 
        \leq\ &\nu_k^{-2}\left(1+\alpha_k^2\right)\Norm{x_k - \olx}_{M_{k-1}}^2 + \left(1+\alpha_k^2\right)\alpha_k^2 \left(\frac{L_k^2\sigma_\phi^2}{m} + \gamma^2\right),
\end{align*}
which completes the proof of (\ref{equ:sippstablebound}). 

\underline{Proof of the estimate (\ref{equ:sippinequhatx1tvalue}).} Now consider the case where Assumption \ref{asm:sppa1} holds and both $\{x_k\}$ and $\{\tlx_k\}$ are contained in some bounded open convex subset $U$ of $V$. Fix any $k\in \mb{Z}_{+}$ and $x\in \mb{R}^d$. Recall that $\eta = L_\infty^2\cdot \oleta$ and $\tau = L_\infty^2\cdot \oltau$. Note that all assumptions required in Lemma \ref{lema:inequalitydistance}(\romannumeral1) are satisfied. Use the same notations as in Lemma \ref{lema:inequalitydistance}(\romannumeral1). Then, by setting $\oleta = \oltau = 0$ in (\ref{equ:inequ1max}) of Lemma \ref{lema:inequalitydistance}(\romannumeral1), applying (\ref{equ:cauchyschwarz}) in Lemma \ref{lema:onestepprogress} with $t = \frac{1}{\alpha_k^2}$ and $\epsilon_k = \gamma \alpha_k^2$, and using $\mb{E}_k[\Vert x_k-x\Vert_{M_k}^2] \leq \rho_k^2 \Vert x_k - x\Vert_{M_{k-1}}^2$ by Assumption \ref{asm:mk}, we can conclude that 
\begin{align}
        \mb{E}_k\left[\Norm{x_{k+1} - x}_{M_k}^2\right] \leq\ 
        &\left(1+\alpha_k^2\right)\cdot\left(\rho_k^2\Vert x_k - x\Vert_{M_{k-1}}^2 - 2\alpha_k\cdot \mb{E}_k\left[\phi\left(\tilde{x}_{k+1}\right) - \phi\left(x\right)\right]\right.\nonumber\\
        &\hspace{9em}+\ \left.\left(\rho_{f,m,U}L_F(U)^2 + \gamma^2\right)\cdot\alpha_k^2\right).\label{equ:sippinequhatx1tvalueoldrho}
\end{align}
Note that $\nu_{k+1} = \rho_k \nu_k$ and $\nu_{k+1} \geq 1$. Then, the estimate (\ref{equ:sippinequhatx1tvalue}) follows by multiplying both sides of (\ref{equ:sippinequhatx1tvalueoldrho}) by $\nu_{k+1}^{-2}$.

\subsection{Proof of Auxiliary Lemmas}\label{apx:proofstabilityasconvergencelemmaused}

\subsubsection{Proof of Lemma \ref{lema:lipoff}} 
By following the proof arguments of \cite[Lemma~4.1]{davis2019stochastic} with suitable modifications, we can prove Lemma \ref{lema:lipoff}. 
\begin{enumerate}
    \item Following an argument nearly identical to that in the proof of the first part of  \cite[Lemma~4.1]{davis2019stochastic}, and using Assumption \ref{asm:stomodel} and \ref{asm:sppaall}, we can show that assertion (\romannumeral1) holds. 
    \item Let Assumption \ref{asm:stomodel}, \ref{asm:sppaall} and \ref{asm:sppa1} hold. Fix any bounded open convex subset $U\subset V$ and let $x\in U$. Then, Assumption \ref{asm:stomodel} combined with Assumption \ref{asm:sppa1} implies 
    \begin{align*}
            F\left(y\right) - F\left(x\right) &
            = \mb{E}_{s\sim P}\left[f\left(y;s\right) - f_x\left(y;s\right) + f_x\left(y;s\right)- f_x\left(x;s\right)\right]\\
            &\leq \mb{E}_{s\sim P}\left[\frac{\oleta}{2}\Norm{x-y}_2^2 + \olL_{f,U}\left(s\right)\Norm{x-y}_2\right]\\
            &\leq \sqrt{\mb{E}_{s\sim P}\left[\olL_{f,U}\left(s\right)^2\right]}\Norm{x-y}_2 + \frac{\oleta}{2}\Norm{x-y}_2^2\\
            &\leq \olL_F\left(U\right)\Norm{x-y}_2 + \frac{\oleta}{2}\Norm{x-y}_2^2\quad \mbox{for all}\ y\in U.
    \end{align*}
    It follows from the preceding inequality that 
    \begin{equation}\label{equ:glipconstant1}
        \begin{split}
            \lim\sup_{y\neq x\to x} \frac{F\left(y\right) - F\left(x\right)}{\Norm{y-x}_2} &\leq \lim\sup_{y\neq x\to x}\frac{\olL_F\left(U\right)\Norm{x-y}_2 + \frac{\oleta}{2}\Norm{x-y}_2^2}{\Norm{y-x}_2} \\
            &= \olL_F\left(U\right)\quad \mbox{for all}\ x\in U.
        \end{split}
    \end{equation}
    Let $U_F \triangleq \{x\in U\ \vert\ F\ \mbox{is differentiable at}\ x\in U\}$. Fix any $x\in U_F$. By setting $y = x + t\nabla F(x)$ with $t\searrow 0$ in (\ref{equ:glipconstant1}) yields 
    \begin{equation}\label{equ:glipconstantdiff}
        \Norm{\nabla F\left(x\right)}_2 = \lim_{t\to 0} \frac{F\left(x + t\nabla F(x)\right) - F\left(x\right)}{\Norm{x + t\nabla F(x) - x}_2} \leq \olL_F\left(U\right)\quad \mbox{for all}\ x\in U_F.
    \end{equation}
    Fix any $x\in U$. Let $\partial_B F\left(x\right)$ be the B-subdifferential of $F$ at $x$. Since $F$ is locally Lipschitz by Assumption \ref{asm:sppaall}, according to the definition of the B-subdifferential, for any $u\in \partial_B F(x)$, there exists sequence $\{x^\nu\}\subset U_F$ such that $\lim_{\nu \to \infty}x^\nu = x$ and $\lim_{\nu \to \infty} \nabla F(x^\nu) = u$, which combined with (\ref{equ:glipconstantdiff}) implies 
    \begin{equation}\label{equ:glipconstantbsubdiff}
        \Norm{u}_2 = \lim_{\nu \to \infty}\Norm{\nabla F\left(x^\nu\right)} \leq \olL_F\left(U\right)\quad \mbox{for all}\ u\in \partial_B F(x)\ \mbox{with}\ x\in U.
    \end{equation}
    Combining (\ref{equ:glipconstantbsubdiff}) with the definition of the Clarke subdifferential yields
    \begin{equation}\label{equ:glipconstantbcsubdiff}
        \Norm{v}_2 \leq \olL_F\left(U\right)\quad \mbox{for all}\ v\in \partial_C F(x)\ \mbox{with}\ x\in U.
    \end{equation}
    Fix any $x,y\in U$. Given that $U$ is open and convex, we can derive from \cite[Theorem~2.4]{clarke1998nonsmooth} that there exists $\theta_1\in (0,1)$ such that 
    \begin{align*}
        F\left(y\right) - F\left(x\right)\in \lrbrackets{\lrangle{v,y-x}\ \vert\ v\in \partial_C F\left(z_{\theta_1}\right)}\quad \mbox{where}\ z_{\theta_1} = \theta_1 x + \left(1-\theta_1\right)y.
    \end{align*}
    Since $z_{\theta_1}\in U$, it follows from (\ref{equ:glipconstantbcsubdiff}) combined with the Cauchy-Schwarz inequality that 
    \begin{align*}
        \lrangle{v,y-x}\leq \Norm{v}_2 \Norm{x-y}_2 \leq \olL_F\left(U\right)\Norm{x-y}_2\quad \mbox{for all}\ v\in \partial_C F\left(z_{\theta_1}\right),
    \end{align*}
    and thus 
    \begin{equation}\label{equ:locallip1}
        F\left(y\right) - F\left(x\right)\leq \olL_F\left(U\right)\Norm{x-y}_2.
    \end{equation}
    Similarly, there exists $\theta_2\in (0,1)$ such that 
    \begin{align*}
        F\left(x\right) - F\left(y\right)\in \lrbrackets{\lrangle{v,x-y}\ \vert\ v\in \partial_C F\left(z_{\theta_2}\right)}\quad \ \mbox{where}\ z_{\theta_2} = \theta_2 x + \left(1-\theta_2\right)y
    \end{align*}
    and thus 
    \begin{equation}\label{equ:locallip2}
        F\left(x\right) - F\left(y\right)\leq \olL_F\left(U\right)\Norm{x-y}_2.
    \end{equation}
    Combining (\ref{equ:locallip1}) with (\ref{equ:locallip2}) proves assertion (\romannumeral2).
    \item Let Assumption \ref{asm:stomodel}, \ref{asm:sppaall} and \ref{asm:generalizedlip2} hold. Fix any $x,y\in V$. Then, Assumption \ref{asm:stomodel} combined with Assumption \ref{asm:generalizedlip2} implies 
    \begin{align*}
            \norm{F\left(y\right) - F\left(x\right)} 
            &\leq \mb{E}_{s\sim P}\left[\norm{f\left(y;s\right) - f_x\left(y;s\right)} + \norm{f_x\left(y;s\right)- f_x\left(x;s\right)}\right]\\
            &\leq \mb{E}_{s\sim P}\left[\frac{\oleta}{2}\Norm{x-y}_2^2 + \olL_f\left(s\right)\glip\left(\Norm{x}_2\right)\Norm{x-y}_2\right]\\
            &\leq \sqrt{\mb{E}_{s\sim P}\left[\olL_f\left(s\right)^2\right]}\glip\left(\Norm{x}_2\right)\Norm{x-y}_2 + \frac{\oleta}{2}\Norm{x-y}_2^2\\
            &\leq \olL_F \glip\left(\Norm{x}_2\right)\Norm{x-y}_2 + \frac{\oleta}{2}\Norm{x-y}_2^2.
    \end{align*}
    This completes the proof of assertion (\romannumeral3).
    \item Let Assumption \ref{asm:stomodel}, \ref{asm:sppaall} and \ref{asm:generalizedlip3} hold. Let $k\in \mb{Z}_{+}$ and $x_k\in V$. Fix any $x\in V$. It follows from Assumption \ref{asm:generalizedlip3} that
    \begin{equation}\label{equ:lipgeneralized2}
        \begin{split}
            &\olvphi_{x_k} \left(y;S_k^{1:m}\right) - \olvphi_{x_k} \left(x;S_k^{1:m}\right)=\ \frac{1}{m}\sum_{i=1}^m \left(\varphi_{x_k} \left(y;S_k^i\right) - \varphi_{x_k} \left(x;S_k^i\right)\right)\\
            \leq\ &\left(\frac{1}{m}\sum_{i=1}^m \olL_\varphi\left(S_k^i\right)\right)\glip\left(\max\lrbrackets{\Norm{x_k}_2,\Norm{x}_2}\right) \cdot \Norm{x - y}_2
            \\
            &\mbox{for all}\ y\ \mbox{in a neighborhood of}\ x.
        \end{split}
    \end{equation}
    Let $V_{\olvphi_{x_k}} \triangleq \{x\in V\ \vert\ \olvphi_{x_k}(\cdot;S_k^{1:m})\ \mbox{is differentiable at}\ x\in V\}$. Using (\ref{equ:lipgeneralized2}) and following the proof argument of (\ref{equ:glipconstantdiff}), we deduce 
    \begin{equation}\label{equ:glip3constantdiff}
        \begin{split}
            &\Norm{\nabla \olvphi_{x_k}\left(x;S_k^{1:m}\right)}_2 
            \leq \left(\frac{1}{m}\sum_{i=1}^m \olL_\varphi\left(S_k^i\right)\right)\glip\left(\max\lrbrackets{\Norm{x_k}_2,\Norm{x}_2}\right)
        \end{split}
    \end{equation}
    holds for all $x\in V_{\olvphi_{x_k}}$. 
    Fix any $x\in V$. Let $\partial_B \olvphi_{x_k}(x;S_k^{1:m})$ be the B-subdifferential of $\olvphi_{x_k}(\cdot;S_k^{1:m})$ at $x$. Since $\olvphi_{x_k}(\cdot;S_k^{1:m})$ is locally Lipschitz by Assumption \ref{asm:generalizedlip3}, by the definition of the B-subdifferential, for any $u\in \partial_B \olvphi_{x_k}(x;S_k^{1:m})$, there exists sequence $\{x^\nu\}\subset V_{\olvphi_{x_k}}$ such that $\lim_{\nu \to \infty}x^\nu = x$ and $\lim_{\nu \to \infty} \nabla \olvphi_{x_k}(x^\nu;S_k^{1:m}) = u$, which combined with (\ref{equ:glip3constantdiff}) implies 
    \begin{equation}\label{equ:glip3constantbsubdiff}
        \begin{split}
            \Norm{u}_2 
            &\leq\lim_{\nu \to \infty}\left(\frac{1}{m}\sum_{i=1}^m \olL_\varphi\left(S_k^i\right)\right)\glip\left(\max\lrbrackets{\Norm{x_k}_2,\Norm{x^\nu}_2}\right)\\
            &= \left(\frac{1}{m}\sum_{i=1}^m \olL_\varphi\left(S_k^i\right)\right)\glip\left(\max\lrbrackets{\Norm{x_k}_2,\Norm{x}_2}\right)\\
            &\hspace{10em}\quad \mbox{for all}\ u\in \partial_B \olvphi_{x_k}\left(x;S_k^{1:m}\right)\ \mbox{with}\ x\in V,
        \end{split}
    \end{equation}
    where the last equality uses the condition that $\glip(\cdot)$ is continuous over $\mb{R}_{+}$. Using the definition of the Clarke subdifferential, we can deduce from (\ref{equ:glip3constantbsubdiff}) that 
    \begin{equation}\label{equ:glip3constantbcsubdiff}
        \begin{split}
            \Norm{v}_2 &\leq \left(\frac{1}{m}\sum_{i=1}^m \olL_\varphi\left(S_k^i\right)\right)\glip\left(\max\lrbrackets{\Norm{x_k}_2,\Norm{x}_2}\right)\\
            &\hspace{10em}\quad \mbox{for all}\ v\in \partial_C \olvphi_{x_k}\left(x;S_k^{1:m}\right)\ \mbox{with}\ x\in V.
        \end{split}
    \end{equation}
    Under Assumption \ref{asm:sppaall}, $\olvphi_{x_k}(\cdot;S_k^{1:m})$ is weakly convex. Then, it follows from \cite[Proposition~4.4.15]{cui2021modern} that $\olvphi_{x_k}(\cdot;S_k^{1:m})$ is Clarke regular. That is, $\partial \olvphi_{x_k}(x;S_k^{1:m}) = \partial_C \olvphi_{x_k}(x;S_k^{1:m})$, which combined with (\ref{equ:glip3constantbcsubdiff}) proves assertion (\romannumeral4).
\end{enumerate}

\subsubsection{Proof of Lemma \ref{lema:onestepprogress}} 
Fix any $k\in \mb{Z}_{+}$. It follows from Assumption \ref{asm:sppaall} that function $\varphi_k(\cdot)\triangleq \olvphi_{x_k}(\cdot;S_k^{1:m})$ is proper, closed and $\oltau$-weakly convex. For any self-adjoint positive semidefinite linear mapping $M$, we have 
\begin{equation}\label{equ:abnormm}
    \lrangle{a,b}_M = \frac{1}{2}\left(\Norm{a}_M^2 + \Norm{b}_M^2 - \Norm{a-b}_M^2\right).
\end{equation}
Combining (\ref{equ:abnormm}) with Assumption \ref{asm:mk} yields
\begin{equation}\label{equ:scvxsqnormm}
    \begin{split}
        &\Norm{ \theta x + \left(1-\theta\right) y - x_k}_{M_k}^2\\
        =\ &\theta \Norm{x-x_k}_{M_k}^2 + \left(1-\theta\right)\Norm{y-x_k}_{M_k}^2 - \theta \left(1-\theta\right)\Norm{x-y}_{M_k}^2\\
        \leq\ &\theta \Norm{x-x_k}_{M_k}^2 + \left(1-\theta\right)\Norm{y-x_k}_{M_k}^2 - \frac{\theta \left(1-\theta\right)}{L_k^2}\Norm{x-y}_2^2\\ &\mbox{for all}\ x,y\in \mb{R}^d\ \mbox{and}\ \theta\in [0,1],
    \end{split}
\end{equation}
which implies $\psi_k(\cdot)\triangleq \frac{1}{2\alpha_k}\Vert \cdot - x_k\Vert_{M_k}^2$ is $\frac{1}{L_k^2\alpha_k}$-strongly convex. Given that $\tau = L_\infty^2\cdot \oltau$ and $\sup_{k\in \mb{Z}_{+}}\alpha_k < \tau^{-1}$, we can deduce from Assumption \ref{asm:mk} that
\begin{align*}
    \frac{1}{L_k^2 \alpha_k} - \oltau = \frac{1-L_k^2\alpha_k \oltau}{L_k^2 \alpha_k} \geq \frac{1-L_\infty^2\alpha_k \oltau}{L_k^2 \alpha_k} = \frac{1-\alpha_k \tau}{L_k^2 \alpha_k} > 0.
\end{align*}
Therefore, $\Phi_k(\cdot)\triangleq \varphi_k(\cdot) + \psi_k(\cdot)$ is proper, closed and $(\frac{1}{L_k^2 \alpha_k} - \oltau)$-strongly convex. Applying \cite[Lemma~2.26]{planiden2016strongly}, we can conclude that the minimizer of $\Phi_k$ on $\mb{R}^d$ exists and is unique. It then follows from the definition of the proximal mapping that $\tlx_{k+1} \triangleq \prox_{\alpha_k \olvphi_{x_k}(\cdot;S_k^{1:m})}^{M_k}(x_k) = \arg\min_{x\in \mb{R}^d} \Phi_k(x)$ is well-defined. Fix any $t\in \mb{R}_{++}$. Combining the Cauchy-Schwarz inequality with criterion (\ref{equ:criteriaa}) yields
\begin{align*}
    \left\Vert x_{k+1} - x\right\Vert_{M_k}^2 &= \left\Vert x_{k+1} - \tilde{x}_{k+1}\right\Vert_{M_k}^2 + \left\Vert\tilde{x}_{k+1} - x\right\Vert_{M_k}^2+ 2\langle x_{k+1} - \tilde{x}_{k+1}, \tilde{x}_{k+1} - x\rangle_{M_k}\\
        &\leq \left(1+t\right)\left\Vert x_{k+1} - \tilde{x}_{k+1}\right\Vert_{M_k}^2 + \left(1 + \frac{1}{t}\right)\left\Vert\tilde{x}_{k+1} - x \right\Vert_{M_k}^2\\
        &\leq \left(1 + \frac{1}{t}\right)\left\Vert\tilde{x}_{k+1} - x \right\Vert_{M_k}^2 + \left(1+t\right)\epsilon_k^2\quad \mbox{for all}\ x\in \mb{R}^d.
\end{align*}
This completes the proof of (\ref{equ:cauchyschwarz}). It only remains to prove (\ref{inequ:onestep}). Using the definition of $\tlx_{k+1}$, we have 
\begin{align*}
    0\in \partial \left(\varphi_k + \psi_k\right)\left(\tlx_{k+1}\right) &\overset{(\star)}{=} \partial \varphi_k\left(\tlx_{k+1}\right) + \nabla \psi_k \left(\tlx_{k+1}\right)\\ &= \partial \varphi_k\left(\tlx_{k+1}\right) + \alpha_k^{-1}M_k \left(\tlx_{k+1}-x_k\right),
\end{align*}
where $(\star)$ follows from \cite[Exercise~8.8(c)]{rockafellar2009variational} combined with the fact that $\varphi_k(\tlx_{k+1})$ is finite and $\psi_k(\cdot)$ is continuously differentiable on $\mb{R}^d$. The preceding estimate implies
\begin{equation}\label{equ:regularsubdiffsubproblem}
    \begin{split}
        \alpha_k^{-1}M_k \left(x_k - \tlx_{k+1}\right)\in \partial \varphi_k\left(\tlx_{k+1}\right)
    \end{split}
\end{equation}
and thus for any $x\in \mb{R}^d$, we have 
\begin{align*}
    &\olvphi_{x_k}\left(x;S_k^{1:m}\right) - \olvphi_{x_k}\left(\tlx_{k+1};S_k^{1:m}\right) = \varphi_k\left(x\right) - \varphi_k\left(\tlx_{k+1}\right)\\
        \geq\ &\lrangle{x-\tlx_{k+1},\alpha_k^{-1}M_k \left(x_k - \tlx_{k+1}\right)} - \frac{\oltau}{2}\Norm{x-\tlx_{k+1}}_2^2\\
        =\ &\frac{1}{2\alpha_k} \left(\Norm{x-\tlx_{k+1}}_{M_k}^2 + \Norm{x_k - \tlx_{k+1}}_{M_k}^2 - \Norm{x-x_k}_{M_k}^2\right) - \frac{\oltau}{2}\Norm{x-\tlx_{k+1}}_2^2\\
        \geq\ &\left(\frac{1}{2\alpha_k} - \frac{\oltau L_\infty^2}{2}\right)\Norm{x-\tlx_{k+1}}_{M_k}^2 + \frac{1}{2\alpha_k} \Norm{x_k -\tlx_{k+1}}_{M_k}^2 - \frac{1}{2\alpha_k} \Norm{x-x_k}_{M_k}^2,
\end{align*}
where the second equality uses (\ref{equ:abnormm}), the first inequality comes from Lemma \ref{lema:weaklyconvex} combined with (\ref{equ:regularsubdiffsubproblem}) and the condition that $\varphi_k(\cdot)$ is $\oltau$-weakly convex, and the last is due to Assumption \ref{asm:mk}. Given that $\tau \triangleq L_\infty^2\cdot \oltau$, the estimate (\ref{inequ:onestep}) follows directly from the inequality derived above. This completes the proof of Lemma \ref{lema:onestepprogress}.

\subsubsection{Proof of Lemma \ref{lema:welldefinedolxk}}\label{apxsubsubsec:moreauenvelopediff}
Let Assumption \ref{asm:stomodel}, \ref{asm:mk} and \ref{asm:sppaall} hold. Fix any $x\in \mb{R}^d$ and $k\in \mb{Z}_{+}$. 
\paragraph*{Proof of assertion (\romannumeral1).} 
Lemma \ref{lema:lipoff}(\romannumeral1) implies that $\phi(\cdot)$ is proper, closed and $(\oleta+\oltau)$-weakly convex. Using (\ref{equ:scvxsqnormm}) in the proof of Lemma \ref{lema:onestepprogress}, we have $\vartheta_{k-1}(\cdot)\triangleq \frac{\olrho}{2}\Vert \cdot - x\Vert_{M_{k-1}}^2$ is real-valued and $\frac{\olrho}{L_{k-1}^2}$-strongly convex. Recall that $\eta = L_\infty^2\cdot \oleta$ and $\tau = L_\infty^2\cdot \oltau$. Using Assumption \ref{asm:mk} and condition $\olrho > \eta + \tau$, we deduce 
\begin{align*}
    \frac{\olrho}{L_{k-1}^2} - \left(\oleta+\oltau\right) = \frac{\olrho - L_{k-1}^2 \left(\oleta+\oltau\right)}{L_{k-1}^2} \geq \frac{\olrho - L_\infty^2 \left(\oleta+\oltau\right)}{L_{k-1}^2} = \frac{\olrho - \left(\eta + \tau\right)}{L_{k-1}^2} > 0.
\end{align*}
It follows that $\phi_{k-1}(\cdot)\triangleq \phi(\cdot) + \vartheta_{k-1}(\cdot)$ is proper, closed and $(\frac{\olrho}{L_{k-1}^2} - \oleta-\oltau)$-strongly convex. By \cite[Lemma~2.26]{planiden2016strongly}, we can thus conclude that the minimizer of $\phi_{k-1}$ on $\mb{R}^d$ exists and is unique, which combined with the definition of the proximal mapping $\prox_{\left(\olrho\right)^{-1}\phi}^{M_{k-1}}(\cdot)$ proves assertion (\romannumeral1). 
\paragraph*{Proof of assertion (\romannumeral2).} Using the condition $\olrho > \eta+\tau$ and Assumption \ref{asm:mk}, we can deduce from (\ref{equ:gradmegenral}) that the Moreau envelope $e_{(\olrho)^{-1}\phi}^{M_{k-1}}(\cdot)$ is differentiable at $x$ with gradient given by 
\begin{equation}\label{equ:meformk}
    \begin{split}
        \nabla e_{\left(\olrho\right)^{-1}\phi}^{M_{k-1}}\left(x\right) = \olrho M_{k-1}\left(x - \prox_{\left(\olrho\right)^{-1}\phi}^{M_{k-1}}\left(x\right)\right).
    \end{split}
\end{equation}
Then, the estimate (\ref{equ:gradnormmkme}) follows immediately from (\ref{equ:meformk}) due to the definition of $M_{k-1}$-norm and $M_{k-1}^{-1}$-norm. 
For completeness, we provide a detailed proof of the estimate (\ref{equ:gradmegenral}) below.
\begin{proof}[Proof of the estimate (\ref{equ:gradmegenral})]
    Let $p\colon \mb{R}^d\to \pminf$ be a proper, closed and $\oltau$-weakly convex function. Suppose that $M\colon \mb{R}^d\to \mb{R}^d$ is a self-adjoint and positive-definite linear mapping satisfying $\Vert x\Vert_2 \leq L_M\Vert x\Vert_M$ for all $x\in \mb{R}^d$. Denote by $\tau \triangleq L_M^2\cdot\oltau$ and set $\alpha\in (0,\tau^{-1})$. Define 
    \begin{align*}
        p_\alpha\left(x\right) \triangleq p(x) + \frac{1}{2\alpha}\Norm{x}_M^2\quad \mbox{for all}\ x\in \mb{R}^d.
    \end{align*}
    By following an argument nearly identical to that in the proof of assertion (\romannumeral1) above, we can show that function $p_\alpha(\cdot)$ is proper, closed and $\frac{1-\alpha \tau}{L_M^2\alpha}$-strongly convex, and the proximal mapping $\prox_{\alpha p}^M(\cdot)$ is single-valued. Using the definition of the Moreau enevlope, we have
    \begin{subequations}
        \begin{align}
            e_{\alpha p}^M\left(x\right) 
            &= \inf_{y\in \mb{R}^d}\lrbrackets{p\left(y\right) + \frac{1}{2\alpha}\Norm{y-x}_M^2}\nonumber\\
            &= \frac{1}{2\alpha}\Norm{x}_M^2 - \sup_{y\in \mb{R}^d}\lrbrackets{\frac{1}{\alpha}\lrangle{Mx,y} - p_\alpha\left(y\right)}\label{equ:minmepalpha}\\
            &= \frac{1}{2\alpha}\Norm{x}_M^2 - p_\alpha^\star\left(\frac{Mx}{\alpha}\right)\quad \mbox{for all}\ x\in \mb{R}^d.\label{equ:meconjugate}
        \end{align}
    \end{subequations}
    where $p_\alpha^\star(\cdot)$ denotes the Fenchel conjugate of $p_\alpha$. Since $p_\alpha(\cdot)$ is proper, closed and $\frac{1-\alpha \tau}{L_M^2\alpha}$-strongly convex, then $p_\alpha^\star(\cdot)$ is differentiable by \cite[Proposition~12.60]{rockafellar2009variational}. Fix any $x\in \mb{R}^d$. We can therefore deduce from the estimate (\ref{equ:meconjugate}) that
    \begin{equation}\label{equ:gradme}
        \nabla e_{\alpha p}^M\left(x\right) = \frac{M}{\alpha}\left(x - \nabla p_\alpha^\star\left(\frac{Mx}{\alpha}\right)\right).
    \end{equation}
    Combining (\ref{equ:minmepalpha}) with the definition of the proximal mapping yields
    \begin{equation}\label{equ:proxpalpha}
        \begin{split}
            \prox_{\alpha p}^M\left(x\right) &= \arg\min_{y\in \mb{R}^d}\lrbrackets{p\left(y\right) + \frac{1}{2\alpha}\Norm{y-x}_M^2} = \arg\min_{y\in \mb{R}^d}\lrbrackets{p_\alpha\left(y\right) - \frac{1}{\alpha}\lrangle{Mx,y}}.
        \end{split}
    \end{equation}
    Consequently, we have
    \begin{align*}
        y= \nabla p_\alpha^\star\left(\frac{Mx}{\alpha}\right)\quad \Leftrightarrow\quad &\frac{Mx}{\alpha}\in \partial p_\alpha\left(y\right)\\
            \Leftrightarrow\quad &0\in \partial \left(p_\alpha - \frac{1}{\alpha}\lrangle{Mx,\cdot}\right)(y)\\
            \Leftrightarrow\quad &y= \prox_{\alpha p}^M\left(x\right),
    \end{align*}
    where the first two equivalence relations use the convexity of $p_\alpha(\cdot)$, and the last is due to (\ref{equ:proxpalpha}). Then, the estimate (\ref{equ:gradmegenral}) follows immediately from the combination of the preceding result with (\ref{equ:gradme}).
\end{proof}

\subsubsection{Proof of Lemma \ref{lema:inequalitydistance}} 
Inspired from the technique used in \cite[Lemma~4.2]{davis2019stochastic}, we provide a detailed proof of Lemma \ref{lema:inequalitydistance} as below. Recall that $B_\rho\triangleq \sup_{k\in\mb{Z}_{+}}\rho_k^2$ and $t_k \triangleq \sqrt{\mb{E}_k[\Vert \tlx_{k+1}-x_k\Vert_{M_k}^2]}$ for all $k\in \mb{Z}_{+}$. Under Assumption \ref{asm:mk}, the estimate (\ref{equ:rhonurelationapx}) leads us to $B_\rho\in [1,\nu_\infty^2]$. To deal with the minibatch case, we provide the following elementary result and refer to \cite[Lemma~B.5]{zhu2025tight} for a detailed proof.
\begin{lema}\label{lema:sumrandom}
    Let $X_1,\cdots,X_m$ be independent and identically distributed random variables with mean $\mb{E}[X_i] = \mu$ and variance $\mr{Var}(X_i) = \sigma^2$ for each $i = 1,\cdots,m$. If $\mb{E}[(X_1)^2]>0$, then 
    \begin{align*}
        \mb{E}\left[\left(\frac{1}{m}\sum_{i=1}^m X_i\right)^2\right] = \frac{1}{m}\sigma^2 + \mu^2 = \frac{1+(m-1)\omega}{m}\mb{E}[\left(X_1\right)^2],
    \end{align*}
    where $\omega \triangleq 1 - \frac{\mr{Var}(X_1)}{\mb{E}[(X_1)^2]}$.
\end{lema}
\paragraph*{Proof of assertion (\romannumeral1).} By assumption, $U$ is a bounded open convex subset of $V$ containing sequences $\{x_k\}$ and $\{\tlx_k\}$. Under Assumption \ref{asm:mk} and \ref{asm:sppa1}, the terms $\omega_{f,U}$ and $\rho_{f,m,U}$ are defined as follows:
\begin{align*}
    &\omega_{f,U} \triangleq 1 - \frac{\mr{Var}\left(L_{f,U}(s)\right)}{\mb{E}_{s\sim P}\left[L_{f,U}(s)^2\right]}\quad \mbox{and}\quad \rho_{f,m,U} \triangleq \left(\sqrt{\frac{1+(m-1)\omega_{f,U}}{m}}+1\right)^2\\
\end{align*} 
where $L_{f,U}(s) \triangleq L_\infty\cdot \olL_{f,U}(s)$ satisfies $\sqrt{\mb{E}_{s\sim P}[L_{f,U}(s)^2]}\leq L_F(U)$ with $L_F(U)\triangleq L_\infty\cdot \olL_F(U)$, and for $P$-almost all $s\in \mc{S}$,
\begin{equation}\label{equ:liplocalmknorm}
    \left\vert f_x(x;s) - f_x(y;s)\right\vert\leq L_{f,U}(s)\Vert x - y\Vert_{M_k}\quad \mbox{for all}\ x,y\in U
\end{equation}
holds for all $k\in \mb{Z}_{+}$. Fix any $k\in \mb{Z}_{+}$. 

\textit{\underline{*Step 1.}} Since $\sup_{k\in \mb{Z}_{+}}\alpha_k < \tau^{-1}$, it then follows from Lemma \ref{lema:onestepprogress} that the proximal point $\tlx_{k+1} \triangleq \prox_{\alpha_k \olvphi_{x_k}(\cdot;S_k^{1:m}}^{M_k}(x_k)$ is well-defined and (\ref{inequ:onestep}) holds. Denote by $\olf(\cdot;S_k^{1:m}) \triangleq \frac{1}{m}\sum_{i=1}^m f(\cdot;S_k^i)$ and $\hatL_{f,U}(S_k^{1:m}) \triangleq \frac{1}{m}\sum_{i=1}^m L_{f,U}(S_k^i)$. Fix any $x\in \mb{R}^d$. Under Assumption \ref{asm:stomodel}, \ref{asm:mk}, \ref{asm:sppaall} and \ref{asm:sppa1}, we have  
\begin{align*}
        &\olf\left(x;S_k^i\right) + r(x) + \frac{1}{2\alpha_k}\Norm{x-x_k}_{M_k}^2 + \frac{\eta}{2}\Norm{x-x_k}_{M_k}^2\\
        \geq\ &\olf_{x_k}\left(x;S_k^i\right) + r(x) + \frac{1}{2\alpha_k}\Norm{x-x_k}_{M_k}^2\\
        \geq\ &\olf_{x_k}\left(\tlx_{k+1};S_k^i\right) + r\left(\tlx_{k+1}\right) + \frac{1}{2\alpha_k}\Norm{\tlx_{k+1}-x_k}_{M_k}^2 + \frac{1-\tau \alpha_k}{2\alpha_k}\Norm{x - \tlx_{k+1}}_{M_k}^2\\
        \geq\ &
        \olf_{x_k}\left(x_k;S_k^i\right) + r\left(\tlx_{k+1}\right) + \frac{1}{2\alpha_k}\Norm{\tlx_{k+1}-x_k}_{M_k}^2 + \frac{1-\tau \alpha_k}{2\alpha_k}\Norm{x - \tlx_{k+1}}_{M_k}^2\\ & - \hatL_{f,U}\left(S_k^{1:m}\right) \cdot \Norm{\tlx_{k+1}-x_k}_{M_k}\\
        =\ &
        \olf\left(x_k;S_k^i\right) + r\left(\tlx_{k+1}\right) + \frac{1}{2\alpha_k}\Norm{\tlx_{k+1}-x_k}_{M_k}^2 + \frac{1-\tau \alpha_k}{2\alpha_k}\Norm{x - \tlx_{k+1}}_{M_k}^2\\
        &- \hatL_{f,U}\left(S_k^{1:m}\right) \cdot \Norm{\tlx_{k+1}-x_k}_{M_k},\\
\end{align*}
where the first inequality comes from (\ref{equ:stomodelmknorm}), the second is due to the estimate (\ref{inequ:onestep}) in Lemma \ref{lema:onestepprogress}, the last uses (\ref{equ:liplocalmknorm}) together with the assumption that $x_k,\tlx_{k+1}\in U$, and the equality follows from Assumption \ref{asm:stomodel}. 
Taking the conditional expectation with respect to $\mc{F}_{k-1}$ on both sides of the previous inequality gives
\begin{equation}\label{inequ:condexp}
    \begin{split}
        &\frac{1-\tau \alpha_k}{2\alpha_k}\mb{E}_k\left[\Norm{\tilde{x}_{k+1} - x}_{M_k}^2\right] - \frac{1+\eta \alpha_k}{2\alpha_k}\mb{E}_k \left[\Norm{x_k - x}_{M_k}^2\right]\\ 
        \leq\ & - \frac{1}{2\alpha_k}\mb{E}_k\left[\Norm{\tlx_{k+1}-x_k}_{M_k}^2\right]+ \mb{E}_k\left[F(x) + r(x) - F(x_k) - r\left(\tlx_{k+1}\right)\right]\\
        &\hspace{0em} + \mb{E}_k \left[\hatL_{f,U}\left(S_k^{1:m}\right) \cdot \Norm{\tlx_{k+1}-x_k}_{M_k}\right]\\
        \leq\ & - \frac{1}{2\alpha_k}\mb{E}_k\left[\Norm{\tlx_{k+1}-x_k}_{M_k}^2\right] + \mb{E}_k\left[\phi(x) - F(x_k) - r\left(\tlx_{k+1}\right)\right]\\
        &\hspace{0em} + \sqrt{\mb{E}_k\left[\left(\hatL_{f,U}\left(S_k^{1:m}\right)\right)^2\right]}\cdot\sqrt{\mb{E}_k\left[\Vert \tilde{x}_{k+1} - x_k\Vert_{M_k}^2\right]}\\
        \leq\ & - \frac{1}{2\alpha_k}\mb{E}_k\left[\Norm{\tlx_{k+1}-x_k}_{M_k}^2\right] + \mb{E}_k\left[\phi(x) - F(x_k) - r\left(\tlx_{k+1}\right)\right]\\
        &\hspace{0em} + \sqrt{\frac{1+(m-1)\omega_{f,U}}{m}}L_F\left(U\right)\cdot \sqrt{\mb{E}_k\left[\Vert \tilde{x}_{k+1} - x_k\Vert_{M_k}^2\right]},
    \end{split}
\end{equation} 
where the second inequality is due to the Cauchy-Schwarz inequality, and the last follows from Lemma \ref{lema:sumrandom} together with the fact that $\sqrt{\mb{E}_{s\sim P}[L_{f,U}(s)^2]}\leq L_F(U)$. 
By Jensen's inequality, the following assertion holds: 
\begin{equation}\label{inequ:jensentk}
    \mb{E}_k\left[\Vert \tilde{x}_{k+1} - x_k\Vert_{M_k}\right]\leq t_k.
\end{equation}
Using the estimate (\ref{inequ:condexp}), we can obtain
\begin{align*}
        &\frac{1-\tau \alpha_k}{2\alpha_k}\mb{E}_k\left[\Norm{\tilde{x}_{k+1} - x}_{M_k}^2\right] - \frac{1+\eta \alpha_k}{2\alpha_k}\mb{E}_k \left[\Norm{x_k - x}_{M_k}^2\right]\\ 
        \leq\ &\mb{E}_k\left[\phi(x) - \left(F\left(\tlx_{k+1}\right) - L_F\left(U\right) \Norm{\tlx_{k+1}-x_k}_{M_k}\right) - r\left(\tlx_{k+1}\right)\right]\\
        &\hspace{0em} - \frac{1}{2\alpha_k}\mb{E}_k\left[\Norm{\tlx_{k+1}-x_k}_{M_k}^2\right] + \sqrt{\frac{1+(m-1)\omega_{f,U}}{m}}L_F\left(U\right)\cdot \sqrt{\mb{E}_k\left[\Vert \tilde{x}_{k+1} - x_k\Vert_{M_k}^2\right]}\\
        \leq\ &\mb{E}_k\left[\phi(x) - \phi\left(\tlx_{k+1}\right)\right] - \frac{1}{2\alpha_k}t_k^2 + \sqrt{\rho_{f,m,U}}L_F\left(U\right) t_k\\
        \leq\ &\mb{E}_k\left[\phi(x) - \phi\left(\tlx_{k+1}\right)\right] + \max_{t\geq 0}\lrbrackets{ \sqrt{\rho_{f,m,U}}L_F\left(U\right) t - \frac{1}{2\alpha_k}t^2}\\
        =\ & \mb{E}_k\left[\phi(x) - \phi\left(\tlx_{k+1}\right)\right] + \frac{\rho_{f,m,U}L_F(U)^2}{2}\alpha_k,
\end{align*}
where the first inequality uses Lemma \ref{lema:lipoff}(\romannumeral2) combined with Assumption \ref{asm:mk} and condition $x_k,\tlx_{k+1}\in U$, and the second follows from (\ref{inequ:jensentk}) together with the definition of $\rho_{f,m,U}$. This completes the proof of both (\ref{equ:inequ1}) and (\ref{equ:inequ1max}). 

\textit{\underline{*Step 2.}} Now consider the case where $\sup_{k\in \mb{Z}_{+}}\alpha_k < (\olrho\max\{2,B_\rho\})^{-1}$ for $\olrho > B_\rho \cdot(\eta + \tau)$. Recall that $\olx_k \triangleq \prox_{(\olrho)^{-1} \phi}^{M_{k-1}}(x_k)$. Setting $x = \olx_k$ in (\ref{equ:inequ1}) and using the definition of the proximal mapping, we have 
\begin{align*}
        &\frac{1-\tau \alpha_k}{2\alpha_k}\mb{E}_k\left[\Norm{\tilde{x}_{k+1} - \olx_k}_{M_k}^2\right]\\ 
        \leq\ &\frac{1+\eta \alpha_k}{2\alpha_k}\mb{E}_k \left[\Norm{x_k - \olx_k}_{M_k}^2\right] + \mb{E}_k \left[-\ \frac{\olrho}{2}\Norm{\olx_k - x_k}_{M_{k-1}}^2 + \frac{\olrho}{2}\Norm{\tlx_{k+1} - x_k}_{M_{k-1}}^2\right]\\
        &\hspace{0em} + \mb{E}_k\left[\left(\phi\left(\olx_k\right) + \frac{\olrho}{2}\Norm{\olx_k - x_k}_{M_{k-1}}^2\right) - \left(\phi\left(\tlx_{k+1}\right) + \frac{\olrho}{2}\Norm{\tlx_{k+1} - x_k}_{M_{k-1}}^2\right) \right]\\
        &\hspace{0em} - \frac{1}{2\alpha_k}t_k^2 + \sqrt{\rho_{f,m,U}}L_F\left(U\right) t_k\\
        \leq\ &\frac{1+\eta \alpha_k}{2\alpha_k}\mb{E}_k \left[\Norm{x_k - \olx_k}_{M_k}^2\right] - \frac{\olrho \rho_k^{-2}}{2}\mb{E}_k \left[\Norm{\olx_k - x_k}_{M_{k}}^2\right] + \frac{\olrho\rho_k^2}{2}t_k^2\\
        &\hspace{0em} - \frac{1}{2\alpha_k}t_k^2 + \sqrt{\rho_{f,m,U}}L_F\left(U\right) t_k\\
        \leq\ &\frac{1+\eta \alpha_k - \olrho \rho_k^{-2}\alpha_k}{2\alpha_k}\mb{E}_k \left[\Norm{x_k - \olx_k}_{M_k}^2\right] + \max_{t\geq 0}\lrbrackets{\sqrt{\rho_{f,m,U}}L_F\left(U\right) t -\frac{1-\olrho\rho_k^2\alpha_k}{2\alpha_k}t^2}\\
        =\ &\frac{1+\eta \alpha_k - \olrho \rho_k^{-2}\alpha_k}{2\alpha_k}\mb{E}_k \left[\Norm{x_k - \olx_k}_{M_k}^2\right] + \frac{\rho_{f,m,U}L_F\left(U\right)^2\alpha_k}{2\left(1-\olrho\rho_k^2\alpha_k\right)},
\end{align*}
where the second inequality uses Assumption \ref{asm:mk}, and the equality is due to $\sup_{k\in \mb{Z}_{+}}\alpha_k < (\olrho\max\{2,B_\rho\})^{-1}$. Note that $1-\tau \alpha_k \geq 1-\olrho\alpha_k > 0$ since $\olrho > B_\rho \cdot(\eta + \tau)$ and $\sup_{k\in \mb{Z}_{+}}\alpha_k < (\olrho\max\{2,B_\rho\})^{-1}$. Rearranging the previous result and using Assumption \ref{asm:mk}, one can conclude that 
\begin{equation}\label{equ:inequ1further}
    \begin{split}
        &\mb{E}_k\left[\Norm{\tilde{x}_{k+1} - \olx_k}_{M_k}^2\right]\\
        \leq\ &\frac{1+\eta \alpha_k - \olrho \rho_k^{-2}\alpha_k}{1-\tau \alpha_k} \mb{E}_k \left[\Norm{x_k - \olx_k}_{M_k}^2\right] + \frac{\rho_{f,m,U}L_F\left(U\right)^2}{\left(1-\tau \alpha_k\right)\left(1-\olrho\rho_k^2\alpha_k\right)}\alpha_k^2\\
        \leq\ &\frac{\rho_k^2\left(1+\eta \alpha_k\right) - \olrho \alpha_k}{1-\tau \alpha_k} \Norm{x_k - \olx_k}_{M_{k-1}}^2 + \frac{\rho_{f,m,U}L_F\left(U\right)^2}{\left(1-\tau \alpha_k\right)\left(1-\olrho\rho_k^2\alpha_k\right)}\alpha_k^2.
    \end{split}
\end{equation}
Fix any $t\in \mb{R}_{++}$. 
Combining (\ref{equ:inequ1further}) and (\ref{equ:cauchyschwarz}) in Lemma \ref{lema:onestepprogress} with $x = \olx_k$ yields 
\begin{equation}\label{equ:sippinequ1}
    \begin{split}
        &\mb{E}_k\left[\Vert x_{k+1} - \olx_k\Vert_{M_k}^2\right]\\ 
        \leq\ &\frac{\left(1+\frac{1}{t}\right)\left[\rho_k^2\left(1+\eta \alpha_k\right) - \olrho \alpha_k\right]}{1-\tau \alpha_k}\Norm{x_k - \olx_k}_{M_{k-1}}^2\\
        &\hspace{0em}+\ \frac{\left(1+\frac{1}{t}\right) \rho_{f,m,U}L_F\left(U\right)^2}{\left(1-\tau \alpha_k\right)\left(1-\olrho\rho_k^2\alpha_k\right)}\cdot\alpha_k^2 + (1 + t)\epsilon_k^2\quad \mbox{for all}\ t\in \mb{R}_{++}.
    \end{split}
\end{equation}
\textit{\underline{*Step 3.}} Now consider the case where $\phi$ is bounded from below. Then, by using the definition of $\olx_k$ and $\olx_{k+1}$, we have
\begin{equation}
    \begin{split}
        &\mb{E}_k\left[\mephi{k+1}\left(x_{k+1}\right) - \phi^{*}\right] = \mb{E}_k \left[\phi\left(\olx_{k+1}\right) + \frac{\olrho}{2}\Norm{\olx_{k+1} - x_{k+1}}_{M_k}^2 - \phi^{*}\right]\\
        \leq\ &\mb{E}_k \left[\phi\left(\olx_{k}\right) + \frac{\olrho}{2}\Norm{\olx_{k} - x_{k+1}}_{M_k}^2 - \phi^{*}\right] \\
        \leq\ &\rho_k^2\cdot\left[\phi\left(\olx_{k}\right) - \phi^{*} + \frac{\olrho}{2}\Norm{\olx_k - x_k}_{M_{k-1}}^2\right]\\
        &\hspace{0em}+\ \frac{\olrho\left(1+\frac{1}{t}\right)\rho_{f,m,U}L_F\left(U\right)^2}{2\left(1-\tau \alpha_k\right)\left(1-\olrho\rho_k^2\alpha_k\right)}\alpha_k^2 + \frac{\olrho\left(1+t\right)}{2}\epsilon_k^2\\
        &\hspace{0em}+\ \frac{\left(1+\frac{1}{t}\right)\left[\rho_k^2\left(1+\eta \alpha_k\right) - \olrho \alpha_k\right] - \rho_k^2\left(1-\tau \alpha_k\right)}{2\olrho\left(1-\tau \alpha_k\right)}\Norm{\olrho \left(x_k - \olx_k\right)}_{M_{k-1}}^2,
    \end{split}
\end{equation}
holds for any $t\in \mb{R}_{++}$, where the last inequality follows from (\ref{equ:sippinequ1}) combined with the fact that $\rho_k \geq 1$ and $\phi(\olx_k) \geq \phi^{*}$. Note that $\nu_{k+1} = \rho_k \nu_k$ and $\nu_{k+1} \geq 1$. Multiplying both sides of the previous inequality by $\nu_{k+1}^{-2}$ and using the definition of the Moreau envelope yields
\begin{equation}\label{equ:sippinequhatx1tvalueuse2}
    \begin{split}
        &\mb{E}_k\left[\nu_{k+1}^{-2}\cdot\left(\mephi{k+1}\left(x_{k+1}\right) - \phi^{*}\right)\right]\\ 
        \leq\ &\nu_k^{-2}\cdot\left(\mephi{k}\left(x_k\right) - \phi^{*}\right) + \frac{\olrho\left(1+\frac{1}{t}\right)\rho_{f,m,U}L_F\left(U\right)^2}{2\left(1-\tau \alpha_k\right)\left(1-\olrho\rho_k^2\alpha_k\right)}\alpha_k^2 + \frac{\olrho\left(1+t\right)}{2}\epsilon_k^2\\
        &\hspace{0em}+\ \frac{\left(1+\frac{1}{t}\right)\left[1+\eta \alpha_k - \olrho \rho_k^{-2}\alpha_k\right] - \left(1-\tau \alpha_k\right)}{2\olrho\left(1-\tau \alpha_k\right)}\cdot \nu_k^{-2}\Norm{\olrho \left(x_k - \olx_k\right)}_{M_{k-1}}^2.\\
    \end{split}
\end{equation}
Let $b \triangleq 4(\olrho\rho_k^{-2} - \eta)(\olrho\rho_k^{-2}-\eta-\tau)$. Condition $\olrho > B_\rho \cdot (\eta+\tau)$ implies that $b>0$. Then, by setting $t = \frac{1}{b\alpha_k^2}$ in (\ref{equ:sippinequhatx1tvalueuse2}), we can conclude that estimate (\ref{equ:mephiinequ1}) holds. 
This completes the proof of assertion (\romannumeral1).

\paragraph*{Proof of assertion (\romannumeral2).} Under Assumption \ref{asm:mk} and Assumption \ref{asm:generalizedlip2}, the terms $\omega_{f}$ and $\rho_{f,m}$ are defined as follows:
\begin{align*}
        &\omega_{f} \triangleq 1 - \frac{\mr{Var}\left(L_{f}(s)\right)}{\mb{E}_{s\sim P}\left[L_{f}(s)^2\right]}\quad \mbox{and}\quad \rho_{f,m} \triangleq \left(\sqrt{\frac{1+(m-1)\omega_{f}}{m}}+1\right)^2\\
\end{align*} 
where $L_{f}(s) \triangleq L_\infty\cdot \olL_{f}(s)$ satisfies $\sqrt{\mb{E}_{s\sim P}[L_{f}(s)^2]}\leq L_F$ with $L_F\triangleq L_\infty\cdot \olL_F$, and for $P$-almost all $s\in \mc{S}$,
\begin{equation}\label{equ:lipgenrealized2mknorm}
    \left\vert f_x(x;s) - f_x(y;s)\right\vert\leq L_{f}(s)\glip\left(\Norm{x}_2\right)\Vert x - y\Vert_{M_k}\quad \mbox{for all}\ x,y\in V
\end{equation}
holds for all $k\in \mb{Z}_{+}$. Fix any $k\in \mb{Z}_{+}$. Since $\dom(r)\subset V$ by Assumption \ref{asm:stomodel}, we can deduce from the definition of $x_k$ and $\tlx_{k+1}$ that $x_k,\tlx_{k+1}\in V$. Then, by following an argument nearly identical to that in the proof of assertion (\romannumeral1) derived above, we can deduce from Lemma \ref{lema:lipoff}(\romannumeral3) and (\ref{equ:lipgenrealized2mknorm}) together with Assumption \ref{asm:mk} that 
\begin{align*}
    &\frac{1-\tau \alpha_k}{2\alpha_k}\mb{E}_k\left[\Vert \tilde{x}_{k+1} - x\Vert_{M_k}^2\right] - \frac{1+\eta \alpha_k}{2\alpha_k}\mb{E}_k\left[\Vert x_k - x\Vert_{M_k}^2\right]\\ \leq\ & \mb{E}_k\left[\phi\left(x\right) - \phi\left(\tilde{x}_{k+1}\right)\right] - \frac{1-\eta\alpha_k}{2\alpha_k}t_k^2 + \sqrt{\rho_{f,m}}L_F\glip\left(\Norm{x_k}_2\right) t_k,
\end{align*}
holds for all $x\in \mb{R}^d$, which proves the estimate (\ref{equ:inequ2}). By using (\ref{equ:inequ2}) and condition $\sup_{k\in \mb{Z}_{+}}\alpha_k < (\olrho\max\{2,B_\rho\} + \eta)^{-1}$, and adopting the proof of (\ref{equ:sippinequ1}) with some modifications, one can deduce that
\begin{equation}\label{equ:sippinequ2}
    \begin{split}
        &\mb{E}_k\left[\Vert x_{k+1} - \olx_k\Vert_{M_k}^2\right]\\ 
        \leq\ &\frac{\left(1+\frac{1}{t}\right)\left[\rho_k^2\left(1+\eta \alpha_k\right) - \olrho \alpha_k\right]}{1-\tau \alpha_k}\Norm{x_k - \olx_k}_{M_{k-1}}^2\\
        &\hspace{0em}+\ \frac{\left(1+\frac{1}{t}\right) \rho_{f,m}L_F^2\glip\left(\Norm{x_k}_2\right)^2}{\left(1-\tau \alpha_k\right)\left(1-\olrho\rho_k^2\alpha_k-\eta\alpha_k\right)}\cdot\alpha_k^2 + (1 + t)\epsilon_k^2\quad \mbox{for all}\ t\in \mb{R}_{++}.
    \end{split}
\end{equation}
Assume further that $\phi$ is bounded from below. By following the proof argument of (\ref{equ:sippinequhatx1tvalueuse2}), we can deduce from (\ref{equ:sippinequ2}) that 
\begin{equation}\label{equ:sippinequhatx1tvalueglip2}
    \begin{split}
        &\mb{E}_k\left[\nu_{k+1}^{-2}\cdot\left(\mephi{k+1}\left(x_{k+1}\right) - \phi^{*}\right)\right] - \nu_k^{-2}\cdot\left(\mephi{k}\left(x_k\right) - \phi^{*}\right) \\ 
        \leq\ & \frac{\olrho\left(1+\frac{1}{t}\right)\rho_{f,m}L_F^2\glip\left(\Norm{x_k}_2\right)^2}{2\left(1-\tau \alpha_k\right)\left(1-\olrho\rho_k^2\alpha_k-\eta\alpha_k\right)}\alpha_k^2 + \frac{\olrho\left(1+t\right)}{2}\epsilon_k^2\\
        &\hspace{0em}+\ \frac{\left(1+\frac{1}{t}\right)\left[1+\eta \alpha_k - \olrho \rho_k^{-2}\alpha_k\right] - \left(1-\tau \alpha_k\right)}{2\olrho\left(1-\tau \alpha_k\right)}\cdot \nu_k^{-2}\Norm{\olrho \left(x_k - \olx_k\right)}_{M_{k-1}}^2.\\
    \end{split}
\end{equation}
Let $b \triangleq 4(\olrho\rho_k^{-2} - \eta)(\olrho\rho_k^{-2}-\eta-\tau)$. Since $\olrho > B_\rho\cdot (\eta+\tau)$, we have $b>0$. Then, estimate (\ref{equ:mephiinequ2}) follows immediately from (\ref{equ:sippinequhatx1tvalueglip2}) with $t = \frac{1}{b\alpha_k^2}$. This completes the proof of assertion (\romannumeral2).

\paragraph*{Proof of assertion (\romannumeral3).} Under Assumption \ref{asm:mk} and Assumption \ref{asm:generalizedlip3}, the terms $\omega_{\varphi}$ and $\rho_{\varphi,m}$ are defined as follows:
\begin{align*}
    &\omega_{\varphi} \triangleq 1 - \frac{\mr{Var}\left(L_{\varphi}(s)\right)}{\mb{E}_{s\sim P}\left[L_{\varphi}(s)^2\right]}\quad \mbox{and}\quad \rho_{\varphi,m} \triangleq \frac{1+(m-1)\omega_{\varphi}}{m}\\
\end{align*} 
where $L_{\varphi}(s) \triangleq L_\infty\cdot \olL_{\varphi}(s)$ satisfies $\sqrt{\mb{E}_{s\sim P}[L_{\varphi}(s)^2]}\leq L_\phi$ with $L_\phi\triangleq L_\infty\cdot \olL_\phi$. Fix any $k\in \mb{Z}_{+}$. Since $x_k\in V$ as claimed in the proof of assertion (\romannumeral2), Lemma \ref{lema:lipoff}(\romannumeral4) implies that 
\begin{equation}\label{equ:lipgenrealized3mknorm}
    \begin{split}
        \Norm{u}_2 \leq\ &\left(\frac{1}{m}\sum_{i=1}^m \olL_\varphi\left(S_k^i\right)\right) \glip\left(\max\lrbrackets{\Norm{x_k}_2,\Norm{x}_2}\right)\\ &\mbox{for all}\ u\in \partial \olvphi_{x_k}\left(x;S_k^{1:m}\right)
    \ \mbox{with}\ x\in V.
    \end{split}
\end{equation}

\textit{\underline{*Step 1.}} As demonstrated in the proof of Lemma \ref{lema:onestepprogress},  Assumption \ref{asm:sppaall} implies that $\olvphi_{x_k}(\cdot;S_k^{1:m})$ is $\oltau$-weakly convex and thus Clarke regular. Therefore, combining \cite[Proposition~4.3.1~and~Definition~4.3.4]{cui2021modern} with the condition $x_k\in V$, we can conclude that there exists a vector $\olvphi_{x_k}^\prime (x_k;S_k^i)$ belonging to $\partial \olvphi_{x_k} (x_k;S_k^{1:m})$. Denote by $\hatL_\varphi(S_k^{1:m}) \triangleq \frac{1}{m}\sum_{i=1}^m L_\varphi(S_k^i)$. It then follows from Lemma \ref{lema:weaklyconvex} combined with the Cauchy-Schwarz inequality that 
\begin{align}
        &\olvphi_{x_k} \left(x_k;S_k^{1:m}\right) - \olvphi_{x_k} \left(\tlx_{k+1};S_k^{1:m}\right)\nonumber \\ 
        \leq\ &\Norm{\olvphi_{x_k}^\prime \left(x_k;S_k^i\right)}_2\cdot \Norm{x_k - \tlx_{k+1}}_2 + \frac{\oltau}{2}\Norm{x_k - \tlx_{k+1}}_2^2\nonumber\\
        \leq\ &\left(\frac{1}{m}\sum_{i=1}^m \olL_\varphi\left(S_k^i\right)\right)\glip\left(\max\lrbrackets{\Norm{x_k}_2,\Norm{x_k}_2}\right)\cdot \Norm{x_k - \tlx_{k+1}}_2 + \frac{\oltau}{2}\Norm{x_k - \tlx_{k+1}}_2^2\nonumber\\
        \leq\ &\hatL_\varphi\left(S_k^{1:m}\right) \glip\left(\Norm{x_k}_2\right)\cdot \Norm{x_k - \tlx_{k+1}}_{M_k} + \frac{\tau}{2}\Norm{x_k - \tlx_{k+1}}_{M_k}^2,\label{equ:normxk1glip}
\end{align}
where the second last inequality follows from (\ref{equ:lipgenrealized3mknorm}) together with the fact that $x_k\in V$, and the last is due to Assumption \ref{asm:mk} and the definition of $\hatL_\varphi(\cdot)$ and $\tau$. Denote by $\olvphi(\cdot;S_k^{1:m}) \triangleq \frac{1}{m}\sum_{i=1}^m \varphi(\cdot;S_k^i)$. Fix any $x\in \mb{R}^d$. Using Assumption \ref{asm:stomodel}, \ref{asm:mk}, \ref{asm:sppaall} and \ref{asm:generalizedlip3}, we have 
\begin{align*}
        &\olvphi\left(x;S_k^{1:m}\right) + \frac{1}{2\alpha_k}\Norm{x-x_k}_{M_k}^2 + \frac{\eta}{2}\Norm{x-x_k}_{M_k}^2\\
        \geq\ &\olvphi_{x_k}\left(x;S_k^{1:m}\right) + \frac{1}{2\alpha_k}\Norm{x-x_k}_{M_k}^2\\
        \geq\ &\olvphi_{x_k}\left(\tlx_{k+1};S_k^{1:m}\right) + \frac{1}{2\alpha_k}\Norm{\tlx_{k+1}-x_k}_{M_k}^2 + \frac{1-\tau \alpha_k}{2\alpha_k}\Norm{x - \tlx_{k+1}}_{M_k}^2\\
        \geq\ &
        \olvphi_{x_k}\left(x_k;S_k^{1:m}\right) + \frac{1}{2\alpha_k}\Norm{\tlx_{k+1}-x_k}_{M_k}^2 \\
        &\hspace{0em} + \frac{1-2\tau \alpha_k}{2\alpha_k}\Norm{x - \tlx_{k+1}}_{M_k}^2 - \hatL_\varphi\left(S_k^{1:m}\right)\glip\left(\Norm{x_k}_2\right) \cdot \Norm{\tlx_{k+1}-x_k}_{M_k}\\
        =\ &
        \olvphi\left(x_k;S_k^{1:m}\right) + \frac{1}{2\alpha_k}\Norm{\tlx_{k+1}-x_k}_{M_k}^2\\
        &\hspace{0em} + \frac{1-2\tau \alpha_k}{2\alpha_k}\Norm{x - \tlx_{k+1}}_{M_k}^2 - \hatL_\varphi\left(S_k^{1:m}\right)\glip\left(\Norm{x_k}_2\right) \cdot \Norm{\tlx_{k+1}-x_k}_{M_k},\\
\end{align*}
where the first inequality comes from (\ref{equ:stomodelmknorm}), the second is due to (\ref{inequ:onestep}) in Lemma \ref{lema:onestepprogress}, the last follows from (\ref{equ:normxk1glip}), and the equality uses Assumption \ref{asm:stomodel}. Taking the conditional expectation with respect to $\mc{F}_{k-1}$ on both sides of the previous inequality gives
\begin{equation}\label{inequ:condexpglip}
    \begin{split}
        &\frac{1-2\tau \alpha_k}{2\alpha_k}\mb{E}_k\left[\Norm{\tilde{x}_{k+1} - x}_{M_k}^2\right] - \frac{1+\eta \alpha_k}{2\alpha_k}\mb{E}_k \left[\Norm{x_k - x}_{M_k}^2\right]\\ 
        \leq\ & - \frac{1}{2\alpha_k}\mb{E}_k\left[\Norm{\tlx_{k+1}-x_k}_{M_k}^2\right] + \mb{E}_k\left[\phi(x)- \phi\left(x_k\right)\right]\\
        &\hspace{0em} + \mb{E}_k \left[\hatL_\varphi\left(S_k^{1:m}\right)\glip\left(\Norm{x_k}_2\right) \cdot \Norm{\tlx_{k+1}-x_k}_{M_k}\right]\\
        \leq\ &- \frac{1}{2\alpha_k}\mb{E}_k\left[\Norm{\tlx_{k+1}-x_k}_{M_k}^2\right] + \mb{E}_k\left[\phi(x) - \phi\left(x_k\right)\right]\\
        &\hspace{0em} + \sqrt{\mb{E}_k\left[\left(\hatL_\varphi\left(S_k^{1:m}\right)\glip\left(\Norm{x_k}_2\right)\right)^2\right]}\cdot\sqrt{\mb{E}_k\left[\Vert \tilde{x}_{k+1} - x_k\Vert_{M_k}^2\right]}\\
        \leq\ &- \frac{1}{2\alpha_k}\mb{E}_k\left[\Norm{\tlx_{k+1}-x_k}_{M_k}^2\right]  + \mb{E}_k\left[\phi(x) - \phi\left(x_k\right)\right]\\
        &\hspace{0em} + \sqrt{\frac{1+(m-1)\omega_{\varphi}}{m}}L_\phi\glip\left(\Norm{x_k}_2\right)\cdot \sqrt{\mb{E}_k\left[\Vert \tilde{x}_{k+1} - x_k\Vert_{M_k}^2\right]},\\
    \end{split}
\end{equation} 
where the second inequality is due to the Cauchy-Schwarz inequality, and the last follows from Lemma \ref{lema:sumrandom} together with the fact that $\sqrt{\mb{E}_{s\sim P}[L_{\varphi}(s)^2]}\leq L_\phi$. 
Recall that $t_k \triangleq \sqrt{\mb{E}_k[\Vert \tlx_{k+1}-x_k\Vert_{M_k}^2]}$. Using the definition of $\rho_{\varphi,m}$, we can deduce from (\ref{inequ:condexpglip}) that 
\begin{align*}
        &\frac{1-2\tau \alpha_k}{2\alpha_k}\mb{E}_k\left[\Norm{\tilde{x}_{k+1} - x}_{M_k}^2\right] - \frac{1+\eta \alpha_k}{2\alpha_k}\mb{E}_k \left[\Norm{x_k - x}_{M_k}^2\right]\\ 
        \leq\ &\mb{E}_k\left[\phi\left(x\right) - \phi\left(x_k\right)\right] - \frac{1}{2\alpha_k}t_k^2 + \sqrt{\rho_{\varphi,m}}L_\phi\glip\left(\Norm{x_k}_2\right) t_k.
\end{align*}
This completes the proof of the estimate (\ref{equ:inequ3}). 

\textit{\underline{*Step 2.}} Recall that $\olx_k = \prox_{(\olrho)^{-1} \phi}^{M_{k-1}}(x_k)$. By setting $x = \olx_k$ in (\ref{equ:inequ3}) and using the definition of the proximal mapping, we deduce
\begin{align*}
        &\frac{1-2\tau \alpha_k}{2\alpha_k}\mb{E}_k\left[\Norm{\tilde{x}_{k+1} - \olx_k}_{M_k}^2\right]\\ 
        \leq\ &\frac{1+\eta \alpha_k}{2\alpha_k}\mb{E}_k \left[\Norm{x_k - \olx_k}_{M_k}^2\right] + \mb{E}_k \left[-\ \frac{\olrho}{2}\Norm{\olx_k - x_k}_{M_{k-1}}^2 \right]\\
        &\hspace{0em} + \mb{E}_k\left[\phi\left(\olx_k\right) + \frac{\olrho}{2}\Norm{\olx_k - x_k}_{M_{k-1}}^2 - \phi\left(x_k\right) \right] - \frac{1}{2\alpha_k}t_k^2 + \sqrt{\rho_{\varphi,m}}L_\phi\glip\left(\Norm{x_k}_2\right) t_k\\
        \leq\ &\frac{1+\eta \alpha_k}{2\alpha_k}\mb{E}_k \left[\Norm{x_k - \olx_k}_{M_k}^2\right] - \frac{\olrho \rho_k^{-2}}{2}\mb{E}_k \left[\Norm{\olx_k - x_k}_{M_{k}}^2\right]\\
        &\hspace{0em}- \frac{1}{2\alpha_k}t_k^2 + \sqrt{\rho_{\varphi,m}}L_\phi\glip\left(\Norm{x_k}_2\right) t_k\\
        \leq\ &\frac{1+\eta \alpha_k - \olrho \rho_k^{-2}\alpha_k}{2\alpha_k}\mb{E}_k \left[\Norm{x_k - \olx_k}_{M_k}^2\right] + \max_{t\geq 0}\lrbrackets{ \sqrt{\rho_{\varphi,m}}L_\phi\glip\left(\Norm{x_k}_2\right) t -\frac{1}{2\alpha_k}t^2}\\
        =\ &\frac{1+\eta \alpha_k - \olrho \rho_k^{-2}\alpha_k}{2\alpha_k}\mb{E}_k \left[\Norm{x_k - \olx_k}_{M_k}^2\right] + \frac{\rho_{\varphi,m}L_\phi^2\glip\left(\Norm{x_k}_2\right)^2\alpha_k}{2},
\end{align*}
where the second inequality uses Assumption \ref{asm:mk}. Note that $1-2\tau\alpha_k\geq 1-\olrho \alpha_k> 0$ since $\olrho > B_\rho\cdot(\eta+2\tau)$ and $\sup_{k\in \mb{Z}_{+}}\alpha_k \leq \sup_{k\in \mb{Z}_{+}}\olalpha_k < (2\olrho)^{-1}$. Rearranging the preceding inequality and using Assumption \ref{asm:mk}, one can conclude that 
\begin{equation}\label{equ:inequ1furtherglip3}
    \begin{split}
        &\mb{E}_k\left[\Norm{\tilde{x}_{k+1} - \olx_k}_{M_k}^2\right]\\
        \leq\ &\frac{1+\eta \alpha_k - \olrho \rho_k^{-2}\alpha_k}{1-2\tau \alpha_k} \mb{E}_k \left[\Norm{x_k - \olx_k}_{M_k}^2\right] + \frac{\rho_{\varphi,m}L_\phi^2\glip\left(\Norm{x_k}_2\right)^2}{1-2\tau \alpha_k}\alpha_k^2\\
        \leq\ &\frac{\rho_k^2\left(1+\eta \alpha_k\right) - \olrho \alpha_k}{1-2\tau \alpha_k} \Norm{x_k - \olx_k}_{M_{k-1}}^2 + \frac{\rho_{\varphi,m}L_\phi^2\glip\left(\Norm{x_k}_2\right)^2}{1-2\tau \alpha_k}\alpha_k^2.
    \end{split}
\end{equation}
Fix any $t\in \mb{R}_{++}$. Combining (\ref{equ:inequ1furtherglip3}) with (\ref{equ:cauchyschwarz}) in Lemma \ref{lema:onestepprogress} leads us to
\begin{equation}\label{equ:sippinequ1glip3}
    \begin{split}
        \mb{E}_k\left[\Vert x_{k+1} - \olx_k\Vert_{M_k}^2\right] &
        \leq \frac{\left(1+\frac{1}{t}\right)\left[\rho_k^2\left(1+\eta \alpha_k\right) - \olrho \alpha_k\right]}{1-2\tau \alpha_k}\Norm{x_k - \olx_k}_{M_{k-1}}^2\\
        &\hspace{0em}+\ \frac{\left(1+\frac{1}{t}\right) \rho_{\varphi,m}L_\phi^2\glip\left(\Norm{x_k}_2\right)^2}{1-2\tau \alpha_k}\cdot\alpha_k^2 + (1 + t)\epsilon_k^2.
    \end{split}
\end{equation}

\textit{\underline{*Step 3.}} If $\phi$ is bounded from below, by following the proof argument of (\ref{equ:sippinequhatx1tvalueuse2}), we can deduce from (\ref{equ:sippinequ1glip3}) that 
\begin{equation}\label{equ:sippinequhatx1tvalueglip3}
    \begin{split}
        &\mb{E}_k\left[\nu_{k+1}^{-2}\cdot\left(\mephi{k+1}\left(x_{k+1}\right) - \phi^{*}\right)\right] - \nu_k^{-2}\cdot\left(\mephi{k}\left(x_k\right) - \phi^{*}\right)\\ 
        \leq\ &\frac{\olrho\left(1+\frac{1}{t}\right)\rho_{\varphi,m}L_\phi^2\glip\left(\Norm{x_k}_2\right)^2}{2\left(1-2\tau \alpha_k\right)}\alpha_k^2 + \frac{\olrho\left(1+t\right)}{2}\epsilon_k^2\\
        &\hspace{0em}+\ \frac{\left(1+\frac{1}{t}\right)\left[1+\eta \alpha_k - \olrho \rho_k^{-2}\alpha_k\right] - \left(1-2\tau \alpha_k\right)}{2\olrho\left(1-2\tau \alpha_k\right)}\cdot \nu_k^{-2}\Norm{\olrho \left(x_k - \olx_k\right)}_{M_{k-1}}^2.\\
    \end{split}
\end{equation}
Let $b \triangleq 4(\olrho\rho_k^{-2} - \eta)(\olrho\rho_k^{-2}-\eta-2\tau)$. Using $\olrho > B_\rho \cdot \left(\eta+2\tau\right)$, one can easily check that $b>0$. Then, the estimate (\ref{equ:mephiinequ3}) follows by applying (\ref{equ:sippinequhatx1tvalueglip3}) with $t = \frac{1}{b\alpha_k^2}$. This completes the proof of assertion (\romannumeral3). 

\section{Proofs for the Results in Section \ref{subsec:sipprateingradme}}\label{secapx:proofofsipprategradme}
In this section, we provide proof details of Claim \ref{claim:mephiinequ1further} and Corollary \ref{coro:sipprate}. 
Recall that function $\varsigma_\beta$ is defined in (\ref{equ:varphi}) as follows:
\begin{align*}
    \varsigma_\beta(t)\triangleq \begin{cases}
        \frac{t^\beta - 1}{\beta}\quad &\mbox{if}\ \beta \neq 0,\\
        \ln(t)\quad &\mbox{if}\ \beta = 0,
    \end{cases}
\end{align*}
for all $t\in \mb{R}_{++}$. We present the following auxiliary lemma summarizing several useful inequalities to be subsequently employed (see \cite[Lemma~C.4]{zhu2025tight} for a detailed proof).
\begin{lema}\label{lema:inequtech}
    The following assertions hold: 
    \begin{enumerate}
        \item Let $\beta\in \mb{R}$. Then the function $\varsigma_\beta$ satisfies $\varsigma_\beta(1) = 0$, and 
        \begin{equation}\label{inequ:varphi}
            \varsigma_\beta(t)\leq \begin{cases}
                \frac{1}{-\beta}\quad &\mbox{if}\ \beta < 0,\\
                \frac{t^\beta}{\beta}\quad &\mbox{if}\ \beta > 0.
            \end{cases}
        \end{equation}
        for all $t\in \mb{R}_{++}$.
        \item The estimate $\frac{k+1}{k} \leq 2$ holds for all $k\in \mb{Z}_{+}$.
        \item Let $\beta\leq 0$ and $j,k\in \mb{N}$ with $k\geq j+1$. Then, we have 
        \begin{equation}\label{equ:sumbetal0}
            \begin{split}
                \varsigma_{\beta + 1}(k+1)-\varsigma_{\beta+1}(j+1)\leq &\sum_{i=j+1}^k i^\beta 
                \leq \frac{1}{2^\beta}\cdot\left[\varsigma_{\beta + 1}(k+1)-\varsigma_{\beta+1}(j+1)\right].
            \end{split}
        \end{equation}
    \end{enumerate}
\end{lema}

\subsection{Proof of Claim \ref{claim:mephiinequ1further}}\label{apx:proofclaimmephiinequ1further}
Fix any $k\in \mb{Z}_{+}$. Condition $\olrho > B_\rho\cdot(\eta + \tau)$ implies $\olrho > \eta + \tau$ since $B_\rho\geq 1$ by (\ref{equ:rhonurelationapx}). Assumption \ref{asm:sppaall} implies $\phi$ is proper, which combined with $\phi^{*}>-\infty$ yields $\phi^{*}\in \mb{R}$. It then follows from Lemma \ref{lema:welldefinedolxk} that $\olx_k$ is well-defined and thus $\mephi{k}(x_k) - \phi^{*}\geq 0$ (see the proof of (\ref{equ:lbmeapx})). Since $\sup_{k\in \mb{Z}_{+}} \alpha_k \leq \alphaou < (\olrho)^{-1}$, Lemma \ref{lema:onestepprogress} implies that $\tlx_{k+1}$ is well-defined. On the event that $\sup_{k\in \mb{Z}_{+}} \Norm{x_k}_2 < \infty$, by following an argument analogous to that in the proof of (\ref{equ:boundtlxkp1}) from Theorem \ref{thm:sippstabilityconvex}, we can deduce from criterion (\ref{equ:criteriaa}) combined with Assumption \ref{asm:mk} that $\sup_{k\in \mb{Z}_{+}} \Norm{\tlx_k}_2 < \infty$. Thus, we denote by $U$ a bounded open convex subset of $V$ such that both $\{x_k\}$ and $\{\tlx_k\}$ are contained in $U$. Use the same notations as in Lemma \ref{lema:inequalitydistance}(\romannumeral1). In this setting, setting $\epsilon_k = \gamma \alpha_k^2$ in (\ref{equ:mephiinequ1}) from Lemma \ref{lema:inequalitydistance}(\romannumeral1) yields 
\begin{equation}\label{equ:mephiinequ1further}
    \begin{split}
        &\mb{E}_k\left[\nu_{k+1}^{-2}\cdot\left(\mephi{k+1}\left(x_{k+1}\right) - \phi^{*}\right)\right] - \nu_k^{-2}\cdot\left(\mephi{k}\left(x_k\right) - \phi^{*}\right)\\ 
        \leq\ &\frac{\olrho\left(1+b\alpha_k^2\right)}{2}\left[\frac{\rho_{f,m,U}L_F\left(U\right)^2}{\left(1-\tau \alpha_k\right)\left(1-\olrho\rho_k^2\alpha_k\right)} + \frac{\gamma^2}{b}\right]\alpha_k^2\\
        &\hspace{0em} - \frac{\left(\olrho\rho_k^{-2}-\eta-\tau\right)\left[1-2\left(\olrho\rho_k^{-2}-\eta\right)\alpha_k\right]^2}{2\olrho\left(1-\tau \alpha_k\right)}\alpha_k \cdot \nu_k^{-2}\Norm{\olrho \left(x_k - \olx_k\right)}_{M_{k-1}}^2,
    \end{split}
\end{equation}
where $b = 4(\olrho\rho_k^{-2} - \eta)(\olrho\rho_k^{-2}-\eta-\tau)$. Condition $\olrho > B_\rho\cdot(\eta+\tau)$ and $\sup_{k\in \mb{Z}_{+}}\alpha_k\leq \alphaou < (\olrho \max\{2,B_\rho\})^{-1}$ combined with (\ref{equ:rhonurelationapx}) confirms that (\ref{equ:medecrease}) holds with postive scalars $\theta_1$ and $\theta_2$ defined in (\ref{def:theta1theta2}).

\subsection{Proof of Corollary \ref{coro:sipprate}}\label{apx:proofcorosipprate}
Recall that $\beta\in (0,1]$. The combination of (\ref{equ:sumbetal0}) with Lemma \ref{lema:inequtech}(\romannumeral1) implies
\begin{equation}\label{equ:ibetasum}
    \begin{split}
        \sum_{i=1}^k i^{-2\beta} &\leq 4^\beta\cdot \left[\varsigma_{1-2\beta}\left(k+1\right) - \varsigma_{1-2\beta}\left(1\right)\right] = 4^\beta \cdot\varsigma_{1-2\beta}\left(k+1\right),\\
        \sum_{i=1}^k i^{-\beta} &\geq \varsigma_{1-\beta}\left(k+1\right) - \varsigma_{1-\beta}\left(1\right) = \varsigma_{1-\beta}\left(k+1\right).
    \end{split}
\end{equation}
Moreover, if $\beta\in (0,1)$, then $1-\beta>0$ and thus
\begin{equation}\label{equ:1beta2}
    \frac{1}{2}\left(k+1\right)^{1-\beta} \geq 1\quad \mbox{for all}\ k\geq K_\beta \triangleq \left\lceil 2^{\frac{1}{1-\beta}}\right\rceil.
\end{equation}
\begin{enumerate}
    \item If $\beta\in (0,\frac{1}{2})$, we have $1-2\beta > 0$ and $1-\beta > 0$, which combined with Lemma \ref{lema:inequtech}(\romannumeral1) implies 
    \begin{align*}
            \varsigma_{1-2\beta}\left(k+1\right) &\leq \frac{\left(k+1\right)^{1-2\beta}}{1-2\beta},\\
            \varsigma_{1-\beta}\left(k+1\right) &= \frac{\left(k+1\right)^{1-\beta}-1}{1-\beta}
            \overset{(\star)}{\geq} \frac{\left(k+1\right)^{1-\beta}-\frac{1}{2}\left(k+1\right)^{1-\beta}}{1-\beta}\\ &= \frac{\left(k+1\right)^{1-\beta}}{2\left(1-\beta\right)}\quad \mbox{for}\ k\ \mbox{sufficiently large},
    \end{align*}
    where $(\star)$ comes from (\ref{equ:1beta2}).
    \item If $\beta = \frac{1}{2}$, we have $1-2\beta = 0$ and $1-\beta > 0$, which combined with Lemma \ref{lema:inequtech}(\romannumeral1) implies
    \begin{align*}
            \varsigma_{1-2\beta}\left(k+1\right) &= \ln\left(k+1\right),\\
            \varsigma_{1-\beta}\left(k+1\right) &\geq \frac{\left(k+1\right)^{1-\beta}}{2\left(1-\beta\right)}\quad \mbox{for}\ k\ \mbox{sufficiently large}.
    \end{align*}
    \item If $\beta\in (\frac{1}{2},1)$, we have $1-2\beta < 0$ and $1-\beta > 0$, which combined with Lemma \ref{lema:inequtech}(\romannumeral1) implies
    \begin{align*}
            \varsigma_{1-2\beta}\left(k+1\right) &\leq \frac{1}{2\beta-1},\\
            \varsigma_{1-\beta}\left(k+1\right) &\geq \frac{\left(k+1\right)^{1-\beta}}{2\left(1-\beta\right)}\quad \mbox{for}\ k\ \mbox{sufficiently large}.
    \end{align*}
    \item If $\beta = 1$, we have $1-2\beta < 0$ and $1-\beta = 0$, which combined with Lemma \ref{lema:inequtech}(\romannumeral1) implies
    \begin{align*}
            \varsigma_{1-2\beta}\left(k+1\right) &\leq \frac{1}{2\beta-1} = 1\quad \mbox{and}\quad 
            \varsigma_{1-\beta}\left(k+1\right) = \ln \left(k+1\right).
    \end{align*}
\end{enumerate}
By Lemma \ref{lema:inequtech}(\romannumeral2), we have $\frac{k+1}{k}\leq 2$ holds for all $k\in \mb{Z}_{+}$. Using (\ref{equ:ibetasum}) and appealing to the results shown above, we can easily deduce the estimate (\ref{equ:asymptoticrate}) from (\ref{asm:epsilonkalphakdiminishingwcvx2}) together with the upper bound (\ref{equ:bdminnormstar}) derived in Theorem \ref{thm:sippmdiminishing2}.

\section{Proofs for the Results in Section \ref{subsec:sipprateinsqdistopt}}\label{secapx:proofofsippratesqdistopt}
In this setion, we provide a detailed proof of Claim \ref{claim:sippinequ1imply}. 
Let $c_1\triangleq \mu_\infty^2\cdot\overline{c}_1$. As noted in Remark \ref{remark:lipqgmknorm}, under Assumption \ref{asm:mk} and \ref{asm:sppa2}, we have
\begin{equation}\label{equ:qgmknorm}
    \phi(x)\geq \phi^{*} + c_1 \mr{dist}_{M_k}(x,\mc{X}^{*})^2\quad \mbox{for all}\ x\in \mb{R}^d
\end{equation}
holds for all $k\in \mb{Z}_{+}$. 

\subsection{Proof of Claim \ref{claim:sippinequ1imply}}\label{subsubsecapx:proofrecursive} 
By assumption, we have $\sup_{k\in \mb{Z}_{+}}\alpha_k \leq \alphaou < (\olrho)^{-1} < \tau^{-1}$, it then follows from Lemma \ref{lema:onestepprogress} that $\tlx_{k+1} \triangleq \prox_{\alpha_k \olvphi_{x_k}(\cdot;S_k^{1:m}}^{M_k}(x_k)$ is well-defined and (\ref{inequ:onestep}) holds. Before proving Claim \ref{claim:sippinequ1imply}, we first demonstrate the following technical identity:
\begin{subequations}
    \begin{align}
         \mb{E}_k\left[\dist_{M_k}\left(\tlx_{k+1},\mc{X}^{*}\right)^2\right] &= \inf_{x^{*}\in \mc{X}^{*}}\mb{E}_k \left[\Norm{\tlx_{k+1}-x^{*}}_{M_k}^2\right]\quad \mbox{and}\label{equ:distexpectation}\\
        \mb{E}_k\left[\dist_{M_k}\left(x_{k+1},\mc{X}^{*}\right)^2\right] &= \inf_{x^{*}\in \mc{X}^{*}}\mb{E}_k \left[\Norm{x_{k+1}-x^{*}}_{M_k}^2\right].\label{equ:distexpectationxk}
    \end{align}
\end{subequations}
With the estimate (\ref{equ:distexpectation}) and (\ref{equ:distexpectationxk}), we can prove the claim. Fix any $k\in \mb{Z}_{+}$. It follows from (\ref{equ:distexpectation}) and (\ref{equ:distexpectationxk}) that 
\begin{equation}\label{equ:cauchyschwarzdistance}
    \begin{split}
        \mb{E}_k\left[\mr{dist}_{M_k}\left(x_{k+1},\mc{X}^{*}\right)^2\right] &= \inf_{x^{*}\in \mc{X}^{*}}\mb{E}_k \left[\Norm{x_{k+1}-x^{*}}_{M_k}^2\right] \\
        &\leq \inf_{x^{*}\in \mc{X}^{*}} \mb{E}_k \left[\left(1 + \frac{1}{t}\right)\left\Vert\tilde{x}_{k+1} - x^{*}\right\Vert_{M_k}^2 + \left(1+t\right)\epsilon_k^2\right]\\
        &= \left(1 + \frac{1}{t}\right)\mb{E}_k\left[\mr{dist}_{M_k}\left(\tilde{x}_{k+1},\mc{X}^{*}\right)^2\right] + \left(1+t\right)\epsilon_k^2
    \end{split}
\end{equation}
holds for any $t\in \mb{R}_{++}$, where the inequality is due to (\ref{equ:cauchyschwarz}) in Lemma \ref{lema:onestepprogress} with $x = x^{*}$ for $x^{*}\in \mc{X}^{*}$. 
Note that both $\{x_k\}$ and $\{\tlx_k\}$ are contained in $U$. Fix any $x^{*}\in \mc{X}^{*}$. Use the same notations as in Lemma \ref{lema:inequalitydistance}(\romannumeral1). Setting $x=x^{*}$ in the estimate (\ref{equ:inequ1max}) of Lemma \ref{lema:inequalitydistance}(\romannumeral1) yields
\begin{align*}
    &\frac{1-\tau \alpha_k}{2\alpha_k}\mb{E}_k\left[\Vert \tilde{x}_{k+1} - x^{*}\Vert_{M_k}^2\right] - \frac{1+\eta \alpha_k}{2\alpha_k}\mb{E}_k\left[\Vert x_k - x^{*}\Vert_{M_k}^2\right]\\ \leq\ & - \mb{E}_k\left[\phi(\tilde{x}_{k+1}) - \phi^{*}\right] + \frac{\rho_{f,m,U}L_F\left(U\right)^2}{2}\alpha_k\\
    \leq\ & - c_1\mb{E}_k\left[\dist_{M_k}\left(\tlx_{k+1},\mc{X}^{*}\right)^2\right] + \frac{\rho_{f,m,U}L_F\left(U\right)^2}{2}\alpha_k,
\end{align*}
where the last inequality comes from (\ref{equ:qgmknorm}) by Assumption \ref{asm:sppa2}. Combining the preceding inequality with (\ref{equ:distexpectation}) yields 
\begin{align*}
    &\frac{1-\tau \alpha_k}{2\alpha_k}\mb{E}_k\left[\mr{dist}_{M_k}\left(\tilde{x}_{k+1},\mc{X}^{*}\right)^2\right] = \frac{1-\tau \alpha_k}{2\alpha_k}\inf_{x^{*}\in \mc{X}^{*}} \mb{E}_k\left[\Norm{\tilde{x}_{k+1} - x^{*}}_{M_k}^2\right]\\
        \leq\ &\inf_{x^{*}\in \mc{X}^{*}}\lrbrackets{\frac{1+\eta \alpha_k}{2\alpha_k}\mb{E}_k\left[\Vert x_k - x^{*}\Vert_{M_k}^2\right] - c_1\mb{E}_k\left[\dist_{M_k}\left(\tlx_{k+1},\mc{X}^{*}\right)^2\right] + \frac{\rho_{f,m,U}L_F\left(U\right)^2}{2}\alpha_k} \\
        \leq\ 
        &\inf_{x^{*}\in \mc{X}^{*}}\lrbrackets{\frac{1+\eta \alpha_k)}{2\alpha_k}\rho_k^2\Vert x_k - x^{*}\Vert_{M_{k-1}}^2 - c_1\mb{E}_k\left[\dist_{M_k}\left(\tlx_{k+1},\mc{X}^{*}\right)^2\right] + \frac{\rho_{f,m,U}L_F\left(U\right)^2}{2}\alpha_k} \\ =\ 
        &\frac{\rho_k^2\left(1+\eta \alpha_k\right)}{2\alpha_k}\dist_{M_{k-1}}\left(x_k,\mc{X}^{*}\right)^2 - c_1\mb{E}_k\left[\dist_{M_k}\left(\tlx_{k+1},\mc{X}^{*}\right)^2\right] + \frac{\rho_{f,m,U}L_F\left(U\right)^2}{2}\alpha_k,
\end{align*}
where the last inequality follows from Assumption \ref{asm:mk}. 
Given that $c_1 > \frac{\tau+\eta}{2}$, rearranging the preceding inequality gives
\begin{align*}
        &\mb{E}_k\left[\mr{dist}_{M_k}\left(\tilde{x}_{k+1},\mc{X}^{*}\right)^2\right]\\ \leq\ &\frac{\rho_k^2\left(1+\eta\alpha_k\right)}{1+2c_1\alpha_k-\tau\alpha_k} \mr{dist}_{M_{k-1}}\left(x_k,\mc{X}^{*}\right)^2 +\rho_{f,m,U}L_F\left(U\right)^2 \cdot\frac{\alpha_k^2}{1+2c_1\alpha_k-\tau\alpha_k}.
\end{align*}
Since $\nu_{k+1} = \rho_k \nu_k$ and $\nu_{k+1} \geq 1$, we can deduce from the estimate above that 
\begin{equation}\label{equ:proximalpointdistanceoptimalset1}
    \begin{split}
        \mb{E}_k\left[\nu_{k+1}^{-2}\mr{dist}_{M_k}\left(\tilde{x}_{k+1},\mc{X}^{*}\right)^2\right]
        \leq\ &\frac{1+\eta\alpha_k}{1+2c_1\alpha_k-\tau\alpha_k} \nu_k^{-2}\mr{dist}_{M_{k-1}}\left(x_k,\mc{X}^{*}\right)^2\\ &\hspace{0em} + \rho_{f,m,U}L_F\left(U\right)^2 \cdot\frac{\alpha_k^2}{1+2c_1\alpha_k-\tau\alpha_k}.\\
    \end{split}
\end{equation}
Recall that $c_{\tau,\eta} \triangleq c_1 - \frac{\tau+\eta}{2}$, $\tilde{c}_1 \triangleq \frac{\olrho c_{\tau,\eta}}{\olrho + \eta}$ and $s\in (0,2)$. Then, condition $c_1 > \frac{\tau+\eta}{2}$ implies $\tilde{c}_1,c_{\tau,\eta}>0$. Applying (\ref{equ:proximalpointdistanceoptimalset1}) to (\ref{equ:cauchyschwarzdistance}) with $t = \frac{1+\eta\alpha_k + sc_{\tau,\eta}\alpha_k}{(2-s)c_{\tau,\eta}\alpha_k}$ and $\epsilon_k = \gamma \alpha_k^{\frac{3}{2}}$, and using the condition $\alpha_k\leq \alphaou < (\olrho)^{-1}$ and $\nu_{k+1}\geq 1$, we have 
\begin{align*}
        &\mb{E}_k\left[\nu_{k+1}^{-2}\mr{dist}_{M_k}\left(x_{k+1},\mc{X}^{*}\right)^2\right]
        \leq \frac{1}{1+s\tilde{c}_1\alpha_k}\cdot \nu_k^{-2}\mr{dist}_{M_{k-1}}\left(x_k,\mc{X}^{*}\right)^2 + \left[\rho_{f,m,U} L_F\left(U\right)^2\right.\\
        &\hspace{5em}\left.+\ \frac{\left(1 + 2c_1\alpha_k-\tau\alpha_k\right)\left(2\left(1+\eta\alpha_k\right) + s\left(2c_1-\tau-\eta\right)\alpha_k\right)}{(2-s)\left(2c_1-\tau-\eta\right)\left(1+\eta\alpha_k\right)}\gamma^2\right]\cdot \frac{\alpha_k^2}{1+s\tilde{c}_1\alpha_k}.
\end{align*}
Note that $\delta_k \triangleq \mb{E}[\nu_k^{-2}\dist_{M_{k-1}}(x_k,\mc{X}^{*})^2]$ for each $k\in \mb{Z}_{+}$. Taking the expectation of both sides of the preceding estimate, applying the law of total expectation and using the definition of $\tilde{C}_{f,m,U,\tau,\eta,c_1,\gamma}(\cdot)$ defined in (\ref{equ:sipprateboundconstant1}), we can prove that Claim \ref{claim:sippinequ1imply} remains valid. The only remaining task is to prove (\ref{equ:distexpectation}) and (\ref{equ:distexpectationxk}).
\begin{proof}[Proof of estimates (\ref{equ:distexpectation}) and (\ref{equ:distexpectationxk})]
    Using the property that 
    \begin{align*}
            \dist_{M_k}\left(\tlx_{k+1},\mc{X}^{*}\right)^2 = \inf_{\olx\in \mc{X}^{*}} \Norm{\tlx_{k+1}-\olx}_{M_k}^2 \leq \Norm{\tlx_{k+1}-x^{*}}_{M_k}^2\quad \mbox{for all}\ x^{*}\in \mc{X}^{*},
    \end{align*}
    we have 
    \begin{equation}\label{equ:distexpectationleq}
        \begin{split}
            &\mb{E}_k\left[\dist_{M_k}\left(\tlx_{k+1},\mc{X}^{*}\right)^2\right] \leq \mb{E}_k\left[\Norm{\tlx_{k+1}-x^{*}}_{M_k}^2\right]\quad \mbox{for all}\ x^{*}\in \mc{X}^{*}\\
            \Rightarrow\quad &\mb{E}_k\left[\dist_{M_k}\left(\tlx_{k+1},\mc{X}^{*}\right)^2\right] \leq \inf_{x^{*}\in \mc{X}^{*}}\mb{E}_k\left[\Norm{\tlx_{k+1}-x^{*}}_{M_k}^2\right].
        \end{split}
    \end{equation}
    On the other hand, for any fixed scalar $\epsilon>0$, the definition of the infimum implies the existence of  $\olx_\epsilon\in \mc{X}^{*}$ such that 
    \begin{align*}
        \dist_{M_k}\left(\tlx_{k+1},\mc{X}^{*}\right)^2 > \Norm{\tlx_{k+1}-\olx_\epsilon}_{M_k}^2 - \epsilon,
    \end{align*}
    which yields
    \begin{align*}
            \mb{E}_k\left[\dist_{M_k}\left(\tlx_{k+1},\mc{X}^{*}\right)^2\right] &\geq \mb{E}_k\left[\Norm{\tlx_{k+1}-\olx_\epsilon}_{M_k}^2\right] - \epsilon \geq \inf_{x^{*}\in \mc{X}^{*}}\mb{E}_k\left[\Norm{\tlx_{k+1}-x^{*}}_{M_k}^2\right] - \epsilon.
    \end{align*}
    Since $\epsilon>0$ is arbitrary, we have 
    \begin{equation}\label{equ:distexpectationgeq}
        \begin{split}
            \mb{E}_k\left[\dist_{M_k}\left(\tlx_{k+1},\mc{X}^{*}\right)^2\right] &\geq \inf_{x^{*}\in \mc{X}^{*}}\mb{E}_k\left[\Norm{\tlx_{k+1}-x^{*}}_{M_k}^2\right].
        \end{split}
    \end{equation}
    Thus, the estimate (\ref{equ:distexpectation}) follows immediately from (\ref{equ:distexpectationleq}) combined with (\ref{equ:distexpectationgeq}). Similarly, we can also show that the estimate (\ref{equ:distexpectationxk}) holds.
\end{proof}

\section{Details of Numerical Experiments}\label{apx:numexp}
In this section, we present several auxiliary results and structural conditions that underpin the practical implementation of the model-based ispPPA for solving the regularized regression problems introduced in Section \ref{sec:sippnumericalexperiments}. Both the $\ell_1$ and MCP regularized regression models can be formulated as the following finite-sum composite optimization problem: 
\begin{equation}\label{equ:finitesum}
    \min_{x\in \mb{R}^d}\ \phi(x) \triangleq F\left(Ax\right) + r(x), 
\end{equation}
where the regularizer $r\colon \mb{R}^d\to \pminf$ is proper, closed and $\oltheta$-weakly convex, the matrix 
\begin{align*}
        A = \begin{pmatrix}
            A_1^\top & \cdots & A_n^\top
        \end{pmatrix}^\top\in \mb{R}^{N\times d}\quad \mbox{with}\ N = \sum_{i=1}^n N_i\ \mbox{and}\ A_i\in \mb{R}^{N_i\times d}\ \mbox{for}\ i\in [n]
\end{align*}
is given, and the function $F\colon \mb{R}^N \to \pminf$ is defined by 
\begin{align*}
    F\left(y\right)\triangleq \frac{1}{n}\sum_{i=1}^n f_i\left(y_i\right)\quad \mbox{for all}\ y = \begin{pmatrix}
        y_1^\top & \cdots & y_n^\top
    \end{pmatrix}^\top\in \mb{R}^N,
\end{align*}
with component functions $f_i\colon \mb{R}^{N_i}\to \pminf$ that are proper, closed and convex for all $i\in [n]$. The remainder of this section provides (\romannumeral1) a principled way to select preconditioners ensuring that Assumption \ref{asm:mk} holds, and (\romannumeral2) based on this choice, sufficient conditions under which criterion (\ref{equ:criteriaa}) is satisfied when solving subproblems inexactly.

\subsection{Proper Choices of Preconditioner}\label{apx:subsecpreconditioners}
We first present a constructive setup for selecting preconditioners that satisfies Assumption \ref{asm:mk} in the context of the finite-sum problem (\ref{equ:finitesum}).
\begin{claim}\label{claim:choicemk}
    Consider the finite-sum optimization problem (\ref{equ:finitesum}). Let the iterates $\{x_k\}$ be generated by (\ref{equ:sippm}) with stepsizes $\alpha_k\in \mb{R}_{++}$, parameters $\epsilon_k\in \mb{R}_{+}$, and positive-definite linear mappings $M_k\colon \mb{R}^d\to \mb{R}^d$. Then, Assumption \ref{asm:mk} holds under the following construction: set $M_0 = I_d$, and for all $k\in \mb{Z}_{+}$,
    \begin{enumerate}
        \item $M_k \triangleq I_d + \alpha_k \tau_k \sum_{i\in S_k^{1:m}}A_i^\top A_i$ with $A_i\in \mb{R}^{N_i\times d}$ for $i\in S_k^{1:m}$,
        \item $\alpha_k \triangleq \alpha_0 k^{-\beta}$ for some $\alpha_0 \in \mb{R}_{++}$ and $\beta\in [0,1]$, and
        \item $\tau_k \triangleq \tau_0 k^{\eta}$ for some $\tau_0 \in \mb{R}_{++}$ and $\eta < \beta - 1$.
    \end{enumerate} 
\end{claim}
\begin{remark}\label{remark:mksubmatrixform}
    Consider the setting of Claim \ref{claim:choicemk}. Recall that $S_k^{1:m} \triangleq \{S_k^i\}_{i=1}^m$. Define the subsampled data matrix $A_{S_k^{1:m},:}\triangleq \begin{pmatrix}
        A_{S_k^1}^\top & \cdots & A_{S_k^m}^\top
        \end{pmatrix}^\top\in \mb{R}^{\tlm\times d}$ and let $\tlF\colon \mb{R}^{\tlm} \to \pminf$ be the function defined by 
    \begin{align*}
        \tlF(\tly) \triangleq \frac{1}{m}\sum_{i=1}^mf_{S_k^i}(y_{S_k^i})\quad \mbox{for all}\ \tly = \begin{pmatrix}
            y_{S_k^1}^\top & \cdots & y_{S_k^m}^\top
        \end{pmatrix}^\top\in \mb{R}^{\tlm},
    \end{align*}
    where $\tlm \triangleq \sum_{i=1}^m N_i$. 
    Using the definition of $A_{S_k^{1:m},:}$, we have $M_k = I_d + \alpha_k \tau_k A_{S_k^{1:m},:}^\top A_{S_k^{1:m},:}$, and for any $x\in \mb{R}^d$,
    \begin{equation}\label{equ:normmk}
          \frac{1}{m}\sum_{i\in S_k^{1:m}} f_i\left(A_ix\right) = \tlF\left(A_{S_k^{1:m},:}x\right)\quad \mbox{and}\quad \Norm{x}_{M_k}^2 = \Norm{x}_2^2 + \alpha_k \tau_k \Norm{A_{S_k^{1:m},:}x}_2^2.
    \end{equation}
    Consequently, the subproblem of the preconditioned isPPA (i.e., Algorithm \ref{algo:ispppa} with model functions $f_x(\cdot;s) =f(\cdot;s)$) at iteration $k$ takes the form of
    \begin{align*}
        x_{k+1} \overset{\epsilon_k}{\approx} \arg\min_{x\in \mb{R}^d}\lrbrackets{\tlF(A_{S_k^{1:m},:}x) + r(x) + \frac{1}{2\alpha_k} \Norm{x-x_k}_2^2 + \frac{\tau_k}{2} \Norm{A_{S_k^{1:m},:}\left(x-x_k\right)}_2^2}.
    \end{align*}
\end{remark}
\begin{proof}[Proof of Claim \ref{claim:choicemk}]
Set $\mu_\infty \triangleq (1+\alpha_0\tau_0\Norm{A}_2^2)^{-\frac{1}{2}}$ and $L_\infty \triangleq 1$. Use the same notations as in Remark \ref{remark:mksubmatrixform}. Fix any $k\in \mb{Z}_{+}$. The definition of $A_{S_k^{1:m},:}$ implies that $A_{S_k^{1:m},:}$ is the submatrix generated by the blocks of $A$ indexed by $S_k^{1:m}$. It then follows from (\ref{equ:normmk}) that
\begin{equation}\label{inequ:norm2mk}
    \begin{split}
        \Norm{x}_2^2 \leq \Norm{x}_{M_k}^2 &= \Norm{x}_2^2 + \alpha_k \tau_k \Norm{A_{S_k^{1:m},:} x}_2^2 \leq \Norm{x}_2^2 + \alpha_k \tau_k \Norm{A x}_2^2\\ &\leq \left(1+\alpha_k\tau_k\Norm{A}_2^2\right) \Norm{x}_2^2\quad \mbox{for all}\ x\in \mb{R}^d.
    \end{split}
\end{equation}
Then, by defining $\mu_k \triangleq (1+\alpha_k\tau_k\Norm{A}_2^2)^{-\frac{1}{2}}$ and $L_k \triangleq 1$, we deduce from $\alpha_k \leq \alpha_0$ and $\tau_k \leq \tau_0$ that $\mu_k \Norm{x}_{M_k} \leq \Norm{x}_2 \leq L_k \Norm{x}_{M_k}$ with $\mu_k \geq \mu_\infty$ and $L_k \leq L_\infty$. 

Note that $\alpha_k \tau_k = \alpha_0 \tau_0 k^{\eta - \beta}$ for each $k\in \mb{Z}_{+}$. Since $\eta - \beta < -1$, then $\sum_{k=0}^\infty \alpha_k\tau_k < \infty$ and
\begin{equation}\label{equ:nonincreasing}
    \mbox{the sequnce}\ \lrbrackets{\alpha_k\tau_k}_{k\geq 0}\ \mbox{is nonincreasing}.
\end{equation}
Fix any $k\in \mb{Z}_{+}$ and $0\neq x\in \mb{R}^d$.
\begin{enumerate}
    \item If $k=1$, one can conclude that 
    \begin{align*}
            \frac{\Norm{x}_{M_k}^2}{\Norm{x}_{M_{k-1}}^2} &\leq \frac{\left(1 + \alpha_k \tau_k \Norm{A}_2^2\right)\Norm{x}_2^2}{\Norm{x}_2^2} = 1 + \alpha_k \tau_k \Norm{A}_2^2 \leq 1 + \alpha_{k-1}\tau_{k-1}\Norm{A}_2^2,
    \end{align*}
    where the first inequality comes from (\ref{inequ:norm2mk}) and the fact that $M_0 \triangleq I_d$, and the last is due to (\ref{equ:nonincreasing}). Similarly, 
    \begin{align*}
            \frac{\Norm{x}_{M_k}^2}{\Norm{x}_{M_{k-1}}^2} &\geq \frac{\Norm{x}_2^2}{\Norm{x}_2^2} = 1
            \geq \left(1 + \alpha_{k-1}\tau_{k-1}\Norm{A}_2^2\right)^{-1}.
    \end{align*}
    \item Now consider the case where $k\geq 2$. Using (\ref{inequ:norm2mk}) and (\ref{equ:nonincreasing}), we have 
    \begin{align*}
            \frac{\Norm{x}_{M_k}^2}{\Norm{x}_{M_{k-1}}^2} &\leq \frac{\left(1 + \alpha_k \tau_k \Norm{A}_2^2\right)\Norm{x}_2^2}{\Norm{x}_2^2} = 1 + \alpha_k \tau_k \Norm{A}_2^2 \leq 1 + \alpha_{k-1}\tau_{k-1}\Norm{A}_2^2,
    \end{align*}
    and 
    \begin{align*}
            \frac{\Norm{x}_{M_k}^2}{\Norm{x}_{M_{k-1}}^2} &\geq \frac{\Norm{x}_2^2}{\left(1+\alpha_{k-1} \tau_{k-1} \Norm{A}_2^2\right)\Norm{x}_2^2} = \left(1 + \alpha_{k-1}\tau_{k-1}\Norm{A}_2^2\right)^{-1}.
    \end{align*}
\end{enumerate}
Set $\rho_0 \triangleq 1$ and $\rho_k \triangleq (1+\alpha_{k-1}\tau_{k-1}\Norm{A}_2^2)^{\frac{1}{2}}$ for each $k\in \mb{Z}_{+}$. It then follows that the sequence $\{\rho_k\}$ satisfies $\rho_k\in [1,\infty)$ and
\begin{align*}
    \rho_k^{-2} \leq \mb{E}_k \left[\frac{\Norm{x}_{M_k}^2}{\Norm{x}_{M_{k-1}^2}}\right] \leq \rho_k^2\quad \mbox{for all}\ 0\neq x\in \mb{R}^d.
\end{align*}
Moreover, the sequence $\{\nu_k\}$ is given by 
\begin{align*}
    \nu_k = \prod_{i=0}^{k-1} \rho_i = \begin{cases}
        1\quad &\mbox{if}\ k = 1,\\
        \prod_{i=1}^{k-1} \left(1+\alpha_{i-1}\tau_{i-1}\Norm{A}_2^2\right)^{\frac{1}{2}} = \prod_{i=0}^{k-2} \left(1+\alpha_{i}\tau_{i}\Norm{A}_2^2\right)^{\frac{1}{2}}\quad &\mbox{if}\ k\geq 2.
    \end{cases}
\end{align*}
Define $\Sigma \triangleq \exp\left(\frac{1}{2}\Norm{A}_2^2\sum_{k=0}^{\infty}\alpha_{k}\tau_{k}\right)$. Since $\sum_{k=0}^\infty \alpha_k\tau_k < \infty$, we have $\Sigma < \infty$, which combined with (\ref{inequ:exp}) implies
\begin{align*}
        \nu_k &\leq \prod_{i=0}^{k-2} \exp\left(\frac{1}{2}\alpha_{i}\tau_{i}\Norm{A}_2^2\right) = \exp\left(\frac{1}{2}\Norm{A}_2^2\sum_{i=0}^{k-2}\alpha_{i}\tau_{i}\right) \leq \Sigma
\end{align*}
for any $k\geq 2$. Note that $\nu_1 \leq \Sigma$ holds trivially. Then, the sequence $\{\nu_k\}$ is nonincreasing and bounded from above, and thus the limit $\nu_\infty \triangleq \prod_{k=0}^\infty (1+\alpha_{k}\tau_{k}\Norm{A}_2^2)^{\frac{1}{2}} < \infty$ exists. This shows that Assumption \ref{asm:mk} holds and thus completes the proof of Claim \ref{claim:choicemk}.
\end{proof}

\subsection{Solve the Inner-Loop Subproblem}\label{apx:subsecinnerproblem}
Having specified the preconditioner structure, we now establish conditions under which the approximate solution of each subproblem satisfies the stopping criterion (\ref{equ:criteriaa}). 
Recall that $e_{\alpha f}$ denotes the Moreau envelope of a proper, closed and $\oltheta$-weakly convex function $f\colon \mb{R}^d\to \pminf$ with parameter $\alpha\in (0,(\oltheta)^{-1})$.
\begin{claim}\label{claim:ippaequialmnew}
    Consider the following composite optimization problem 
    \begin{align*}
            \min_{x\in \mb{R}^d}\ \Phi(x) \triangleq H(Ax) + p(x),
    \end{align*}
    where $A\in \mb{R}^{\tlm\times d}$ is given, $H\colon \mb{R}^{\tlm}\to \pminf$ is proper, closed and convex, and $p\colon \mb{R}^d\to \pminf$ is proper, closed and $\oltheta$-weakly convex. Fix any $\olx\in \mb{R}^d$, $\alpha \in (0,(\oltheta)^{-1})$, $\tau \in \mb{R}_{++}$ and $\epsilon\in \mb{R}_{+}$. Define 
    \begin{equation}\label{equ:psidualpredisppa}
        \begin{split}
            \Psi(\xi) &\triangleq \frac{\tau}{2}\Norm{\tau^{-1}\xi + A\olx}_2^2 - e_{\tau^{-1}H}\left(\tau^{-1}\xi+A\olx\right)\\
            &\hspace{1.5em}+ \frac{1}{2\alpha}\Norm{\olx - \alpha A^\top \xi}_2^2 - e_{\alpha p}\left(\olx - \alpha A^\top \xi\right)\quad \mbox{for all}\ \xi\in \mb{R}^{\tlm}.
        \end{split}
    \end{equation}
    Let $M\triangleq I_d + \alpha \tau A^\top A$. If there exists some $\tlxi\in \mb{R}^{\tlm}$ such that 
    \begin{equation}\label{equ:ippaequialmnew}
        H\left(A\tlx\right) + p\left(\tlx\right) + \frac{1}{2\alpha}\Norm{\tlx - \olx}_M^2 + \Psi\left(\tlxi\right) - \frac{1}{2\alpha}\Norm{\olx}_M^2 \leq \frac{1-\oltheta\alpha}{2\alpha}\epsilon^2
    \end{equation}
    for some $\tlx\in \mb{R}^d$, then $\tlx$ satisfies
    \begin{equation}\label{equ:scasufficient}
        \tilde{x} \approx \arg\min_{x\in \mb{R}^d}\lrbrackets{H(Ax) + p(x) + \frac{1}{2\alpha}\Norm{x - \olx}_M^2}\quad \mbox{where}\quad \Norm{\tilde{x} - \prox_{\alpha \Phi}^M(\olx)}_M \leq \epsilon.
    \end{equation}
\end{claim}
\begin{remark}
    To solve the finite-sum optimization problem (\ref{equ:finitesum}) via the preconditioned isPPA, we use the preconditioners defined in Claim \ref{claim:choicemk}. By invoking Claim \ref{claim:ippaequialmnew} with the subsampled matrix $A = A_{S_k^{1:m},:}$, fucntions $H(\cdot)$ and $p(\cdot)$ specified as
    \begin{align*}
            H\left(\tly\right) &\triangleq \tlF\left(\tly\right)\ \mbox{for all}\ \tly\in \mb{R}^{\tlm}\quad \mbox{and}\quad 
            p(x) \triangleq r(x)\ \mbox{for all}\ x\in \mb{R}^d,
    \end{align*}
    stepsizes $\alpha = \alpha_k$ and $\tau = \tau_k$, and preconditioners $M = M_k$, we can derive a sufficient condition for $x_{k+1}$ to satisfy the stopping criterion (\ref{equ:criteriaa}). Specifically, we need to find an approximate minimizer $\xi_{k+1}$ of $\Psi(\cdot)$ defined in (\ref{equ:psidualpredisppa}) such that the termination condition  
    \begin{align*}
        &H\lrrbrackets{A_{S_k^{1:m},:} x_{k+1}} + p\lrrbrackets{x_{k+1}} + \frac{1}{2\alpha_k}\Norm{x_{k+1} - x_k}_{M_k}^2 + \Psi\lrrbrackets{\xi_{k+1}} - \frac{1}{2\alpha_k}\Norm{x_k}_{M_k}^2 \leq \hat{\epsilon}_k
    \end{align*}
    is satisfied for $x_{k+1} = \prox_{\alpha_k p}(x_k - \alpha_k A_{S_k^{1:m},:}^\top \xi_{k+1})$, where $\hat{\epsilon}_k\triangleq \frac{1-\oltheta\alpha_k}{2\alpha_k}\epsilon_k^2$.
\end{remark}
\begin{remark}\label{remark:predid}
    In Claim \ref{claim:ippaequialmnew}, the analysis considers the case where the preconditioner $M \triangleq I_d + \alpha \tau A^\top A$ for $\alpha,\tau > 0$. If the trivial preconditioner $M \triangleq I_d$ is used, we may define
    \begin{align*}
            \Psi(\xi) \triangleq H^\star\left(\xi\right) + \frac{1}{2\alpha}\Norm{\olx - \alpha A^\top \xi}_2^2 - e_{\alpha p}\left(\olx - \alpha A^\top \xi\right)\quad \mbox{for all}\ \xi\in \mb{R}^{\tlm}.
    \end{align*}
    Following the proof arguments of Claim \ref{claim:ippaequialmnew}, we can show that if there exist $\tlxi\in \mb{R}^{\tlm}$ and $\tlx\in \mb{R}^d$ such that condition (\ref{equ:ippaequialmnew}) holds, then $\tlx$ also satisfies (\ref{equ:scasufficient}). 
\end{remark}
\begin{proof}[Proof of Claim \ref{claim:ippaequialmnew}]
Since $M \triangleq I_d + \alpha \tau A^\top A$, we have $\Vert x\Vert_2 \leq \Vert x\Vert_M$ for all $x\in \mb{R}^d$. By assumption, $H$ is proper, closed and convex, and $p$ is proper, closed and $\oltheta$-weakly convex. Then, using condition $\alpha\in (0,(\oltheta)^{-1})$ and following the proof argument of Lemma \ref{lema:onestepprogress}, we can conclude that $\Phi(\cdot)+\frac{1}{2\alpha}\Vert \cdot-\olx\Vert_M^2$ is proper, closed and $\frac{1-\oltheta\alpha}{\alpha}$-strongly convex, which satisfies
\begin{equation}\label{equ:oneprogresssubproblem}
        \Phi(x) + \frac{1}{2\alpha}\Norm{x - \olx}_M^2 - \left(\Phi\left(\hat{x}\right) + \frac{1}{2\alpha}\Norm{\hat{x} - \olx}_M^2\right) \geq \frac{1-\oltheta \alpha}{2\alpha}\Norm{x-\hat{x}}_M^2
\end{equation}
for all $x\in \mb{R}^d$. 
The minimization problem for obtaining $\hat{x}$ can be written as
\begin{equation}\label{equ:primalproblem}
    \begin{split}
        \min_{x\in \mb{R}^d}\ &\Phi\lrrbrackets{x} + \frac{1}{2\alpha}\Norm{x-\olx}_M^2\\
        \Leftrightarrow\quad \min_{(x,y)\in \mb{R}^d\times \mb{R}^{\tlm}}\ &H\lrrbrackets{y} + p\lrrbrackets{x} + \frac{1}{2\alpha}\Norm{x-\olx}_2^2 + \frac{\tau}{2}\Norm{y-A\olx}_2^2\\
        \mbox{s.t.}\quad &Ax-y = 0.
    \end{split}
\end{equation}
Using the definition of $\Psi(\cdot)$ in (\ref{equ:psidualpredisppa}) and following the proof argument of \cite[(6)]{lin2024highly}, one can conclude that the dual problem of (\ref{equ:primalproblem}) takes the form of 
\begin{equation*}
    \begin{split}
        \max_{\xi\in \mb{R}^{\tlm}}\ -\Psi\lrrbrackets{\xi} + \frac{1}{2\alpha}\Norm{\olx}_M^2.
    \end{split}
\end{equation*}
By applying (\ref{equ:oneprogresssubproblem}) with $x = \prox^M_{\alpha\Phi}(\olx)$ and following an argument nearly identical to that in the proof of \cite[Theorem~2.1]{lin2024highly}, we can show that Claim \ref{claim:ippaequialmnew} remains valid.

\end{proof}

\subsection{Equivalence of the KKT Residual}\label{apx:kktresidualcompute}
The following result is crucial for discussing the convergence rates measured by the KKT residual in numerical experiments.
\begin{claim}\label{claim:kktresidualcompute}
    Let $F\colon \mb{R}^{d}\to \mb{R}$ be convex and $\lipsmf$-Lipschitz smooth, and let $r\colon \mb{R}^d\to \pminf$ be proper, closed and $\oltheta$-weakly convex. Define $\phi(\cdot)\triangleq F(\cdot)+r(\cdot)$. Fix any $\alpha\in (0,(\oltheta)^{-1})$. Then, we have 
    \begin{align*}
        \left(1-\frac{\alpha \lipsmf}{1-\alpha \oltheta}\right)\Norm{x-\prox_{\alpha r}\left(x - \alpha \nabla F\left(x\right)\right)}_2 &\leq \Norm{x-\prox_{\alpha \phi}\left(x\right)}_2\\ &\leq \left(1+\frac{\alpha \lipsmf}{1-\alpha \oltheta}\right)\Norm{x-\prox_{\alpha r}\left(x - \alpha \nabla F\left(x\right)\right)}_2
    \end{align*}
    holds for all $x\in \mb{R}^d$.
\end{claim}
\begin{proof}[Proof of Claim \ref{claim:kktresidualcompute}]
Fix any $x \in \mb{R}^d$ and $\alpha\in (0,(\oltheta)^{-1})$. Note that $\phi$ and $r$ are $\oltheta$-weakly convex by assumption. According to \cite[Proposition~3.3]{hoheisel2010proximal}, the proximal mappings $\prox_{\alpha r}(\cdot)$ and $\prox_{\alpha \phi}(\cdot)$ are single-valued and $(1-\alpha \oltheta)^{-1}$-Lipschitz continuous. Denote by $z \triangleq \frac{1}{\alpha}(x-\prox_{\alpha r}(x-\alpha \nabla F(x)))$ and $w \triangleq \nabla F(x-\alpha z) -\nabla F(x)$. Then, it holds that
\begin{equation}\label{equ:proxkkt}
    \begin{split}
        &x-\alpha z = \prox_{\alpha r}\left(x-\alpha \nabla F\left(x\right)\right)\\ \overset{(\star)}{\Leftrightarrow}\quad &\frac{1}{\alpha}\left(x-\alpha \nabla F\left(x\right) - \left(x-\alpha z\right)\right) \in \partial r\left(x-\alpha z\right)\\
        \overset{(\triangle)}{\Leftrightarrow}\quad &\frac{1}{\alpha}\left(x + \alpha\left(\nabla F\left(x-\alpha z\right) - \alpha \nabla F\left(x\right)\right)\right) - \left(x-\alpha z\right)\in \partial \phi\left(x-\alpha z\right)\\
        \overset{(\square)}{\Leftrightarrow}\quad &x-\alpha z = \prox_{\alpha \phi}\left(x+\alpha w\right)
        \quad \Leftrightarrow\quad z = \frac{1}{\alpha}\left(x-\prox_{\alpha \phi}\left(x+\alpha w\right)\right),
    \end{split}
\end{equation}
where $(\star)$ and $(\square)$ uses (\ref{equ:regularsubdiffsubproblem}), and $(\triangle)$ follows from \cite[Exercise~8.8(c)]{rockafellar1998convex} combined with the fact that $r(x-\alpha z)$ is finite and $F(\cdot)$ is continuously differentiable on $\mb{R}^d$. The Lipschitz smoothness of $F$ implies that $\Vert w\Vert_2 \leq \alpha \lipsmf\Vert z\Vert_2$. It then follows from (\ref{equ:proxkkt}) that
\begin{align*}
        \norm{\Norm{z}_2 - \Norm{\frac{1}{\alpha}\left(x-\prox_{\alpha \phi}\left(x\right)\right)}_2} &\leq \Norm{z - \frac{1}{\alpha}\left(x-\prox_{\alpha \phi}\left(x\right)\right)}_2\\
        &= \Norm{\frac{1}{\alpha}\left(x-\prox_{\alpha \phi}\left(x+\alpha w\right)\right) - \frac{1}{\alpha}\left(x-\prox_{\alpha \phi}\left(x\right)\right)}_2\\
        &\leq \frac{1}{\alpha}\cdot \frac{1}{1-\alpha \oltheta}\Norm{\alpha w}_2 \leq \frac{\alpha \lipsmf}{1-\alpha \oltheta}\Norm{z}_2,
\end{align*}
which combined with the definition of $z$ implies 
\begin{align*}
        &\left(1-\frac{\alpha \lipsmf}{1-\alpha \oltheta}\right)\Norm{x-\prox_{\alpha r}\left(x - \alpha \nabla F\left(x\right)\right)}_2
        \leq \Norm{x-\prox_{\alpha \phi}\left(x\right)}_2\\ \leq\ &\left(1+\frac{\alpha \lipsmf}{1-\alpha \oltheta}\right)\Norm{x-\prox_{\alpha r}\left(x - \alpha \nabla F\left(x\right)\right)}_2.
\end{align*}
This completes the proof of Claim \ref{claim:kktresidualcompute}.
\end{proof}

\end{document}